\documentclass[10pt]{amsart}
\usepackage[unicode=true,linktoc=pages]{hyperref}

\usepackage{amssymb,amsthm,amscd,array,xypic,stmaryrd}
\usepackage{enumitem}
\usepackage[all]{xy}
\usepackage{lmodern}
\usepackage[utf8]{inputenc}

\usepackage[french,english]{babel}

\makeatletter
\newenvironment{abstracts}{%
	\ifx\maketitle\relax
	\ClassWarning{\@classname}{Abstract should precede
		\protect\maketitle\space in AMS document classes; reported}%
	\fi
	\global\setbox\abstractbox=\vtop \bgroup
	\normalfont\Small
	\list{}{\labelwidth\z@
		\leftmargin3pc \rightmargin\leftmargin
		\listparindent\normalparindent \itemindent\z@
		\parsep\z@ \@plus\p@
		
		\itemsep\medskipamount
	}%
}{%
	\endlist\egroup
	\ifx\@setabstract\relax \@setabstracta \fi
}

\newcommand{\abstractin}[1]{%
	\otherlanguage{#1}%
	\item[\hskip\labelsep\scshape\abstractname.]%
}
\makeatother
\numberwithin{equation}{section}

\newtheorem{thm}{Theorem}[section]
\newtheorem{lem}[thm]{Lemma}

\newtheorem{cor}[thm]{Corollary}
\newtheorem{prop}[thm]{Proposition}
\newtheorem{defn}[thm]{Definition}
\theoremstyle{definition} % upright text
\newtheorem{ex}[thm]{Example}
\newtheorem{rem}[thm]{Remark}
\newtheorem*{rem*}{Remark}

\def\a{\alpha}  % warning: Don't use with package beamer
\def\b{\beta}
\def\g{\gamma}
\def\Ga{\Gamma} % section

\let\euler=\epsilon
\let\epsilon=\varepsilon
\def\k{\kappa}
\def\l{\lambda}

\def\m{\mu}
\def\om{\omega}
\def\Om{\Omega}

\def\[{\begin{equation}}  % use $$ ... $$ for equations without numbers
\def\]{\end{equation}}
\newcommand{\<}{\langle}
\renewcommand{\>}{\rangle}

\DeclareMathOperator\ad{ad}
\newcommand{\Ber}{{\mathrm{Ber}}}  %  Berezinian Bundle

\newcommand{\bbC}{{\mathbb{C}}}
\def\C{\bbC}

\newcommand{\iCinfty}{{\mathrm{C}^\infty}}
\newcommand{\Cinfty}{{R}}
\newcommand{\smooth}{{\mathrm{C}^\infty}}

\newcommand{\Cl}{{\mathrm{Cl}}}
\newcommand{\bbCl}{{\mathrm{\bbC l}}}

\def\conn_#1{\nabla_{\!#1\,}}  % better spacing for nabla with only a lower argument
\def\tildeconn_#1{\widetilde\nabla_{\!#1\,}}
\newcommand\CTor[1][\nabla]{{\mathrm{C}_{#1}}}  % Torsion of a Courant algebroid with E-connection
\newcommand\cycl{{\text{cycl}}}  % cyclic permutations
\def\ud{{\mathrm{d}}}
\newcommand{\D}{{\mathcal{D}}}
\DeclareMathOperator\Div{Div}
\newcommand{\e}{{\mathrm{e}}}   % Exponential
\DeclareMathOperator\End{End}
\newcommand\GG{{\mathcal{G}}}  % bundle of quadratic Lie algebras
\DeclareMathOperator\Gr{gr}
\newcommand{\mg}{{\mathrm{g}}}  % metric
\newcommand{\half}{{\tfrac12}}

\newcommand{\calH}{{\mathcal{H}}}          %for connection/torsion on regular Courant algd
\newcommand{\bi}{{\mathsf{i}}}

\newcommand{\cI}{{\mathcal{I}}}
\DeclareMathOperator\Id{id}
\newcommand\ins{\iota}  % contraction of a form with a vector

\DeclareMathOperator\Lie{L}  % Lie derivative, has proper spacing
\newcommand{\ccL}{{\mathcal{L}}}
\newcommand{\cM}{{\mathcal{M}}}
\newcommand\M{\cM}
\newcommand{\cN}{{\mathcal{N}}}
\newcommand{\Enabla}{{\boldsymbol{\nabla}}}  % deprecated
\def\Econn_#1{\Enabla_{\!#1\,}}
\newcommand{\bbN}{{\mathbb{N}}}
\newcommand\N{\bbN}
\newcommand{\cO}{{\mathcal{O}}}
\def\gsmooth{\cO}

\newcommand{\cQ}{{\mathcal{Q}}}
\newcommand{\Qodd}{\gamma}  %{\mathcal{C}}
\def\br{{\mathbf{r}}}

\newcommand{\bbR}{{\mathbb{R}}}
\def\R{\bbR}

\DeclareMathOperator\rk{rk}

\newcommand{\cS}{{\mathcal{S}}}  % graded symmetric product/ power of a vector space/ module

\newcommand{\fS}{{\mathfrak{S}_\nabla}}

\DeclareMathOperator\Span{span}  % linear span of some vectors
\newcommand{\sS}{{\mathbb{S}}}

\newcommand{\cT}{{\mathcal{T}}}
\def\xto{\xrightarrow}
\newcommand\Top{{\text{top}}\,}
\DeclareMathOperator\Tr{Tr}

\newcommand{\ve}{\varepsilon}  % deprecated, use \epsilon instead
\newcommand{\Vect}{\mathfrak{X}}

\newcommand\Vol[1][M]{|{\wedge}^\Top T^*#1|}
\newcommand{\WQ}{{\mathcal{WQ}}}  % this is the one we define in 3
\newcommand{\bbZ}{{\mathbb{Z}}}
\newcommand{\aA}{{A[1]}}
\newcommand{\aaA}{{A^* [1]}}
\DeclareMathOperator{\pr}{pr}

\long\def\sideremark#1{\ifvmode\leavevmode\fi\vadjust{\vbox to0pt{\vss% the remark
 \hbox to 0pt{\hskip\hsize\hskip1em%                          will appear only
 \vbox{\hsize3cm\tiny\raggedright\pretolerance10000%          on the side
 \noindent #1\hfill}\hss}\vbox to8pt{\vfil}\vss}}}%

\begin{document}

\baselineskip=15pt

%%%%%%%%%%%%%%%%%%%%%%%%%%%%%%%%%%%%%%%%%%%%%%%%%%%%%%%%%%%%
%%%%%%%%%%%%%%%%%%%%%%%%% TITLE %%%%%%%%%%%%%%%%%%%%%%%%%%%%
\title{Weyl quantization of  degree 2 symplectic graded manifolds}
%%%%%%%%%%%%%%%%%%%%%%%%%%%%%%%%%%%%%%%%%%%%%%%%%%%%%%%%%%%%
%%%%%%%%%%%%%%%%%%%%%%%%%%%%%%%%%%%%%%%%%%%%%%%%%%%%%%%%%%%%

\author{Melchior Gr\"utzmann}
\address{Department of Physics, University of Jena, Germany}
\email{melchiorG@gmail.com}

\author{Jean-Philippe~Michel}
\address{Department of Mathematics, University of Li\`ege, Belgium}
\email{jpmichel82@gmail.com}

\author{Ping Xu}
\address{Department of Mathematics, Penn State University, United States}
\email{ping@math.psu.edu}

%\date{\today}

%% Classification math�matique  (2000)
%%%%Funds
\thanks{Research partially funded by the NSF grants
DMS-2001599, DMS-1707545, DMS-1406668, DMS-1101827, the Luxembourgian NRF
and the Interuniversity Attraction  Poles Program initiated by the Belgian Science Policy Office.}

	\begin{abstracts}
	\abstractin{french}

%\begin{otherlanguage}{french} 
%\begin{abstract}
Soit $S$ un fibré spinoriel d'un fibré vectoriel pseudo-Euclidien $(E,\mg)$ de rang pair.
Nous introduisons une nouvelle filtration de l'algèbre  $\D(M,S)$ des opérateurs différentiels 
sur $S$. Pour cette filtration, l'algèbre graduée associée $\Gr \D(M,S)$ s'avère être isomorphe à l'algèbre 
$\cO(\cM)$ des fonctions lisses sur $\cM$, la variété graduée symplectique de degré $2$
canoniquement associée à $(E,\mg)$.
En conséquence, nous construisons la quantification de Weyl de $\cM$ comme une application
$\WQ_\hbar:\cO(\cM)\to\D(M,S)$, et montrons que $\WQ_\hbar$ satisfait toutes les propriétés
voulues d'une quantification. En application, nous obtenons une bijection entre les structures
d'algébroïdes de Courant $(E,\mg,\rho,\llbracket\cdot,\cdot\rrbracket)$, caractérisées de manière
équivalentes par des fonctions hamiltoniennes génératrices sur la variété graduée symplectique $\cM$,
et les opérateurs de Dirac générateurs anti-symétriques $D\in\D(M,S)$.
L'opérateur $D^2$ est un nouvel invariant de $(E,\mg,\rho,\llbracket\cdot,\cdot\rrbracket)$,
qui généralise la norme au carré de la $3$-forme de Cartan d'une algèbre de Lie quadratique.
Nous étudions cet invariant en détail dans le cas particulier où $E$ est le double d'un bi-algébroïde
de Lie  $(A,A^*)$.

	\abstractin{english}
Let $S$ be a spinor bundle of a
pseudo-Euclidean vector bundle $(E,\mg)$ of even rank.
We introduce  a new filtration on the algebra $\D(M,S)$
of differential operators on $S$. As a  main property, the
associated  graded algebra $\Gr \D(M,S)$
is  proved to be isomorphic to the algebra  $\cO(\cM)$ of smooth functions
on  $\cM$, where $\cM$ is the degree $2$
symplectic graded manifold
canonically associated to $(E,\mg)$.
Accordingly, we establish  the Weyl quantization of $\cM$
 as a map $\WQ_\hbar:\cO(\cM)\to\D(M,S)$,
and prove that $\WQ_\hbar$ satisfies all the
 desired  properties of quantizations.
As an application, we obtain a bijection
 between Courant algebroid structures
$(E,\mg,\rho,\llbracket\cdot,\cdot\rrbracket)$
 equivalently characterized by  Hamiltonian generating  functions on $\cM$,
and skew-symmetric Dirac generating operators $D\in\D(M,S)$.
The operator $D^2$
gives a new invariant of $(E,\mg,\rho,\llbracket\cdot,\cdot\rrbracket)$
which generalizes the square norm of the Cartan $3$-form of a quadratic Lie algebra.
We study this  invariant
 in detail in the particular case of $E$ being the double of
a Lie bialgebroid $(A,A^*)$.
\end{abstracts}

\thispagestyle{empty}
\maketitle

\tableofcontents

%%%%%%%%%%%%%%
\section{Introduction}
%%%%%%%%%%%%%%

This paper is devoted to the study of Weyl quantization of
 degree 2 symplectic graded manifolds, and its application to
Courant algebroids.

In searching for the Lie algebroid analogue
of Drinfeld's  double of Lie bialgebras \cite{Dri86},
Liu--Weinstein--Xu introduced the notion of Courant
algebroids \cite{LWX97}. The Courant algebroid axioms were later
reformulated by Roytenberg  \cite{Roy99}
in terms of Dorfman brackets. However these axioms  still  remain
mysterious, and  various attempts have been made   in order
to understand Courant algebroids in a  more conceptual
and transparent  way.  One approach
  was through a degree 3  Hamiltonian function in a
  degree 2 symplectic graded manifold.
When a   Courant algebroid $E$ is  the  double $A\oplus A^*$ of
a Lie bialgebroid $(A, A^*)$, Roytenberg \cite{Roy99} proved that
the  Courant algebroid structure on $E$ is indeed equivalent to
a degree 3  Hamiltonian generating function $\Theta$ on
the symplectic graded manifold $T^*[2](A[1])$
satisfying the equation $\{\Theta, \Theta \}=0$.
Following an idea of Weinstein, Roytenberg \cite{Roy02} and
 \v{S}evera \cite{Sev7}
extended this result
to arbitrary Courant algebroids and proved that
there is essentially a bijection between Courant algebroids and
degree 3 Hamiltonian generating functions on
 degree 2 symplectic graded manifolds. By a graded manifold
we always mean a  $\bbN$-graded  manifold.

Around the same time, in 2001, in an unpublished
manuscript \cite{AX03}, Alekseev--Xu took a different
approach in terms of Dirac generating operators,
an analogue of Kostant's cubic Dirac operators \cite{Kos99}.
Alekseev--Xu's approach was motivated
by the following basic example due to Cabras--Vinogradov \cite{CVi92}.
Let $M$ be a manifold, and ${\mathfrak X}(M) \cong \Gamma(TM)$
be the Lie algebra of vector fields on $M$. The idea is to
extend the Lie bracket on ${\mathfrak X}(M)$ to sections of %
the bundle $E=TM\oplus T^*M$.
Observe that sections of $E$ act on the space of differential
forms $\Om (M)$ by contraction and by  exterior multiplication,
respectively,
$$%\begin{equation}
%\label{eq:action}
(X + \a ) \cdot  \mu := \iota_X \mu + \a  \wedge \mu,
$$%\end{equation}
where $\a  \in \Gamma(T^*M) = \Omega^1 (M)$,
  $X \in {\mathfrak X}(M)$ and
$\mu\in \Om (M)$. This action turns
$\Om (M)$ into a Clifford module of the Clifford
bundle $\Cl(E)$, where $E$ is equipped with the standard bilinear
form:
$$%\begin{equation}
%\label{eq:standard}
(X_1+\a_1 , X_2+\a_2) =
\langle \a_1, X_2   \rangle + \langle \a_2 ,  X_1\rangle,
$$%\end{equation}
and the Clifford  generating relation is $xy+yx=(x, y)$.

Using the  de Rham differential $\ud: \Omega^\bullet (M) \rightarrow
 \Omega^{\bullet+1}(M)$,
one can form the {\em derived bracket} \cite{Kos96} on sections of $E$:
\begin{equation}
\label{eq:pig}
\llbracket e_1, e_2\rrbracket := [[\ud, e_1], e_2],
\qquad \forall  e_1 , e_2 \in \Gamma (E),
\end{equation}
where both sides are viewed as operators on $ \Om (M)$.
 It is straightforward to check that Eq.~\eqref{eq:pig}
coincides with  the Dorfman bracket of the standard Courant
algebroid $TM \oplus T^*M$.

Observe that  $\Om (M)$ is indeed  a real
spinor bundle of $\big(TM\oplus T^*M , (\cdot, \cdot )\big)$.
For a general Courant algebroid $(E,\mg,\rho,\llbracket\cdot,\cdot\rrbracket)$,
Alekseev--Xu  \cite{AX03} proved that  there exist cubic Dirac type
generating operators, acting on a certain  spinor bundle of $(E,\mg)$, %$\D(M, S)$,
which play exactly the same role as the  de Rham differential operator
$\ud$  does in the Cabras--Vinogradov's approach to the standard Courant
algebroid\footnote{
To the best of our knowledge, Cabras-Vinogradov's contributions
to Courant brackets seem to have been overlooked by the community so far.}.
%  Here $S$ is a certain  spinor bundle of $(E,\mg)$
%i.e. $\End S\cong \bbCl (E)$, where $\bbCl (E)=\Cl(E)\otimes\bbC$ denotes the
%complex Clifford bundle.
%and  $\D(M, S)$  denotes the algebra of differential operators on $S$.
Thus a natural question arises: is  there any   relation
between Hamiltonian generating functions and Dirac generating operators
for a Courant algebroid?

To answer this question, we are naturally led to the study
of Weyl quantization on symplectic graded manifolds of  degree 2.
The classical Weyl quantization formula  is the prototypical example
of a quantization map:
%due to Weyl in connection of his attempt to understand the
%transition from classical mechanics to quantum mechanics \cite{Weyl}.
to each polynomial
function on the classical phase space $T^*\bbR^n$, it assigns
a differential operator on $\bbR^n$. More precisely,
%for polynomials in the coordinates $(x^i,p_i)$,
the Weyl quantization maps polynomials in the coordinates $(x^i,p_i)$
to differential operators on $\bbR^n$
with coordinates $(x^i )$, and is defined as the symmetrization
map such that $p_i \mapsto
\frac{\hbar}{\bi}\frac{\partial}{\partial x^i}$
and $x^i \mapsto m_{x^i}$ (multiplication by $x^i$).
% For this reason, Weyl quantization
%is also referred to as the  quantization of symmetric ordering.
Weyl quantization has been  generalized to various contexts
and plays an important  role in different branches of mathematics, e.g.\
harmonic analysis \cite{Fol89}, pseudo-differential symbolic
 calculus \cite{Hor79}, formal  and strict  deformation quantizations
\cite{BFFLS78, Rie89}.
More specifically, we are concerned with the extension of Weyl quantization
to smooth manifolds and supermanifolds.
Using an affine connection on $M$, Underhill  built a Weyl quantization on
the symplectic manifold $T^*M$ \cite{Und78}. 
More generally, by choosing an additional linear connection
on the  vector bundle $V\to M$, Widom \cite{Wid80}  obtained
 a quantization map of the form
$\cQ_\hbar^M:\cO(T^*[2]M)\otimes_\Cinfty\Ga(\End V)\to\D(M,V)$,
where $\D(M,V)$ is the algebra of differential operators on $V$ (see also \cite{Bok69}).
In the case when
 $V$ is the spinor bundle $S$ of a Riemannian spin manifold $(M,\mg)$,
Getzler \cite{Get83} used such a map $\cQ_\hbar^M$ to obtain
 a Weyl quantization on the symplectic
supermanifold $T^*M\oplus\Pi TM$ and then proved
the index theorem for the Riemannian Dirac operator.
See \cite{Vor99} for the case where $S=\Om(M)$.
Other related quantization schemes have been applied to even
 symplectic supermanifolds.
See e.g.\ \cite{Bor00,Mic10a,Mic12}.
%We refer to \cite{Mic10a} for the geometric quantization of $T^*M\oplus\Pi TM$
%and to \cite{Mic12} for its conformally equivariant quantization.
%Besides, Bordemann \cite{Bor00} has performed the Fedosov quantization
%of even symplectic supermanifolds, which are classified
%by Rothstein \cite{Rot90}.

In our situation, as a first step, by taking a metric preserving
linear connection on $E$,
we  can identify any  degree 2 symplectic graded manifold
as  $T^*[2]M\oplus E[1]$, where $E$ is a vector bundle
over~$M$ equipped with a pseudo-metric $g$. According to \cite{Rot90},
the degree 2 symplectic form on $T^*[2]M\oplus E[1]$ can be written
explicitly in terms of the canonical symplectic structure on
$T^*[2]M$, the  pseudo-metric $g$ and the chosen connection.
Our first main result is to construct a Weyl type quantization
for this graded symplectic manifold  $T^*[2]M\oplus E[1]$
by a  combination of Clifford quantization and classical
Weyl quantization. This generalizes the previous works \cite{Get83,Vor99}.
More precisely, we assume that $(E,\mg)$ is of even rank and admits
a spinor bundle $S$. That is,
 $\End S\cong \bbCl (E)$ with $\bbCl (E)=\Cl(E)\otimes\bbC$
being the complex Clifford bundle.
We introduce an increasing and exhaustive filtration of the algebra
 $\D(M,S) =\bigcup_{k\in\bbN}\D_k(M,S)$:
\begin{equation}
\label{eq:Lyon}
\D_0(M,S)\subset\D_1(M,S)\subset\cdots\subset\D_k(M,S)\subset\cdots
\end{equation}
satisfying the conditions: for any $k,l\in\bbN$, 
\begin{align*}
&\D_k(M,S)\cdot\D_l(M,S) \subseteq \D_{k+l}(M,S),%\qquad\forall k,l\in\bbN.\\
&\left[\D_k(M,S),\D_l(M,S)\right] \subseteq \D_{k+l-2}(M,S).%\qquad\forall k,l\in\bbN.
\end{align*}
%
%\[\nonumber
%\begin{split}
%D_k(M,S)\cdot\D_l(M,S) &\subseteq \D_{k+l}(M,S),\\
%\left[\D_k(M,S),\D_l(M,S)\right] &\subseteq \D_{k+l-2}(M,S),
% \end{split}
%\qquad\forall k,l\in\bbN.
%\]
As a consequence, the associated graded algebra
$
\Gr\D(M,S)=\bigoplus_{k\in\bbN}\D_{k}(M,S)/\D_{k-1}(M,S),
$
is a graded commutative Poisson algebra of degree $-2$, which we
prove to be isomorphic to $\cO^\bbC(T^*[2]M\oplus E[1])$.
Such an isomorphism enables us to define the principal symbol
$\sigma_k: \D_k(M,S)\to \cO^\bbC_k(T^*[2]M\oplus E[1])$, $\forall k\in\bbN$,
exactly in the same way as in the classical case.
We prove that the Weyl quantization map
 $\WQ_\hbar$ establishes an isomorphism
$\gsmooth^\bbC(T^*[2]M\oplus E[1])\stackrel{\sim}{\to} \D(M, S)$, which is
the right inverse of the principal symbol map, and satisfies
all the  desired properties. See Theorem \ref{thm:WQ}.
Note that the filtration \eqref{eq:Lyon}  on $\D(M,S)$ is different
 from both the usual filtration by the orders of 
differential operators and Getzler's filtration  \cite{Get83}.
In a certain sense, this is the only filtration which turns $\Gr\D(M,S)$ into a
graded commutative Poisson algebra,
with a  non-degenerate Poisson bracket. Roughly speaking, it assigns degree $2$
to derivations and degree $1$ to sections in $\Ga(E)\subset\Ga(\bbCl(E))$.

The second part of the paper is devoted to the application
of Weyl quantization to Courant algebroids.
The consideration of a specific spinor bundle $\sS$,
obtained by twisting $S$ by a certain line bundle, allows us
to define a conjugation map and an adjoint operation on~$\D(M,\sS)$.
A Dirac generating operator is then defined as a real operator in $\D_3(M,\sS)$,
which is odd and squares to a function on the base manifold.
Following an idea of \v{S}evera \cite{Sev7}, we prove that, for a
given Courant algebroid, there exists a unique
 skew-symmetric Dirac generating operator,
and moreover  the Weyl quantization map $\WQ$
establishes a bijection between Hamiltonian
generating functions  and skew-symmetric Dirac generating operators.
This is our second main result.

As an application, by considering the square of the unique
skew-symmetric Dirac generating operator,
 we obtain a new Courant algebroid invariant.
This new invariant, as a function
on the base manifold, is a natural extension
of the square norm of the Cartan 3-form of
a quadratic Lie algebra.
As another consequence,
we recover a result of Chen--Sti\'enon \cite{CS09}
regarding Dirac generating operators for Lie bialgebroids,
which gives an equivalent description of the Lie bialgebroid
compatibility condition. In this case, $E=A\oplus A^*$ and
there are two natural twisted spinor bundles,
each of them admitting a  Dirac generating operator.
For the Lie bialgebroid
arising from  a generalized complex structure,
they coincide with
the $\partial$ and $\overline{\partial}$-operators
of the generalized complex structure \cite{Che09}.

Some remarks are in order.
 We learned recently that Li-Bland and Meinrenken also
 obtained  similar results
in their study of Dirac generating operators~\cite{LMe12}.

\subsection*{Notations}
Finally, we list the notations used throughout the paper.

$\bbN=\{0,1,2,\ldots\}$ denotes the set of non-negative integers,
$\bbN^\times=\{1,2,\ldots\}$ the set of positive integers and $\bi=\sqrt{-1}$.
Tensor products over the algebra of real numbers are denoted by
$\otimes$ or $\otimes_\bbR$, whereas tensor products
over the algebra $\Cinfty:=\smooth(M)$
%of smooth functions on a manifold $M$
are denoted by  $\otimes_\Cinfty$.
For a vector space (or a vector bundle) $V$, symmetric and skew-symmetric
tensor products  are denoted by $\cS V$ and $\wedge V$, respectively.
We use the Einstein's summation convention without further comments.

\subsection*{Acknowledgments}
We would like to thank several institutions for their hospitality
 while work on this project was being done:
University of Luxembourg (Gr\"utzman),
Penn State University (Michel),
Universit\'e Catholique de Louvain (Xu),  and Universit\'e Paris Diderot (Xu).
We also wish to thank many people for useful discussions and comments,
 including Rajan Mehta, Dmitry Roytenberg, Mathieu Sti\'enon,  Yannick Voglaire, and in particular
 Pavol \v{S}evera for sharing many of his unpublished
letters. Finally, we are very grateful to Anton Alekseev, who contributes
to numerous key concepts and ideas in this paper through the unpublished
manuscript  \cite{AX03}, and should be considered
as a “virtual co-author” of this paper.

%%2 %%%%%%%%%%%%%%%%%%%%%%%%%%%%%%%%%%%%%%%%%%%%%%%%
%%%%%%%%%%%%%%%%%%%%%%%%%%%%%%%%%%%%%%%%%%%%%%%%%%%%
\section{Symplectic graded manifolds of degree $2$}
%%%%%%%%%%%%%%%%%%%%%%%%%%%%%%%%%%%%%%%%%%%%%%%%%%%%
%%%%%%%%%%%%%%%%%%%%%%%%%%%%%%%%%%%%%%%%%%%%%%%%%%%%

We recall in this section some
 standard materials concerning symplectic graded manifolds.
Our presentation is mainly  based  on \cite{Roy02,Roy06},
 enriched with the work of Rothstein on symplectic supermanifolds \cite{Rot90}. \nocite{Gru09}

%%% 2.1 %%%%%%%%%%%%%%%%%%%%%%%%%%%%%%%%%%%%%%%%%%%%
\subsection{Definition}
%%%%%%%%%%%%%%%%%%%%%%%%%%%%%%%%%%%%%%%%%%%%%%%%%%%%

A \emph{graded manifold} $\cM$ is a smooth manifold $M$ endowed with
 a sheaf of $\bbN$-graded algebras $\cO$ such that, for every
contractible open set  $U\subset M$, the algebra of functions $\cO(U)$
 is isomorphic to $\iCinfty(U)\otimes \cS V$,
for a fixed $\bbN^\times$-graded vector space
 $V=\bigoplus_{i\in\bbN^\times}V_i$.
Here $\cS V$
%=\bigotimes_{i\in 2\bbN^\times}\cS V_i\otimes\wedge V_{i-1}$
denotes the graded symmetric tensor algebra of~$V$
and the grading of $\cO(U)$ is induced by the one of $V$.
In particular, the degree $0$ component of $\cO(U)$ is $\iCinfty(U)$.
The algebra of global sections of $\cO$ is
called the algebra of functions on $\cM$ and denoted by~$\cO(\cM)$.
 Coordinates on $U$ together with a graded basis of  $V$
form  a local coordinate system $(x^i)$ on $\cM$.
The grading on the algebra of functions, $\cO(\cM)=\bigoplus_{k\in\bbN}\cO_k(\cM)$,
is then given by the decomposition of $\cO(\M)$ into
the  eigenspaces of the
 Euler vector field on $\cM$:
$$%\[
  \euler = w(x^i)x^i\frac{\partial}{\partial {x^i}},
$$%\]
where $w(x^i)\in\bbN$ denotes the degree of the coordinate $x^i$,
i.e., $w(x^i)=0$ if $x^i$ is a coordinate on $U$ and $w(x^i)=j$
if $x^i\in V_j$.

 Given  a smooth vector bundle  $E\to M$ and
a positive integer $k\in\N^\times$,
by $E[k]$, we  denote the graded manifold with base $M$,
 whose algebra of functions is
 $$\cO(E[k])=  \begin{cases}  \Gamma(\wedge E^*)  \quad\text{for }k\text{ odd},\\
   \Gamma(\cS E^*)  \quad\text{for }k\text{ even.}
 \end{cases}$$
Here sections of $E^*$ are assigned the degree $k$.
The graded manifold $E[k]$ is also called a shifted vector bundle.
We will mostly consider graded manifolds built out of shifted vector bundles.

A \emph{symplectic graded manifold} of degree $n$ is a graded manifold
$\cM$ endowed with a symplectic form of degree $n$, i.e.\ a
 closed non-degenerate 2-form $\om$, whose Lie derivative along
the Euler vector field $\euler$ satisfies $\Lie_\euler\om=n\om$.
The algebra of functions $\cO(\cM)$ admits a Poisson bracket of degree $-n$,
 i.e. $\{\cO_k(\cM),\cO_l(\cM)\}\subset\cO_{k+l-n}(\cM)$  for
all $k,l\in\bbN$. The Poisson  bracket is graded skew-symmetric
and satisfies the graded Jacobi identity, namely
\begin{align*}
\{F,G\} &= -(-1)^{(k-n)(l-n)}\{G,F\},\\
\{F,\{G,H\}\} &= \{\{F,G\},H\}+(-1)^{(k-n)(l-n)}\{G,\{F,H\}\},
\end{align*}
for all $F\in\cO_k(\cM)$, $G\in\cO_l(\cM)$ and $H\in\cO(\cM)$.

%%% 2.2 %%%%%%%%%%%%%%%%%%%%%%%%%%%%%%%%%%%%%%%%%%%%
\subsection{Example}
%%%%%%%%%%%%%%%%%%%%%%%%%%%%%%%%%%%%%%%%%%%%%%%%%%%%

In \cite{Rot90}, Rothstein gives a description of symplectic
 supermanifolds in terms of the following data: a pseudo-Euclidean
vector bundle $(E,\mg)$ over an ordinary symplectic manifold
and a metric connection $\nabla$ on $(E,\mg)$, i.e., a (linear)
 connection on $E$ satisfying $\nabla\mg=0$.
In the sequel, we will adapt his construction
to the graded context and  obtain a symplectic structure of
degree~$2$ on the Whitney sum $T^*[2]M\oplus E[1]$  over $M$.

\begin{prop}[\cite{Rot90}]\label{PropRoth}
Let $E\to M$ be a smooth vector bundle, endowed with a pseudo-Euclidean metric $\mg$
and a metric connection $\nabla$. Then, the graded manifold $T^*[2]M\oplus  E[1]$
admits an exact symplectic $2$-form of degree $2$:
\[\label{SymplForm}
\om_{\mg,\nabla}:=\ud\a \qquad \text{where} \qquad \a=\pi_1^*\a_0+\pi_2^*\b.
\]
Here $\pi_1$ and $\pi_2$ are the canonical projections on $T^*[2]M$ and $E[1]$
respectively, $\a_0$ is the Liouville $1$-form on $T^*[2]M$
  and $\b\in\Om^1(E[1])$ is the $1$-form on $E[1]$ which
 annihilates  the horizontal subspace of $T (E[1])$, corresponding to
the connection $\nabla$, and satisfies
 $\b_e(v_e)=\frac{1}{2}\mg_{\pi(e)}(v_e,e)$ for all 
vertical tangent vectors $v_e\in T_e (E[1])$.
\end{prop}

In what follows, we make  repeated use of the identification $E\cong E^*$,
induced by the metric $\mg$,
and  of  the identifications below,  without mentioning them explicitly:
\begin{align*}
\cO_0(T^*[2]M\oplus E[1])&\cong\iCinfty(M),\\
\cO_1(T^*[2]M\oplus E[1])&\cong\Gamma(E),\\
\cO_2(T^*[2]M\oplus E[1])&\cong\Gamma(TM\oplus\wedge^2 E).
\end{align*}
The spaces $\iCinfty(M)$, $\Ga(E)$ and $\Vect(M)\cong \cO_2(T^*[2]M)$  generate the
 entire algebra of functions on~$T^*[2]M\oplus E[1]$. Hence,
 by the Leibniz rule, the  symplectic structure
on $T^*[2]M\oplus E[1]$ can be characterized in terms of the following
 Poisson brackets:
\begin{align}
 \{X,f\}   &=X(f),      & \{h,f\}   &=0, \nonumber\\
 \{X,\xi\} &=\conn_X\xi,& \{f,\xi\} &=0, \label{PoissonRothstein}\\
 \{X,Y\}   &=[X,Y]+R^E (X,Y), &\qquad  \{\xi,\eta\} &=\mg(\xi,\eta), \nonumber
\end{align}
where $X,Y\in{\mathfrak X}(M)$, $\xi,\eta\in\Ga(E)$,  $f,h\in\iCinfty(M)$,
and $R^E$ denotes the curvature of the connection $\nabla$.
Since $\nabla$ is a metric connection,  $R^E (X,Y)$ defines an element
in $\Ga(\wedge^2 E)$, i.e. a degree $2$ function on $T^*[2]M\oplus E[1]$.
From \eqref{PoissonRothstein},
 it is simple to see that the bracket
 $\{\cdot,\cdot\}$ is indeed of degree $-2$.

\begin{rem}\label{rem:PB_E}
According to  \eqref{PoissonRothstein},
$E[1]$ is a Poisson submanifold of $T^*[2]M\oplus E[1]$,
whose Poisson bracket is induced by the pseudo-metric $\mg$.
\end{rem}

For a vector bundle $E$ over $M$, the cotangent bundle $T^*[2](E[1])$
 is naturally a
degree $2$ symplectic graded manifold.
Let 
$$\pi_E:T^*[2](E[1])\to E[1]$$ and 
$$\pi_{E^*}:T^*[2](E[1])\cong T^*[2] (E^*[1])\to E^*[1]$$
be the natural  projections.
They  combine into a map:
\[\label{def:tilde_pi}
 \tilde\pi:T^*[2](E[1])\longrightarrow (E\oplus E^*)[1].
\]
A connection $\nabla$ on $E$ gives rise to a horizontal distribution
on the tangent bundle $TE$, which in turn induces a surjective submersion
\[\label{def:pi_nabla}
\pi_\nabla:T^*[2](E[1])\longrightarrow T^*[2]M.
\]
It is simple to see that $\pr\circ \tilde\pi=\pr\circ\pi_\nabla$,
where, by abuse of notation, $\pr$ denotes both natural projections
$\pr:  (E\oplus E^*)[1]\to M$ and $T^*[2]M\to M$.
Putting together the maps $\tilde\pi$ and $\pi_\nabla$, we obtain a map:
\begin{equation}
\label{eq:HKG}
\widetilde{\Xi}_{\nabla}:  T^*[2] (E[1]) \longrightarrow T^*[2]M\oplus (E\oplus E^*)[1].
\end{equation}
It is simple to check that $\widetilde{\Xi}_{\nabla}$ is
indeed  a diffeomorphism.
Note that $E\oplus E^*$  is a
pseudo-Euclidean bundle over $M$ with the  duality pairing.
 The connection $\nabla$ on $E$ induces a connection on~$E\oplus E^*$, which
is compatible with the duality pairing.
 According to Proposition~\ref{PropRoth}, these data induce
 a degree $2$ symplectic structure on $T^*[2]M\oplus (E\oplus E^*)[1]$.

\begin{lem}\label{lem:splitting}
The map $\widetilde{\Xi}_{\nabla}$ in \eqref{eq:HKG}
 is a symplectic diffeomorphism.
\end{lem}
\begin{proof}
Since $\widetilde{\Xi}_{\nabla}$ is a diffeomorphism, it suffices
to prove that  it is a Poisson map.
As $\cO(T^*[2]M\oplus (E\oplus E^*)[1])$ is
generated by $\iCinfty(M)$, $\Ga(E\oplus E^*)$ and ${\mathfrak X}(M)$,
it  is thus sufficient to check that  $\widetilde{\Xi}_{\nabla}^*$
preserves  the  Poisson brackets between elements of these  spaces.

It is simple to see that
\begin{eqnarray}\nonumber
(\widetilde{\Xi}_{\nabla})^*f=(\pi_E)^*f,& &\quad (\widetilde{\Xi}_{\nabla})^*\eta=(\pi_E)^*\eta,\\ \label{Eq:fxietaX}
(\widetilde{\Xi}_{\nabla})^*X=(\pi_\nabla)^*{X},& &\quad (\widetilde{\Xi}_{\nabla})^*\xi=(\pi_{E^*})^*\xi,
\end{eqnarray}
for any $f\in\iCinfty(M)$, $\eta\in\Ga(E)$, $\xi\in\Ga(E^*)$ and $X\in\Vect(M)$.
Here, both  $(\pi_E)^*f$ and $(\pi_E)^*\eta$
are fiberwise constant functions on $T^*[2](E[1])$,
while   both $(\pi_\nabla)^*{X}$ and $(\pi_{E^*})^*\xi$ are fiberwise
 linear functions
 on $T^*[2](E[1])$, corresponding to the vector fields
$\conn_X$ and $\iota_\xi=\mg(\xi,\cdot)$ on $E[1]$, respectively.
Hence, the Poisson brackets in~$T^*[2](E[1])$ of  the four 
types of functions in  \eqref{Eq:fxietaX}
%in~$T^*[2]E[1]$ of
%$(\pi_E)^*f, (\pi_E)^*\eta, (\pi_{E^*})^*\xi, (\pi_\nabla)^*{X}$
%$(\pi_E)^*f$, $(\pi_E)^*\eta$, $(\pi_{E^*})^*\xi$ and $(\pi_\nabla)^*{X}$
 are easily computed.  On the other hand,
for the graded symplectic manifold ~$T^*[2]M\oplus (E\oplus E^*)[1]$,
the Poisson brackets of %$f,\eta,\xi,X$
$f$, $\eta$, $\xi$ and $X$ are given by  \eqref{PoissonRothstein},
with $E$, $\mg$ and $\nabla$ being replaced by $E\oplus E^*$,
the duality pairing and the induced connection, respectively.
The conclusion follows from a straightforward verification.
\end{proof}

%% 2.3 %%%%%%%%%%%%%%%%%%%%%%%%%%%%%%%%%%%%%%%%%%%%%
\subsection{Classification of symplectic graded manifolds of degree $2$}\label{Par:SymplGM}
%%%%%%%%%%%%%%%%%%%%%%%%%%%%%%%%%%%%%%%%%%%%%%%%%%%%

According to \cite{Roy02},
any  pseudo-Euclidean vector bundle $(E,\mg)$ determines a degree $2$
 symplectic graded manifold $\cM$,  which is defined to be the fiber product
 $\big(T^*[2](E[1]) \big) \times_{ (E\oplus E^*)[1]} E[1]$. Thus,
we have    the following commutative  diagram:
\[\label{M-Roytenberg}
 \xymatrix{
  \M\ar[r]^{i_\M\quad}\ar[d]^{\pi} & T^*[2]E[1]\ar[d]^{\tilde\pi}  \\
  E[1] \ar[r]^{i\quad} & (E\oplus E^*)[1]
 }
\]
Here $\tilde\pi$ is defined as in \eqref{def:tilde_pi},
 and $i$ is the diagonal-like embedding
 $\psi\mapsto\psi\oplus\tfrac12\mg(\psi,\cdot)$.
%essentially the diagonal map:$E[1]\to (E\oplus E^*)[1]$.
Since $i$ is an
isometric embedding, $i_\M$ must be  an embedding as well.
It is simple to see that the restriction
to $\M$ of the canonical
symplectic form on $T^*[2](E[1])$ is %indeed
non-degenerate.
Therefore $\M$ is a symplectic submanifold of
$T^*[2] (E[1])$.
It is known that, up to an isomorphism, every degree $2$ symplectic
graded manifold indeed arises in this way. We refer the interested
reader to \cite{Roy02} for details.

Let $\nabla$ be a  connection on $E$.
Then composing $i_\M$ with the map $\pi_\nabla$, as defined
 in \eqref{def:pi_nabla}, we obtain a map
$$
\pi_\nabla\circ i_\M:\M\longrightarrow  T^*[2]M.
$$
Together with the natural projection
 $\pi:\cM\rightarrow E[1]$,   we obtain a map
$$
\Xi_\nabla:\M\longrightarrow T^*[2]M\oplus E[1].
$$
It is simple to check that $\Xi_{\nabla}$ is  a diffeomorphism.
The following result is known to experts and sketched in \cite{Roy02}.

\begin{thm}
\label{thmRoyM}
Let $\nabla$ be a metric connection on $(E,\mg)$.
Then the map $\Xi_\nabla$ is a symplectic diffeomorphism,
where $T^*[2]M\oplus E[1]$ is endowed with the symplectic structure \eqref{SymplForm}.
\end{thm}

\begin{proof}
%The metric connection $\nabla$ on $E$ induces a
%connection on $E\oplus E^*$, which is compatible
%with the natural duality  pairing.
%According to Proposition \ref{PropRoth},
%$T^*[2]M\oplus(E\oplus E^*)[1]$ is a degree 2 graded symplectic
%manifold. On the other hand, $T^*[2]M\oplus E[1]$
%is also a degree 2 graded symplectic manifold.
Consider the map
$$
i_T=\Id_{T^*[2]M}\oplus i :T^*[2]M\oplus E[1]\longrightarrow T^*[2]M\oplus(E\oplus E^*)[1].
$$
It is simple to check that
$i_T$ is a symplectic embedding.
By definition, we have
\[\label{eq:xi}
\widetilde{\Xi}_{\nabla}=\pi_{\nabla}\oplus\tilde\pi\quad\text{ and }\quad
\Xi_{\nabla}=\pi_{\nabla}\circ i_\M\oplus\pi.
\]
According to the relation
 $\tilde\pi\circ i_\M=i\circ\pi$  (see diagram \eqref{M-Roytenberg}),
we have the following commutative diagram:
$$
\xymatrix{
  T^*[2] E[1]  \ar^{\widetilde{\Xi}_{\nabla}\qquad}[rr]& & T^*[2]M\oplus (E\oplus E^*)[1]  \\
\M \ar^{i_\M}[u] \ar^{\Xi_{\nabla}\qquad}[rr]   &  & T^*[2]M\oplus E[1] \ar_{i_T}[u]
}
$$
Since both $i_\M$ and $i_T$ are  symplectic embeddings
and $\widetilde{\Xi}_\nabla$ is a symplectic diffeomorphism,
it follows that
$\Xi_\nabla$ must be a symplectic diffeomorphism.
\end{proof}

As a consequence, the symplectic graded manifolds
 $(T^*[2]M\oplus E[1],\om_{\mg,\nabla})$, associated to
 different metric connections, are all isomorphic.
They provide minimal symplectic realization
 of the Poisson manifold $E[1]$  \cite{Roy02}.

%% 3 %%%%%%%%%%%%%%%%%%%%%%%%%%%%%%%%%%%%%%%%%%%%%%%%%%%
%%%%%%%%%%%%%%%%%%%%%%%%%%%%%%%%%%%%%%%%%%%%%%%%%%%%%%%%
\section{The algebra $\D(M,S)$ of spinor differential operators}
%%%%%%%%%%%%%%%%%%%%%%%%%%%%%%%%%%%%%%%%%%%%%%%%%%%%%%%%
%%%%%%%%%%%%%%%%%%%%%%%%%%%%%%%%%%%%%%%%%%%%%%%%%%%%%%%%

In this section, after recalling some basic materials regarding
Clifford algebras and spinor bundles (see e.g.\ \cite{BGV92, Mei13, Tra08}),
we introduce a new filtration on the algebra $\D(M,S)$
of spinor differential operators
 and determine its associated graded Poisson algebra.
By considering a specific spinor bundle~$\sS$, we obtain, in addition,
two involutions on $\D(M,\sS)$.

%%% 3.1 %%%%%%%%%%%%%%%%%%%%%%%%%%%%%%%%%%%%%%%%%%%%%%%%
\subsection{Clifford and spinor bundles}\label{Par:CliffSpin}
%%%%%%%%%%%%%%%%%%%%%%%%%%%%%%%%%%%%%%%%%%%%%%%%%%%%%%%%

Let $(E,\mg)$ be a pseudo-Euclidean vector bundle of even rank
over a smooth manifold $M$. The real \emph{Clifford bundle}
$\Cl(E)$ is a bundle of associative algebras, whose fiber at $x\in M$
is isomorphic to the real Clifford algebra
$$
\Cl(E_x):=\left(\bigotimes E_x/\cI\right),
$$
where $\cI$ is the ideal generated by
$\xi_1(x)\otimes\xi_2(x)+\xi_2(x)\otimes\xi_1(x)-\mg(\xi_1(x),\xi_2(x))$,
$\forall \xi_1,\xi_2\in\Ga(E)$. In the sequel, we mainly use
the complex Clifford bundle $\bbCl(E):=\Cl(E)\otimes\bbC$
and refer to it simply as the Clifford bundle.

The Clifford bundle inherits a natural filtration,
$\bbCl(E)=\bigcup_{k\in\bbN}\bbCl_k(E)$,
and a natural  $\bbZ_2$-grading,
$$
\bbCl(E)=\bbCl^{+}(E)\oplus \bbCl^{-}(E),
$$
where $\Ga(\bbCl_k(E))$ is spanned by products of at most $k$ sections of
 $E$, and $\bbCl^{+}(E)$ (resp.\ $\bbCl^{-}(E)$) is spanned by  products of  even
(resp.\ odd) number of sections of $E$.
 Through the metric~$\mg$, the algebra $\Ga(\wedge E\otimes \bbC)$
 can be identified with $\cO^\bbC(E[1])$, the algebra of complex valued
functions on $E[1]$. Moreover, $\mg$ induces a Poisson bracket
 $\{\cdot,\cdot\}$ on $E[1]$ (see Remark \ref{rem:PB_E}).
There is a  standard $\bbC$-linear isomorphism,
called the Clifford quantization map:
\[\label{Chevalley}
\g :\cO^\bbC(E[1])\longrightarrow\Ga(\bbCl(E)).
\]
It extends the canonical embedding $\Ga(E)\hookrightarrow\Ga(\bbCl(E))$
by skew-symmetrization and satisfies
\[\label{Eq:gammaequiv}
\g(\{\mu,\cdot\})=[\g(\mu),\g(\cdot)], \qquad\forall\mu\in\Ga(E\oplus\wedge^2 E).
\]
%where $\{\cdot,\cdot\}$ stands for the Poisson bracket on $E[1]$ induced by $\mg$.
We refer the interested reader to  \cite{BGV92, Mei13, Mic12} for more details.

By a \emph{spinor bundle},
we mean   a complex vector bundle $S$ over $M$ such that
 $\End S\cong \bbCl(E)$.

\begin{ex}\label{ex:wedgeA}
Let $A$ be a vector bundle, and  $E=A\oplus A^*$. Let
$\mg$ be  the following pairing on $E$:
$$
\mg( \zeta_1+\eta_1, \zeta_2+\eta_2)=\langle \zeta_1,\eta_2 \rangle+\langle \zeta_2,\eta_1 \rangle,
\qquad \forall \zeta_1,\zeta_2\in\Ga(A), \eta_1,\eta_2\in\Ga(A^*),
$$
where $\langle\cdot,\cdot\rangle$ denotes the duality pairing.
Then  $\wedge A^*\otimes\bbC$ (or $\wedge A\otimes\bbC$)
is  a  spinor bundle of~$(E,\mg)$,
and the Clifford action is given by
\[\label{Eq:CliffordAction1}
\g(\zeta)\phi= \iota_\zeta\phi
 \quad\text{and}\quad \g(\eta)\phi=\eta\wedge \phi, \qquad\forall \zeta\in \Gamma(A), \eta\in  \Gamma(A^*), \phi\in \Gamma (\wedge A^* \otimes\bbC).
\]
Note that $\wedge A^*$ is a real spinor bundle, i.e., 
$\End (\wedge A^* )\cong\Cl(E)$.
\end{ex}

\begin{rem}
%Not every pseudo-Euclidean vector bundle $(E, \mg)$ admits
%a spinor bundle.
There are certain topological obstructions to the existence
of a spinor bundle, which depend on the signature of the metric $\mg$.
The existence of a real spinor bundle imposes a further topological constraint
on the pseudo-Euclidean vector bundle $(E,\mg)$.
\end{rem}

 In what follows, we  will always assume that a spinor bundle exists,
 and  use the algebra isomorphisms below  without mentioning
them explicitly
$$
\Ga(\wedge E\otimes \bbC)\cong\cO^\bbC(E[1]) \qquad\text{and} \qquad\Ga(\End S)\cong\Ga(\bbCl(E)).
$$
%The existence of a spinor bundle on $M$ imposes in general restrictions on the topology of $E$.

Let $\nabla^E$ be a connection on $E$. It induces a connection
on the exterior bundle $\wedge E$  compatible with
the wedge product, i.e.,
$\conn_X^E(\xi_1\wedge\xi_2)=(\conn_X^E\xi_1)\wedge\xi_2+\xi_1\wedge(\conn_X^E\xi_2)$,
$\forall X\in{\mathfrak X}(M),\xi_1,\xi_2\in\Ga(\wedge E)$.
If~$\nabla^E$ is a metric connection, it also induces
a connection on the Clifford bundle $\bbCl(E)$  compatible
with the Clifford multiplication. The quantization map \eqref{Chevalley}
 intertwines the induced connections on $\wedge E\otimes \bbC$ and $\bbCl(E)$.
 That is,
$$
\g(\conn_X^E\eta)=\conn_X^E\g(\eta), \qquad\forall
\eta\in\cO^\bbC(E[1]), X\in{\mathfrak X}(M).
$$

\begin{defn}
Assume that $(E,\mg)$ admits a spinor bundle $S$.
A connection $\nabla^S$ on $S$ is called a spinor connection if there
 exists a metric connection $\nabla^E$ on $E$ satisfying the following compatibility condition
\[\label{compatibilitySE}
[\conn_X^S,\g(\eta)]=\g(\conn_X^E\eta),
\]
for all $X\in\Vect(M)$ and $\eta\in\cO^\bbC(E[1])$. In this case, $(\nabla^E, \nabla^S)$ is called a compatible pair of connections.
\end{defn}

For compatible connections $(\nabla^E, \nabla^S)$, the induced connections on $\bbCl(E)$ and $\End S$ are intertwined  by the isomorphism $\bbCl(E)\cong\End S$.

\begin{lem}
\label{lem:CompatibleConnections}
Assume that $(E,\mg)$ admits a spinor bundle $S$. Then,
\begin{enumerate}[label=(\roman*)]
\item there always exist compatible connections $(\nabla^E, \nabla^S)$;
\item there exists  $r\in\Om^2(M)\otimes\bbC$ such that the curvature
$R^E$ of $\nabla^E$ and the curvature $R^S$ of $\nabla^S$ satisfy
\[\label{compatibilityR}
R^S=\g(R^E)+r\Id_S;
\]
\item any two pairs of compatible connections $(\tilde{\nabla}^E, \tilde{\nabla}^S)$
and  $(\nabla^E, \nabla^S)$  are related as follows
$$
\tilde{\nabla}^E-\nabla^E=\{\varpi,\cdot\}\quad\text{ and }\quad \tilde{\nabla}^S-\nabla^S=\g(\varpi)+\nu\Id_S,
$$
where $\varpi\in\Om^1(M,\wedge^2 E)$ and $\nu\in\Om^1(M)\otimes\bbC$.
Here, $\{\cdot,\cdot\}$ stands for the Poisson bracket on $E[1]$.
% and we have $\{\varpi(X),\xi\}=-\iota_\xi \varpi(X)$, $\forall X\in{\mathfrak X}(M),\xi\in\Ga(E)$.
\end{enumerate}
\end{lem}
\begin{proof}
The existence of a metric connection  $\nabla^E$ is classical.
Such a connection induces a connection
on $\bbCl(E)\cong \End S$.
According to \cite{Tra08}, the latter is always
 induced by a connection $\nabla^S$ on $S$. Thus
$(\nabla^E, \nabla^S)$ is a pair of compatible connections.

Let $X,Y\in{\mathfrak X}(M)$. Since $\nabla^E$ preserves the metric,
its curvature can be identified with $R^E\in\Om^2(M,\wedge^2 E)$,
via the following relation
$$
\left(\conn_X^E\conn_Y^E-\conn_Y^E\conn_X^E-\conn_{[X,Y]}^E\right)\xi=\{R^E(X,Y),\xi\},
\qquad\forall \xi\in\Ga(E).
$$
The same relation holds by replacing $\xi\in\Ga(E)$
with  $\eta\in\Ga(\wedge E\otimes \bbC)$.
Using the compatibility condition \eqref{compatibilitySE}, we deduce that
$$
[R^S(X,Y),\g(\eta)]-\g\left(\{R^E(X,Y),\eta\}\right)=0, \qquad \forall \eta\in\Ga(\wedge E\otimes \bbC).
$$
Hence, by Eq.\ \eqref{Eq:gammaequiv}, we have
$$
[R^S(X,Y)-\g(R^E(X,Y)),\g(\eta)]=0, \qquad \forall \eta\in\Ga(\wedge E\otimes \bbC).
$$
The center of $\Ga(\bbCl(E))$ being the algebra $\iCinfty(M)\otimes \bbC$,
Eq.\ \eqref{compatibilityR} thus  follows.

We now compare two pairs of compatible connections
$(\tilde{\nabla}^E, \tilde{\nabla}^S)$ and  $(\nabla^E, \nabla^S)$.
First, we have $\tilde{\nabla}^E-\nabla^E=\mu$,
for some $\mu\in\Om^1(M,\End E)$. As both connections preserve the metric~$\mg$,
so does $\mu$ and we have $\m=\{\varpi,\cdot\}$ with $\varpi\in\Om^1(M,\wedge^2 E)$.
Similarly, we have $\tilde{\nabla}^S-\nabla^S=A$ with $A\in\Om^1(M,\End S)$.
The compatibility condition \eqref{compatibilitySE} implies that $[A,\cdot]=[\g(\varpi),\cdot]$
on~$\Ga(\End S)$. This yields $A=\g(\varpi)+\nu$
for some $\nu\in\Om^1(M)\otimes\bbC$.
\end{proof}

%%% 3.2 %%%%%%%%%%%%%%%%%%%%%%%%%%%%%%%%%%%%%%%%%%%%%%%%
\subsection{The filtered algebra $\D(M,S)$ and its associated graded algebra}\label{ParSymbol}
%%%%%%%%%%%%%%%%%%%%%%%%%%%%%%%%%%%%%%%%%%%%%%%%%%%%%%%%

Assume that  $(E,\mg)$ admits a spinor bundle $S$.
The algebra $\D(M,S)$ of differential operators on $S$ is
 a subalgebra of $\End(\Ga(S))$
 generated by $\Ga(\End S)\cong\Ga(\bbCl(E))$ and the covariant derivatives $\conn_X$,
where $\nabla$ is any connection on $S$ and $X$ ranges over all vector fields on $M$.
There is a $\bbZ_2$-grading on $\D(M,S)$, inherited from $\bbCl(E)$,
\[\label{ParityD}
\D(M,S)=\D^{+}(M,S)\oplus\D^{-}(M,S),
\]
where $\D^{\pm}(M,S)$ is generated by the covariant derivatives
 $\conn_X$ and $\Ga(\bbCl^{\pm}(E))$. By $[\cdot,\cdot]$, we denote the
 graded commutator:
$$
 [A,B] = A\circ B-(-1)^{|A|\,|B|}B\circ A,
$$
 for any $A,B\in\D(M,S)$ of  the parity $|A|,|B|\in\{0,1\}$, respectively.

We now  introduce a sequence of subspaces of $\D (M,S)$ inductively as follows.
 Set $\D_{-1}(M,S):=\{0\}$, and
 \[\label{Groth}
\D_k(M,S):=\{D\in\D(M,S)| \; \forall \xi\in\Ga(E), [D,\g(\xi)]\in\D_{k-1}(M,S)\}.
\]
Note that the ordinary filtration on $\D(M,S)$, given by the differential order,
%the order of differential operator,
can be defined in a similar fashion  replacing $\Ga(E)$ by $\iCinfty(M)$.

\begin{prop}\label{Prop:Dk}
Let $\nabla^S$ be a spinor connection on $S$. For all $k\in\bbN$, we have
\[\label{filtration}
  \D_{k}(M,S) = \Span\{\g(\eta)\circ\conn_{X_1}^S\cdots\conn_{X_m}^S\,|\,  2m+\deg(\eta)\leq k \},
\]
where $X_1,\ldots,X_m\in\Vect(M)$ and $\deg(\eta)$ denotes the degree of $\eta\in\cO^\bbC(E[1])$.
\end{prop}
\begin{proof}
The proof is by induction. Set
$$
\widetilde{\D_k}(M,S)=\Span\{\g(\eta)\circ\conn_{X_1}^S\cdots\conn_{X_m}^S\,|\,  2m+\deg(\eta)\leq k \}.
$$
By definition, $\D_{-1}(M,S)=\widetilde{\D_{-1}}(M,S)=\{0\}$.
Assume that $\D_k(M,S)=\widetilde{\D_k}(M,S)$ holds for a given $k\geq -1$.
Since $\nabla^S$ is a spinor connection, there exists a metric connection $\nabla^E$
on $E$ such that the following identities hold:
\[\label{CommRel:D}
 [\conn_X^S,\g(\xi)]=\g(\conn_X^E\xi) \quad \text{and}
\quad [\g(\eta),\g(\xi)]=\mg(\eta,\xi),  \qquad\forall \xi,\eta\in\Ga(E),
X\in\Vect(M).
\]
Hence, if $D\in\widetilde{\D_{k+1}}(M,S)$, we have
$[D,\g(\xi)]\in\widetilde{\D_{k}}(M,S)$
for all $\xi\in\Ga(E)$.
The inclusion $\widetilde{\D_{k+1}}(M,S)\subset \D_{k+1}(M,S)$ thus follows.
To prove the converse inclusion we need a lemma.

Let $(x^i)$ be a local coordinate system on $M$.
For any  multi-index $\a=(\a_1,\ldots,\a_k)$,
$\conn_\a^S$ stands for the composition
$\conn_{\a_1}^S\circ\cdots\circ\conn_{\a_k}^S$,
where $\conn_i^S:=\conn_{\frac{\partial}{\partial x^i}}^S$.

\begin{lem}\label{lem:DecompoD} %\JP{New Lemma}
Any operator $D\in\D(M,S)$ locally admits a unique linear decomposition
$$
D=\sum_{\k,\a} \g(\eta^\a_\k)\conn_\a^S,
$$
where $\eta^\a_\k\in\Ga(\wedge^\k E\otimes \bbC)$. In the sum above, $\k$ runs over
$\bbN$ and $\a$ runs over ordered multi-indices of arbitrary
 length $|\a|=k\in\bbN$,
i.e. $\a=(\a_1,\ldots,\a_k)$ and $1\leq \a_1\leq\cdots\leq\a_k\leq \dim M$.
\end{lem}
\begin{proof}
It is well-known that $\D(M,S)$ is a locally free $\iCinfty(M)$-module,
generated as an algebra by $\conn_i^S$ and $\g(\eta)$, with $\eta\in\Ga(\wedge^\k E\otimes \bbC)$, $\k\in\bbN$.
The result follows  from Eqns~\eqref{compatibilitySE}--\eqref{compatibilityR}.
\end{proof}

We assume now that $D\in\D_{k+1}(M,S)$.    %\JP{New proof}
There exists a minimum $l\in \bbN$  such that
$D\in\widetilde{\D_{l}}(M,S)$. We prove that $l-1\leq k$.
By Lemma \ref{lem:DecompoD}, we have
$D=\sum_{\k+2|\a|\leq l} \g(\eta^\a_\k)\conn_\a^S$.
From the definition of $\D_{k+1}(M,S)$, we deduce that
$[D,\g(\xi)]\in\widetilde{D_{k}}(M,S)$ and
$\g(\xi)\cdot [D,f]=[D,\g(\xi)f]-[D,\g(\xi)]f\in\widetilde{D_{k}}(M,S)$.
These relations read respectively as
\begin{align}\label{Eq:D1}
\sum_{\k+2|\a|\leq l}\,\big[\g(\eta^a_\k),\g(\xi)\big]\cdot\conn_\a^S+\g(\eta^\a_\k)\cdot\big[\conn_\a^S,\g(\xi)\big]&\in\widetilde{D_{k}}(M,S),\qquad\forall\xi\in\Ga(E),\\ \label{Eq:D2}
\sum_{\k+2|\a|\leq l}\, \g(\xi)\g(\eta^\a_\k)\cdot\big[\conn_\a^S,f\big]&\in\widetilde{D_{k}}(M,S),\qquad\forall\xi\in\Ga(E),f\in\iCinfty(M).
\end{align}
Consider $\k$ and $\a$ such that $\eta_\k^\a$
is non-vanishing, $\k+2|\a|=l$ and $\a$ is maximal.
The term $\big[\g(\eta^a_\k),\g(\xi)\big]\cdot\conn_\a^S$
is then clearly independent of the other terms in  \eqref{Eq:D1},
hence it pertains to $\widetilde{D_{k}}(M,S)$.
If $\k\neq 0$, there exists $\xi\in\Ga(E)$ such that
$0\neq \big[\g(\eta^a_\k),\g(\xi)\big]\in\Ga(\wedge^{\k-1} E)$,
and we deduce that $l-1\leq k$. If $\k=0$, the same conclusion
follows from ~\eqref{Eq:D2}.
\end{proof}

According to  Proposition \ref{Prop:Dk}, we have
 $\g(\xi)\in\D_1(M,S)$ and
$\conn_X^S\in\D_2(M,S)$,  $\forall \xi\in\Ga(E)$ and $X\in{\mathfrak X}(M)$.
Moreover, the operators in $\D_0(M,S)$ are multiplication
by functions in $\iCinfty(M)\otimes\bbC$.
The following proposition
is a direct consequence of ~\eqref{Groth}--\eqref{filtration}.

\begin{prop}\label{Prop:filtration} %\edz{modified proposition}
The subspaces $\{\D_k(M,S)\}_{k\in\bbN}$ define an increasing 
 and exhaustive filtration 
of the algebra $\D(M,S)$:
\begin{align*}
\D_0(M,S)\subset\D_1(M,S)&\subset\cdots\subset\D_k(M,S)\subset\cdots,\\
\D(M,S)&=\bigcup_{k\in\bbN}\D_k(M,S).
\end{align*}
Moreover, for any $k,l\in\bbN$, we have
\begin{align*}
&\D_k(M,S)\cdot\D_l(M,S) \subseteq \D_{k+l}(M,S),%\qquad\forall k,l\in\bbN.\\
&\left[\D_k(M,S),\D_l(M,S)\right] \subseteq \D_{k+l-2}(M,S).%\qquad\forall k,l\in\bbN.
\end{align*}
\end{prop}
\begin{rem}
\label{remFiltration}
Note that in \cite{Get86,Mic10a,Mic12}, the same choice of filtration on $\D(M,S)$ was made,
 whereas in \cite{Get83} and \cite{Vor99}, a different choice
was used: both Clifford generators $\g(\xi)$, $\xi\in\Ga(E)$, and covariant derivatives
$\conn_X^S$, $X\in{\mathfrak X}(M)$, are of order $1$.
In the more standard filtration, however,  the Clifford generators
normally have order $0$.
\end{rem}

As a consequence of Proposition \ref{Prop:filtration}, the graded algebra
$$
\Gr\D(M,S)=\bigoplus_{k\in\bbN}\D_{k}(M,S)/\D_{k-1}(M,S),
$$
is a graded commutative Poisson algebra of degree $-2$.
%Moreover, it is isomorphic to $\D(M,S)$ as $\iCinfty(M)$-module. %\edz{New sentence}
By $\varsigma$, we denote the canonical projection  $\varsigma:\D(M,S)\rightarrow\Gr\D(M,S)$.

\begin{thm}\label{thm:IsoPoissonAlg}
Let $S$ be a spinor bundle of $(E,\mg)$
and $(\nabla^E,\nabla^S)$ a pair of compatible connections.
There exists a unique isomorphism of graded commutative Poisson $\bbC$-algebras:
$$\fS:\Gr\D(M,S)\stackrel{\sim}{\longrightarrow}\cO^\bbC(T^*[2]M\oplus E[1]),$$
 satisfying
\[\label{Sigma012}
(\fS\circ\varsigma)(f)=f, \qquad (\fS\circ\varsigma)\left(\g(\xi)\right)= \xi,
\qquad (\fS\circ\varsigma)(\conn_X^S)=X,
\]
for all $f\in\iCinfty(M)$, $\xi\in\Ga(E)$ and $X\in{\mathfrak X}(M)$.
Moreover, the isomorphism $\fS$  only depends  on the connection
$\nabla^E$.
\end{thm}

\begin{proof}%\edz{New proof}
Both algebras $\cO^\bbC(T^*[2]M\oplus E[1])$ and $\Gr\D(M,S)$
are generated by $\Ga(E)$ and ${\mathfrak X}(M)$ over  $\iCinfty(M)$
(see Proposition \ref{Prop:Dk}).
Hence, a $\bbC$-algebra morphism $\fS$ satisfying  \eqref{Sigma012}
is an isomorphism and  must be unique if it exists.
As a result, it suffices to prove the existence of $\fS$ locally.
We use canonical coordinates $(x^i,p_i)$ on $T^*M$,
and write $\conn_i^S:=\conn_{\frac{\partial}{\partial x^i}}^S$.
We  define a map $\fS$ locally by setting
$$
\fS\circ\varsigma\left(\g(\eta)\conn_{\a}^S\right)=\eta\, p_{\a},
$$
where $\eta\in\Ga(\wedge E\otimes \bbC)$,
$\a=(\a_1,\ldots,\a_k)$ is an ordered multi-index %(i.e.\ $\a_1\leq\cdots\leq\a_k$),
%$\conn_\a^\sS=\conn_{\a_1}^\sS\circ\cdots\circ\conn_{\a_k}^\sS$
and $p_{\a}=p_{\a_1}\cdots p_{\a_k}$.
According to Lemma \ref{lem:DecompoD}, the map $\fS$ is well-defined.
Moreover, by  \eqref{CommRel:D}, this is indeed a $\bbC$-algebra morphism
which satisfies  \eqref{Sigma012}.

To prove that $\fS$ is a Poisson map, it suffices to check
that $\fS$ preserves the Poisson brackets  on generators, i.e.
 on elements in  $\smooth(M)\oplus\Ga(E)\oplus{\mathfrak X}(M)$.
For any $f, g\in \smooth(M)$, $\xi, \eta
\in \Ga(E)$, and $X, Y
\in {\mathfrak X}(M)$, we have the following graded commutators
of differential operators
\begin{align}\nonumber
  &[\conn_X^S,f]=X(f),&                   & [g,f] =0,\\  \label{Commutator}
  &[\conn_X^S,\g(\xi)]=\g(\conn_X^E\xi),& &[g,\g(\xi)] =0,\\ \nonumber
  &[\conn_X^S,\conn_Y^S]=\conn_{[X,Y]}^S+R^S(X,Y),&&  [\g(\xi),\g(\eta)] =\mg(\xi,\eta),
\end{align}
where we have used Eq.\ \eqref{compatibilitySE}.
Writing the curvature $R^S$ as in Eq.\ \eqref{compatibilityR} and
applying the map  $\fS\circ\varsigma$,
we  immediately   see that  these commutators  coincide with
the generating relations of the Poisson bracket as
 given in \eqref{PoissonRothstein}. Therefore,
$\fS$ is indeed a Poisson map.

According to Lemma~\ref{lem:CompatibleConnections},
two spinor connections, compatible with the same metric connection $\nabla^E$,
satisfy $\widetilde{\nabla}^S-\nabla^S=\nu\in\Om^1(M)\otimes\bbC$.
However,  different $\nu$ do not
 modify the defining relations for $\fS$ in  \eqref{Sigma012}.
This concludes the proof of the theorem.
\end{proof}

\begin{defn}\label{def:symbols}
We define the principal symbol map
$\sigma_k: \D_k(M,S)\to \cO^\bbC_k(T^*[2]M\oplus E[1])$, $\forall k\in\bbN$,
as the composition  $\fS\circ\varsigma$.
%Thus we have   the following commutative diagram:
%\[\label{Diagram:sigma}
%\xymatrix{
 % \D_k(M,S)\ar[rd]_{\sigma_k}\ar[r]^{\varsigma\quad\qquad} & \D_k(M,S)/\D_{k-1}(M,S) \ar[d]^{\fS}  \\
%& \cO^\bbC_k(T^*[2]M\oplus E[1]),
%}
%\]
%where $\varsigma$ is the canonical projection and $\fS$ is as in
%Theorem \ref{thm:IsoPoissonAlg}.
\end{defn}

Note that the principal symbol  map $\sigma_k$ depends
 on the choice of the connection $\nabla^E$.
%Accordingly, the elements in $\Gr\,\D(M,S)$ or $\cO^\bbC(T^*[2]M\oplus E[1])$
\begin{prop}\label{Thm:Symbol}
Let $k,l\in\bbN$. The principal symbol maps satisfy the following properties:
\begin{align}\nonumber% \label{SymbolProduct}
\sigma_{k+l}(AB) &= \sigma_k(A)\sigma_l(B),\\
\label{SymbolPoisson}
\sigma_{k+l-2}([A,B]) &= \{\sigma_k(A),\sigma_l(B)\},
\end{align}
for all $A\in\D_k(M,S)$ and $B\in\D_l(M,S)$.
\end{prop}
\begin{proof}
This is a direct consequence of Theorem~\ref{thm:IsoPoissonAlg}.
\end{proof}

\begin{rem}
It is helpful to compare our  filtration with those in
 \cite{Get83, Vor99}, as shown in the  table below:\par\medskip
\centerline{\begin{tabular}{|c||c|c|c|}
			\hline
Filtration 	& order  $\g(\xi)$  & order $\conn_X^S$  & $\Gr\D(M,S)$  \\[3pt]\hline\hline
 standard filtration	& 0                 & 2                & $\cO(T^*[2]M)\otimes_{\Cinfty}\Ga(\End S)$  \\[3pt]\hline
filtration in \cite{Get83,Vor99}       & 2                 & 2      & $\left(\Ga(\bigotimes TM) \otimes_{\Cinfty}\cO^\bbC(E[1])\right)/\mathcal{J}$  \\[3pt] \hline
filtration in \eqref{Groth} & 1  & 2         & $\cO^\bbC(T^*[2]M\oplus E[1])$ \\[3pt]\hline
\end{tabular}\\
\medskip
}
\noindent where $\xi\in\Ga(E)$, $X,Y\in{\mathfrak X}(M)$ and $\mathcal{J}=\left(X\otimes Y-
Y\otimes X - R^E(X,Y) \right)$. As before, we have $R^E(X,Y)\in\Ga(\wedge^2 E)$.
The  filtration \eqref{Groth} is the only one for which $\Gr\D(M,S)$ is
 a graded commutative algebra.
\end{rem}

Recall that, according to Theorem \ref{thmRoyM}, there exists
a symplectic diffeomorphism $\Xi_\nabla:\M\to T^*[2]M\oplus E[1]$.
Let $\rho_\nabla=(\Xi_\nabla)^*\circ\fS$.
We  have  the following

\begin{prop}\label{prop: cM-grD}
The map $ \rho_\nabla: \ \Gr\D(M,S)\rightarrow\gsmooth^\bbC(\cM)$
 is an isomorphism of  graded Poisson algebras,
which is   independent of the metric connection $\nabla^E$ used in defining
~$\Xi_\nabla$ and~$\fS$.
\end{prop}
\begin{proof}

Let $(\nabla^E, \nabla^S)$ be  a pair of   compatible connections,
 and $\Xi_\nabla$, $\fS$ their associated maps as defined in Theorems %\edz{Ref. added}
\ref{thmRoyM} and \ref{thm:IsoPoissonAlg}.
%Let $\rho_\nabla=(\Xi_\nabla)^*\circ\fS$.
 Since both
 $(\Xi_\nabla)^*$ and $\fS$  are isomorphisms of graded Poisson algebras,
so is  their composition $\rho_\nabla$.
It remains to prove that $\rho_\nabla$ is independent of the
choice of the connections $(\nabla^E, \nabla^S)$.
Since $\rho_\nabla$ is an algebra isomorphism, it suffices to prove
the assertion  on generators of $\Gr\D(M,S)$, which
are of the following three types:
 $f$, $\varsigma\circ\g(\xi)$ and $\varsigma(\conn_X^S)$,  $\forall f
\in \iCinfty(M)$, $\xi\in \Ga(E)$, $X\in{\mathfrak X}(M)$.
By definition,  both maps  $(\Xi_\nabla)^*$ and $\fS$
 are independent of the choice of  connections when
applying  to elements $f$ and $\varsigma\circ\g(\xi)$.
Thus it remains to show that $\rho_{\nabla}$ is
also independent of the choice of connections on the elements
 $\varsigma(\conn_X^S)$,  $\forall X\in{\mathfrak X}(M)$.

Assume that $(\widetilde{\nabla}^E, \widetilde{\nabla}^S)$ is
 another pair of compatible connections.
 From \eqref{eq:xi} and \eqref{Sigma012}, it follows that
\[ \label{rho-varsigma}%\begin{align}\nonumber
\begin{split}
 \rho_\nabla\circ\varsigma(\conn_X^S) &= (\pi_\nabla\circ i_{\M})^* l_X, \\
  \rho_{\widetilde{\nabla}}\circ\varsigma(\tildeconn_X^S)
    &= (\pi_{\widetilde\nabla}\circ i_{\M})^*l_X,
\end{split}
\]%\end{align}
where $l_X$ denotes the fiberwise linear function on $T^*[2]M$
 corresponding to $X\in{\mathfrak X}(M)$.
According to Lemma~\ref{lem:CompatibleConnections},
 there exists  $\varpi\in\Om^1(M,\wedge^2 E)$
 and $\nu\in\Om^1(M)\otimes\bbC$ such that
\begin{align}\label{Eq:nabla-tilde}
\tildeconn_X^S-\conn_X^S&=\g(\varpi(X))+\nu(X), \\ \label{eq:nabla-tilde1}
 \tildeconn_X^E-\conn_X^E&=\{\varpi(X),\cdot\}.
\end{align}
 Eq.\ \eqref{eq:nabla-tilde1} implies that
\[\label{Eq:FX}
(\pi_{\widetilde\nabla}\circ i_{\M})^*l_X-(\pi_\nabla\circ i_{\M})^* l_X=\pi^*\varpi(X).
\]
Therefore, by Eqns \eqref{rho-varsigma} and \eqref{Eq:FX}, we have
$$
\rho_{\widetilde{\nabla}}\circ\varsigma(\tildeconn_X^S) =  \rho_\nabla\circ\varsigma(\conn_X^S)+\pi^*\varpi(X).
$$
On the other hand,  Eq.\ \eqref{Eq:nabla-tilde} implies that
$$
\rho_{\nabla}\circ\varsigma\big(\tildeconn_X^S) = \rho_{\nabla}\circ\varsigma\big(\conn_X^S\big)+\pi^*\varpi(X).
$$
As a consequence,  we conclude that
the equality $\rho_{\widetilde{\nabla}}=\rho_\nabla$ holds on every element
 $\varsigma(\tildeconn_X^S)$, $\forall X\in{\mathfrak X}(M)$.
This concludes the proof of the proposition.
\end{proof}

%%% 3.3 %%%%%%%%%%%%%%%%%%%%%%%%%%%%%%%%%%%%%%%%%%%%%%%%%%%%%%%%
\subsection{Two involutions on $\D(M,\sS)$}\label{ParReal}
%%%%%%%%%%%%%%%%%%%%%%%%%%%%%%%%%%%%%%%%%%%%%%%%%%%%%%%%%%%%%%%%

 Assume that the line bundle $(\det S^*)^{1/N}$  exists,
with $N$ being the rank of $S$.
Consider the \emph{twisted spinor bundle}
 $\sS:=S\otimes(\det S^*)^{1/N}\otimes\Vol^{1/2}$, where
 $\Vol^{1/2}$ denotes the  half-density   line bundle of $M$.
A connection  on $TM$ induces a connection on 
$\Vol^{1/2}$.
%all density bundles.
Together with a spinor connection $\nabla^S$ on $S$,
it yields an \emph{induced spinor connection} on the twisted spinor bundle
$\sS$, denoted by $\nabla^\sS$.
Since $\sS$ is a spinor bundle,
the algebra $\D(M,\sS)$ of differential operators  on $\sS$,
satisfies in particular all the results proved in Section \ref{ParSymbol}
for the algebra of spinor differential operators.
Considering $\sS$ rather than an arbitrary spinor bundle
 allows us to have  additional structures on $\D(M,\sS)$.

First, we introduce an adjoint operation on the algebra $\D(M,\sS)$.
Let $U\subset M$ be a contractible open subset.
A pseudo-Hermitian pairing on $\sS|_U$ is a fiberwise non-degenerate sesquilinear map
\[\label{PairingS}
\<\cdot,\cdot\>_U:\Ga(\sS|_U)\times\Ga(\sS|_U)\longrightarrow \Ga(|{\wedge^{\mathrm{top}}T^*U}|)\otimes\bbC
\]
 satisfying $\overline{\<\phi,\psi\>_U}=\<\psi,\phi\>_U$ for all $\phi,\psi\in\Ga(\sS|_U)$.

\begin{prop}\label{Prop:pairing_S}
Up to multiplication by a scalar $a\in\bbR^\times$, there exists
 a unique smooth pseudo-Hermitian pairing $\<\cdot,\cdot\>_U$ on $\sS|_U$
satisfying the following properties:
\begin{enumerate}[label=(\roman*)]
\item \label{1:nabla} $\conn_X(\<\phi,\psi\>_U)=\<\conn_X^\sS\phi,\psi\>_U+\<\phi,\conn_X^\sS\psi\>_U,$
\item \label{2:symmetric} $\langle \g(\xi)\phi,\psi\rangle_U=\langle \phi,\g(\xi)\psi\rangle_U,$
\end{enumerate}
for any $ X\in\Ga(TU)$, $\xi\in\Ga(E|_U)$
and $\phi,\psi\in\Ga\left( \sS|_{U} \right)$.

Moreover, $\<\cdot,\cdot\>_U$   is independent of
 the choice of connections $\nabla$ and $\nabla^\sS$
involved in \ref{1:nabla}.
%It is called a \emph{spinor pairing}.
\end{prop}
\begin{proof}
Choosing a local frame $(\xi^i)$ of  $E|_U$  such that
$\mg (\xi^i, \xi^j)=\pm \delta_{ij}$,
we obtain  a trivialization of the bundle $E|_U\cong U\times \bbR^n$,
  on which $\mg$ becomes  a constant  metric. This in turn  induces
a trivialization of $\bbCl (E|_U)\cong  \End (S|_U)$, for which $\g(\xi^i) \in
\Gamma (\End (S|_U ))\cong C^\infty (U, \End   \bbC^N)$
are  constant sections, for all  $i=1, \cdots n$.
Let  $S|_U\cong U\times \bbC^N$ be a trivialization defined by
  a local frame $(\phi_a)$ of $S$. The action of $\End (S|_U)$ on
$S|_U$ translates into a family of faithful representations 
$\rho_x:\End   \bbC^N \times \bbC^N \longrightarrow \bbC^N$, depending smoothly on $x\in U$.
Each representation defines an automorphism of $\End   \bbC^N$.
Since  automorphisms of $\End   \bbC^N$ are all inner automorphisms,
 there exists a unique element $\alpha_x\in \text{PGL}(N, \bbC)$
such that $\rho_x=\text{Ad}_{\alpha_x}$. It is clear that $x\to \alpha_x$ is
a smooth map from $U$ to  $\text{PGL}(N, \bbC)$. Shrinking $U$ if necessary,
we can choose 	 a smooth lift	 $x\to A_x$  from $U$ to 
 $\text{GL}(N, \bbC)$.
Then $\rho_x(M)\cdot \phi_0=A_x M A_x^{-1}\cdot \phi_0$, for all $x\in U$, $M\in\End   \bbC^N$ and $\phi_0\in\bbC^N$.
Thus  $(A_x\phi_a)$ is another local frame of $S$, under which
 the trivialization of $S|_U$  is consistent with the
trivialization of $\End (S|_U )$ induced by the 
trivialization of $\bbCl (E|_U)$ at the  beginning of the proof.
% is coherent with the one chosen in $\End (S|_U )$ : constant endomorphism over $U$ have constant actions.
 We pick up such a  trivialization of $S|_U$ for the rest of
 the proof, so that the operators $\g(\xi^i) \in
\Gamma (\End (S|_U ))\cong C^\infty (U, \End   \bbC^N)$
acts as constant endomorphisms on $S|_U\cong U\times \bbC^N$, for all  $i=1, \cdots n$.
%By choosing a nowhere vanishing section $\phi\in\Gamma(S|_U)$,
%then we obtain a trivialization  $S|_U\cong U\times \bbC^N$, under
%which the operators $\g(\xi^i) \in
%\Gamma (\End (S|_U ))\cong C^\infty (U, \End   \bbC^N)$
%are  constant endomorphisms, for all  $i=1, \cdots n$.
Consider the induced trivialization of the line bundle
$(\det S^*|_U)^{1/N}$ and pick up a trivialization
of $|\wedge^{\mathrm{top}}T^*U|$.
We can identify   the space of sections of
$|\wedge^{\mathrm{top}}T^*U|$ with
$C^\infty (U)$,  and the  space of sections of
 $\sS$   with $C^\infty (U, \bbC^N)$.
Under such identifications,
the connections $\nabla$ and $\nabla^\sS$ can be written as follows:
\begin{align}\nonumber
\nabla_X v&= X v+\nu_0(X)v, \qquad \forall v\in\Ga(|\wedge^{\mathrm{top}}T^*U|)
\cong  C^\infty (U) ,\\ \label{Eq:nablasS}
\conn_X^\sS\phi&=X\phi+\g(\varpi(X))\phi+\frac{1}{2}\nu_0(X)\phi, \qquad \forall \phi\in\Ga(\sS|_U)\cong C^\infty (U, \bbC^N),
\end{align}
where $\nu_0\in\Om^1(U)$ and $\varpi\in\Om^1(U,\wedge^2E|_U)$
(see  Lemma \ref{lem:CompatibleConnections}).
%Under  such a  trivialization of $\sS_{|U}$,
If  Property  \ref{2:symmetric} holds,
it is then simple to check that  Property \ref{1:nabla} is  equivalent to
\begin{equation}
\label{eq:constant}
X(\<\phi,\psi\>_U)=\<X\phi,\psi\>_U+\<\phi,X\psi\>_U, \quad\forall \phi, \psi
\in C^\infty (U, \bbC^N), X\in \Gamma (TU).
\end{equation}
The latter equation means  that $\<\cdot,\cdot\>_U$ is a constant pairing
in the trivialization $\sS|_U\cong U\times \bbC^N$.

%Consider $\<\cdot,\cdot\>_U\vert_x$ the restriction of a pseudo-Hermitian $\<\cdot,\cdot\>_U$
%over a point $x\in U$.
Let $x\in U$ and $\beta$ be the semi-linear involution of the complex Clifford algebra
$\bbCl (E_x)\cong  \End \sS_x$ defined by
$$
\beta\left(\g(\xi^{i_1})\g(\xi^{i_2})\cdots \g(\xi^{i_\kappa})\right) = \g(\xi^{i_\kappa})\cdots \g(\xi^{i_2})\g(\xi^{i_1}),
\qquad \forall i_1, i_2,\ldots, i_\kappa \in \{1,2,\ldots,n\}.
$$
Then, $\<\cdot,\cdot\>_U$ satisfies Property \ref{2:symmetric} over $x$ if and only if
its adjoint operation on $\End \sS_x\cong\bbCl (E_x)$ is given by $\beta$.
The existence of a pseudo-Hermitian pairing on the complex spinor space $\sS_x$
with adjoint operation $\beta$ is classical.
 See e.g.\ \cite[Theorem 7.14]{Por95}.
Two such pairings differ by multiplication by a non-vanishing scalar,
because they have the same adjoint operation.
Using the trivialization $\sS|_U\cong U\times \sS_x$, we deduce that there exists
a constant pseudo-Hermitian pairing $\<\cdot,\cdot\>_U$ satisfying
Property \ref{2:symmetric} over $U$, and that any other such pairing is of the form
$f\,\<\cdot,\cdot\>_U$ with $f\in\iCinfty(U)$. According to Eq.\ \eqref{eq:constant},
the pairing $f\,\<\cdot,\cdot\>_U$ satisfies  Property \ref{1:nabla} if and only if
$f$ is a constant function.
%According to \cite{Por69},
%there exists a smooth fiberwise  pairing on $\sS|_{U}$ satisfying Property
% \ref{2:symmetric}, and  any two such pairings
%differ by multiplication of a smooth function $f\in C^\infty (U)$.
%By Eq.~\eqref{eq:constant}, this function $f$
%must be a constant function.

As Property \ref{2:symmetric} and Eq.\ \eqref{eq:constant} are independent
of the choice of connections $\nabla$ and $\nabla^\sS$, so
is the  pseudo-Hermitian pairing $\<\cdot,\cdot\>_U$.
 This concludes the proof.
\end{proof}

The pairing $\<\cdot,\cdot\>_U$ in Proposition \ref{Prop:pairing_S} 
 is called a \emph{spinor pairing}.
Since one can integrate $1$-densities,
a spinor pairing yields a pseudo-Hermitian scalar product on the space
of compactly supported sections~$\Ga_c(\sS|_U)$:
\[\label{HermitianProduct}
(\phi,\psi)_U=\int_U \<\phi,\psi\>_U, \qquad  \forall\phi,\psi\in\Ga_c(\sS|_U).
\]
This is called a \emph{spinor scalar product}. %Two such scalar products
%two of them must differ by multiplication of a constant $a\in\bbR^\times$.
It is not clear if $(\cdot,\cdot)_U$ can be extended to a  global
scalar product on $M$.
However, as we see below, the spinor scalar products over $\sS_{|U}$
enable us to introduce  a globally defined  adjoint operation
on $\D(M,\sS)$.

\begin{lem}\label{lem:*exist}
For any  $D\in\D(M,\sS)$, there exists a unique differential operator
 $D^*\in\D(M,\sS)$
such that
\[\label{adjoint}
(D\phi,\psi)_U=(\phi,D^*\psi)_U, \qquad \forall\phi,\psi\in\Ga_c(\sS|_U),
\]
where $U$ is any contractible open subset of $M$ and
$(\cdot,\cdot)_U$ is any spinor scalar product on $\sS_{|U}$.
\end{lem}

\begin{proof}
Let $D\in\D(M,\sS)$.
Denote by $U\subset M$ a contractible open subset
and by $D|_U\in\D(U,\sS|_U)$ the restriction of the operator $D$ to $U$.
Choose a spinor scalar product  $(\cdot,\cdot)_U$.
% as in Eq.\ \eqref{HermitianProduct}.
As is well-known, there exists a unique operator $(D|_U)^*\in\D(U,\sS|_U)$
 satisfying  \eqref{adjoint}. Clearly, this operator also
satisfies   \eqref{adjoint}
for the spinor scalar product $a(\cdot,\cdot)_U$, with $a\in\bbR^\times$.
By Proposition \ref{Prop:pairing_S},
$(D|_U)^*$ satisfies   \eqref{adjoint} for any choice of
  spinor scalar products on~$\sS|_{U}$.
% as in \eqref{HermitianProduct}.

Since the operators $(D|_U)^*$ are uniquely defined
by   \eqref{adjoint},
they must glue  together into
a globally defined differential operator $D^*\in\D(M,\sS)$.
\end{proof}

The map $D\mapsto D^*$ is called the \emph{adjoint operation},
and admits the usual  properties of the standard adjoint operation.
%which are defined from a single scalar product.

\begin{prop}\label{prop:*}
The adjoint operation $D\mapsto D^*$ satisfies the following properties:
\begin{enumerate}[label=(\roman*)]
\item it is  an involutive map on $\D(M,\sS)$ preserving the filtration;
\item it is a $\bbC$-antilinear antiautomorphism:
\[\label{*anti}
\begin{cases}
(\l_1 D_1+ \l_2 D_2)^*=\overline{\l_1}D_1^*+\overline{\l_2}D_2^* ,\\
(D_1D_2)^*=D_2^*D_1^*,
\end{cases}
\qquad \forall \l_1, \l_2 \in
\bbC,  D_1,D_2\in\D(M,\sS),
\]
uniquely determined by the following relations
\[\label{Eq:D*}
f^*=\overline{f},\qquad \g(\xi)^*=\g(\xi)\quad\text{and}\quad
(\conn_{X}^\sS)^*=-\conn_X^\sS-\Tr\nabla X,
\]
for all $f\in\iCinfty(M)\otimes \bbC$, $\xi\in\Ga(E)$ and
$X\in{\mathfrak X}(M)$. Here, $\Tr\nabla X\in\iCinfty(M)$ is the trace
of $\nabla X\in\Ga(T^*M\otimes TM)$; %\cong \Ga (\End TM)
\item it  does not depend on the choice of  connections
$\nabla$ on $TM$ and $\nabla^\sS$ on $\sS$.
\end{enumerate}
\end{prop}
\begin{proof}
It suffices to prove all the  properties above over
any contractible open subset $U\subset M$.

Choose  a spinor scalar product $(\cdot,\cdot)_U$. % as in \eqref{HermitianProduct}.
Eq. \eqref{adjoint}  uniquely characterizes the operation $*$.
As a consequence, Eq.\  \eqref{*anti} holds as well as
the relations  $(D^*)^*=D$, $f^*=\overline{f}$.
%From Eq.\ \eqref{adjoint}, we deduce Eq.\  \eqref{*anti} and the relations  $(D^*)^*=D$, $f^*=\overline{f}$.
Using in addition Proposition \ref{Prop:pairing_S}, we obtain that
$\g(\xi)^*=\g(\xi)$ and that 
the  operation  $*$ is independent of  any choice of connections.

Any $v\in\Ga_c (\Vol)$ satisfies the identities $L_X v= \nabla_X v+
 (\Tr\nabla X)v$ and
$\int_U L_X v=0$. For any $\phi$ and $\psi\in\Ga_c(\sS|_U)$, their pairing $\<\phi,\psi\>_U$
pertains to $\Ga_c (\Vol)$, so the above identities hold
for $v=\<\phi,\psi\>_U$. Using in addition Proposition \ref{Prop:pairing_S},
we obtain that
$$
\int_U L_X \<\phi,\psi\>_U = \int_U \<\conn_X^\sS\phi,\psi\>_U + \<\phi,\conn_X^\sS\psi+(\Tr\nabla X)\psi\>_U =0.
$$
Hence, the equation $(\conn_{X}^\sS)^*=-\conn_X^\sS-\Tr\nabla X$ holds.

Since $\D(M,\sS|_U)$ is generated by $\iCinfty(U)$, $\Ga(E|_U)$ and $\Ga(TU)$,
Eqns \eqref{*anti}-\eqref{Eq:D*} completely determine
the adjoint operation~$*$. Finally, one  easily
sees that  if $D\in\D_k(M,\sS|_{U})$ then $D^*\in\D_k(M,\sS|_{U})$.
This concludes the proof of the proposition.
\end{proof}

Denote by $\overline{\, \cdot \,}$ the complex conjugation on $\bbC$
and its natural extension to any complexified $\bbR$-vector space
 $V\otimes\bbC$.
In what follows, we  will extend this operation to the algebra $\D(M,\sS)$.

\begin{prop}\label{Prop:Conjugation}
There exists a unique $\bbC$-antilinear algebra morphism
$$\overline{\, \cdot \,}:\D(M,\sS)\longrightarrow\D(M,\sS),$$
 which coincides
with the complex conjugation on $\D_0(M,\sS)\cong\iCinfty(M)\otimes\bbC$
and  satisfies the following properties:
\[\label{Eq:Conjugation}
\overline{\g(\xi)}=\g(\overline{\xi})
\quad\text{and}\quad\overline{\conn_X^\sS}=\conn_{\overline{X}}^\sS,
\qquad \forall \xi\in\Ga(E)\otimes\bbC, X\in{\mathfrak X}(M)\otimes\bbC,
\]
for any induced spinor connection $\nabla^\sS$.
%, induced by connections on $TM$ and $S$.
%
%Moreover, the map $\overline{\, \cdot \,}$ is an involution, and
%is independent of the choice of  connections on
% $TM$ and $S$ inducing $\nabla^\sS$.
\end{prop}

\begin{proof}
Choose a spinor connection on  $S$ and a connection on $TM$.
Let  $\nabla^\sS$ be the  induced  connection on
$\sS$.
Since $\D(M,\sS)$ is generated by $f$, $\g(\xi)$, $\conn_X^\sS$,
with $f\in\iCinfty(M)\otimes\bbC$, $\xi\in\Ga(E)$, $X\in{\mathfrak X}(M)$,
 a $\bbC$-antilinear algebra morphism $\overline{\,\cdot\,}$
satisfying  \eqref{Eq:Conjugation} must be unique
if it exists. Hence, it suffices to prove the existence of
 $\overline{\,\cdot\,}$ locally.
Over a contractible open subset of $M$,  we extend the conjugation
of complex functions  by setting
$$
\overline{f\g(\eta)\conn_{\a}^\sS}:=
\overline{f}\g(\eta)\conn_{\a}^\sS,
$$
with $f\in\iCinfty(M)\otimes \C$, $\eta\in\Ga(\wedge E)$ and
$\a=(\a_1,\ldots,\a_k)$ an ordered multi-index.
%(i.e.\ $\a_1\leq\cdots\leq\a_k$) and $\conn_\a^\sS=\conn_{\a_1}^\sS\circ\cdots\circ\conn_{\a_k}^\sS$.
By Lemma~\ref{lem:DecompoD}, the map $\overline{\,\cdot\,}$ is well-defined.
It is obvious to see that
$\overline{\,\cdot\,}$ is indeed a $\C$-antilinear involutive map
satisfying  \eqref{Eq:Conjugation}.
To check that $\overline{\,\cdot\,}$ is   an algebra morphism,
 one may use the commutation relations \eqref{Commutator}, and note
that the curvature $2$-form of $\nabla^\sS$ is valued in the real
 Clifford bundle $\Cl(E)$ (see  \eqref{Eq:nablasS}).

According to  \eqref{Eq:nablasS},
different choices of connections on $TM$ and $S$
induce connections on $\sS$ which differ by a real term.
 Therefore  the map $\overline{\, \cdot \,}$ is
independent of choices of connections. This concludes the proof.
\end{proof}

The map $\overline{\, \cdot \,}$ introduced above is
called the {\em conjugation map}.
By $\D(M,\sS)_\bbR$, we denote the real subalgebra of $\D(M,\sS)$
generated by sections of the real Clifford bundle
 $\Cl(E)$ and the covariant derivatives $\conn_X^\sS$, with
$X\in {\mathfrak X}(M)$.  The next proposition gives an intrinsic
description of this  real subalgebra.

\begin{prop}\label{Dreal}
The real subalgebra $\D(M,\sS)_\bbR$ is the fixed point set of $\overline{\, \cdot \,}$.
Therefore we have
 $\D(M,\sS)\cong\D(M,\sS)_\bbR\otimes\bbC$.
\end{prop}

\begin{proof}
The proof is straightforward and is left to the reader.
%Since the conjugation map $\overline{\,\cdot\,}$ is an algebra morphism,
%its fixed point set is a subalgebra of $\D(M,\sS)$. It clearly contains $\sD(M,\sS)$.
%Since $\bbC$ and $\sD(M,\sS)$ generate $\D(M,\sS)$, we deduce  that
%$\sD(M,\sS)$ is exactly the fixed point set. %of the map $\overline{\,\cdot\,}$.
%Hence, it satisfies $\D(M,\sS)\cong\sD(M,\sS)\otimes_\bbR\bbC$.
%As the conjugation map does not depend on the chosen connections
%on $TM$ and $S$, the same holds for $\sD(M,\sS)$.
\end{proof}

\begin{rem}\label{rem:realspinor}
Assume there exists a real spinor bundle $\sS_\bbR$ such that
$\sS\cong\sS_\bbR\otimes\bbC$. % and $\End\sS_\bbR\cong\Cl(E)$.
Then, the conjugation map on $\D(M,\sS_\bbR\otimes\bbC)$
is induced by the natural conjugation on $\sS_\bbR\otimes\bbC$,
so that $\D(M,\sS)_\bbR\cong\D(M,\sS_\bbR)$.
However, $\sS_\bbR$ may not exist.
\end{rem}

The following result is obvious.

\begin{prop}
For any $k,l\in\bbN$, let $m=2k+l$.   Assume that the principal
symbol of $D\in\D_m(M,\sS)$
satisfies $\sigma_m(D)\in\cO_k(T^*[2]M)\otimes_{\Cinfty}\cO^\bbC_l(E[1])$.
Then we have
$$
\sigma_m(\overline{D})=\overline{\sigma_m(D)}\quad \text{ and } \quad
\sigma_m(D^*)=(-1)^k(-1)^{\frac{l(l-1)}{2}}(\overline{\sigma_m(D)}).
$$
\end{prop}

%% 4 %%%%%%%%%%%%%%%%%%%%%%%%%%%%%%%%%%%%%%%%%%%%%%%%%%%
%%%%%%%%%%%%%%%%%%%%%%%%%%%%%%%%%%%%%%%%%%%%%%%%%%%%%%%%
\section{Quantization of symplectic graded manifolds of degree $2$}
%%%%%%%%%%%%%%%%%%%%%%%%%%%%%%%%%%%%%%%%%%%%%%%%%%%%%%%%
%%%%%%%%%%%%%%%%%%%%%%%%%%%%%%%%%%%%%%%%%%%%%%%%%%%%%%%%

First, we recall some  basic  materials on the  Weyl quantization
on $T^*\bbR^n$ and its extension to  arbitrary cotangent bundles $T^*M$
 \cite{Und78,Wid80}. Then, we introduce a similar construction
of  quantizations on  symplectic graded manifolds of degree $2$.

%%% 4.1 %%%%%%%%%%%%%%%%%%%%%%%%%%%%%%%%%%%%%%%%%%%%%%%%%%%%%%%%
\subsection{Weyl quantization on $T^*\bbR^n$}
%%%%%%%%%%%%%%%%%%%%%%%%%%%%%%%%%%%%%%%%%%%%%%%%%%%%%%%%%%%%%%%%

 Let $V\rightarrow\bbR^n$ be a complex vector bundle over
$\bbR^n$, and denote by
  $\D(\bbR^n,V)$ the algebra of differential operators on $V$.
 Choosing a trivialization of~$V$, one readily gets that
$\D(\bbR^n,V)\cong\D(\bbR^n)\otimes_{\Cinfty}\Ga(\End V)$.
Recall that the normal order quantization is defined,
in terms  of the canonical coordinate system  $(x^i,p_i)$ on $T^*[2]\bbR^n$, by
\begin{align*}
\cN_\hbar :\cO(T^*[2]\bbR^n)\otimes_{\Cinfty}\Ga(\End V)&\longrightarrow\D(\bbR^n,V)\\
F^{i_1\cdots i_k}(x)p_{i_1}\cdots p_{i_k}&\longmapsto
F^{i_1\cdots i_k}(x)\big( \frac{\hbar}{\bi}\frac{\partial}{\partial x^{i_1}}\big)
\cdots \big( \frac{\hbar}{\bi}\frac{\partial}{\partial x^{i_k}}\big),
\end{align*}
where  $F^{i_1\cdots i_k}\in\Ga(\End V)$
and  $\hbar\in\bbC^\times$  is a parameter\footnote{In this paper,
we use the notation $\hbar$ to denote a nonzero  variable in $\bbC$.
In most situation, $\hbar$ can be interpreted as
the Planck constant. However, we sometimes let $\hbar=\bi$.}.
%Clearly, the map $\cN_\hbar$ is a $\bbC$-linear isomorphism.
Let $\Div=\frac{\partial^2}{\partial x^{i}\partial p_i}$ be the
divergence-like operator acting on the symbol algebra $\cO(T^*[2]\bbR^n)$.

\begin{defn}
On $T^*[2]\bbR^n$, the $\D(\bbR^n,V)$-valued
Weyl quantization is  a $\bbC$-linear isomorphism,
$\cQ^{\bbR^n}_\hbar :\cO(T^*[2]\bbR^n)\otimes_{\Cinfty}\Ga(\End V)\to\D(\bbR^n,V)$,
indexed by $\hbar\in\bbC^\times$ and defined by
\[\label{WQ-ND}
\cQ^{\bbR^n}_\hbar := \cN_\hbar\circ\exp\left(\tfrac{\hbar}{2\bi}\Div\otimes\Id\right).
\]
\end{defn}
In particular, the Weyl quantization satisfies
$\cQ^{\bbR^n}_\hbar\!(p_i)=\frac{\hbar}{\bi}\frac{\partial}{\partial x^i}$
and $\cQ^{\bbR^n}_\hbar\!(F)= F$, for any $F\in\Ga(\End V)$.
In the classical case when $V$ is  a trivial  line bundle,
 $\Ga(\End V)$  is
identified with $\iCinfty(\bbR^n)$, and $\cQ^{\bbR^n}_\hbar\!(F)$
is the multiplication operator by $F$.
For polynomials in the coordinates~$(x^i,p_i)$,
the Weyl quantization is just the symmetrization map
satisfying $\cQ^{\bbR^n}_\hbar\!(p_i)=\frac{\hbar}{\bi}\frac{\partial}{\partial x^i}$
and $\cQ^{\bbR^n}_\hbar\!(x^i)=x^i$.

The Weyl quantization can also be defined by the integral formula~\eqref{WQ-Rn}
below  \cite{Fol89}.
\begin{prop}
The Weyl quantization satisfies
\[
\label{WQ-Rn}
(\cQ^{\bbR^n}_\hbar\!(F)\,\psi)(x) = \tfrac1{(2\pi\hbar)^n}\int_{T^*_x\R^n\oplus T_x\bbR^n}\e^{-\frac{\bi}{\hbar}\<p,v\>}\; F(x+v/2,p)\cdot\psi(x+v)\;\ud p\,\ud v,
\]
for all $\psi\in\Ga(V)$ and $F\in\cO(T^*[2]\bbR^n)\otimes_\Cinfty\Ga(\End V)$.
\end{prop}

%%% 4.2 %%%%%%%%%%%%%%%%%%%%%%%%%%%%%%%%%%%%%%% FINDDD
\subsection{Exponential map and parallel transport}\label{Sec:Exp}
%%%%%%%%%%%%%%%%%%%%%%%%%%%%%%%%%%%%%%%%%%%%%%%

%The integral formula~\eqref{WQ-Rn} of the Weyl quantization relies
%on the affine structure of $\bbR^n$ and the trivialization of $V$.
%A generalization to arbitrary manifolds and vector bundles
%requires the choice of connections.
To generalize the above Weyl quantization to
vector bundles over  any smooth manifolds,
we will need connections. We fix some notations concerning exponential maps
and parallel transports below.

A connection $\nabla$ on $TM$ induces an exponential map, indexed by $x\in M$,
\begin{align*}
 \exp_x:U_x &\longrightarrow M, %\\
% v&\mapsto \exp_x v:=\gamma_{x,v}(1),
\end{align*}
where $U_x$ is an open neighborhood of $0\in T_xM$.
%Since $\exp_x$ is a local diffeomorphism,  one can choose $U_x$ such that $\exp_x$ is injective.
We choose $U_x$ such that $\exp_x$ is a diffeomorphism
onto its image. In addition, let $\nabla^V$ be a connection
on a vector bundle $V\rightarrow M$.
Then we have parallel transport maps:
$$    %\[\label{d:Txy}  \begin{split}
  \cT_{x,y}^V :V_x \longrightarrow V_y , %\\
%  \phi &\longmapsto \tilde\gamma_{(x,v,\phi)}(1).
% \end{split}
$$  %\]
for any pairs of points $x,y\in M$ such that $y=\exp_x v$ with $v\in U_x$.
Exponential map and parallel transport together induce
 a local isomorphism of vector bundles
\begin{align}\nonumber
  \cT_x :  U_x \times V_x&\longrightarrow  V\\ \label{map:T}
  (v,\phi) &\longmapsto  (y, \cT^V_{x,y}\,\phi),
\end{align}
where $y=\exp_x \!v$.
Consider a cut-off function $\chi\in\iCinfty(TM)$, i.e.
a function which equals to $1$ in a neighborhood of the zero section
of $TM$ such that the support of $\chi(x,\cdot)$ is included in $U_x$.
%The product of the pull-back map $\cT^*_x$ with $\chi(x,\cdot)$
Define a map
\[\label{map:T*}
\chi(x,\cdot)\cT^*_x :\Ga(V)\longrightarrow \iCinfty(T_xM)\otimes_\bbR V_x,
\]
by setting
$$
\left(\chi(x,\cdot)\cT_x^*\psi\right)(v)=\chi(x,v)\left(\cT_x^*\psi\right)(v),\qquad\forall v\in T_xM,\psi\in\Ga(V).
$$
Note that by definition, if $v\in U_x$ and $y=\exp_x \!v$, we have
\[\label{def:pull-back}
\left(\cT_x^*\psi\right)(v)=\cT_{y,x}^V\left( \psi(y)\right).
\]
The connections $\nabla$ and $\nabla^V$ induce a connection
on the vector bundle $\cS TM\otimes\End V$, and therefore we have a map
analogous to \eqref{map:T*}. Since $\cO(T^*[2]M)\cong\Ga(\cS TM)$,
this map reads as
$$
\chi(x,\cdot)\cT^*_x:\cO(T^*[2]M)\otimes_\Cinfty\Ga(\End V)\longrightarrow \cO\big(T^*[2](T_xM)\big) \otimes_\bbR \End(V_x).
$$
If $(v,p)\in T^*U_x$ and $y=\exp_x \!v$, for any
$F\in \cO(T^*[2]M)\otimes_\Cinfty\Ga(\End V)$,
 we have the following
 explicit formula
\[\label{def:pull-backT}
\left(\cT_x^*F\right)(v,p):=\cT^V_{y,x}\circ F(y,\cT_{x,y} \, p)\circ\cT^V_{x,y},
\]
where $\cT_{x,y}:T^*_xM\to T^*_yM$ is the parallel transport map
induced by $\nabla$.

%%% 4.3 %%%%%%%%%%%%%%%%%%%%%%%%%%%%%%%%%%%%%%%
\subsection{Weyl quantization on $T^*[2]M$}
%%%%%%%%%%%%%%%%%%%%%%%%%%%%%%%%%%%%%%%%%%%%%%%

The  Weyl quantization integral formula \eqref{WQ-Rn}
has been extended to arbitrary cotangent bundles $T^*M$
with the help of a connection on $TM$ \cite{Und78}.
%It was further generalized by Widom, who introduced a total symbol map for
It has been further generalized to differential operators
acting on any vector bundle $V$ over $M$,
using an additional connection on the vector bundle $V$ \cite{Wid80}.
We recall the construction briefly below.

\begin{defn}
By a Weyl quantization map, we mean a map
$$\cQ^M_\hbar : \cO(T^*[2]M)\otimes_\Cinfty\Ga(\End V)\rightarrow\End(\Ga(V))$$
defined by
\[\label{Def:cQ}
(\cQ^M_\hbar \!(F) \,\psi)(x)=\tfrac1{(2\pi\hbar)^n}\int_{T^*_xM\oplus T_xM}  \e^{-\frac{\bi}{\hbar}\<p,v\>}\;\big(\cT^{*}_{x}F\big)(v/2,p)\cdot\big(\cT^*_x\psi\big)(v)\;\chi(x,v)\ud p\,\ud v,
\]
for all $x\in M$, $\psi\in\Gamma(V)$ and $F\in\cO(T^*[2]M)\otimes_\Cinfty\Ga(\End V)$.
\end{defn}

The map $\cQ^M_\hbar$ depends on the choice of connections  on $TM$ and $V$
defining
 the pull-backs $\cT^*_x\psi$ and $\cT^{*}_{x}F$.
Note that a priori, $\cQ^M_\hbar$ also depends  on the cut-off
 function $\chi\in\iCinfty(TM)$.

For any $x\in M$, the projection $\bbR^n\times V_x\rightarrow \bbR^n$
defines a trivial vector bundle on $\bbR^n$. Denote by
$$
\cQ^{\bbR^n}_\hbar :\cO(T^*[2]\bbR^n)\otimes_{\bbR}\Ga(\End V_x)\longrightarrow\D(\bbR^n,\bbR^n\times V_x)\\
$$
the corresponding Weyl quantization map, given
by \eqref{WQ-ND} or equivalently  \eqref{WQ-Rn}.
The maps $\cQ^{\bbR^n}_\hbar$ and  $\cQ^M_\hbar$ are related
by parallel transportation as indicated below.

\begin{lem}\label{prop:cQ:Q}
For any $x\in M$, $\psi\in\Ga(V)$ and $F\in\cO(T^*[2]M)\otimes_\Cinfty\Ga(\End V)$, we have
$$
\cT^*_x\left(\cQ^M_\hbar \!(F)\,\psi \right)(w)= \left(\cQ^{\bbR^n}_\hbar \!(\cT_x^*F)\,\cT_x^*\psi\right)(w), \qquad \forall w\in U_x,
$$
where $U_x$ is an open neighborhood of $0\in T_xM$,
on which $\exp_x$ is a diffeomorphism onto its image.
\end{lem}

\begin{proof}
 For all $w\in U_x$,
Eqns \eqref{def:pull-back} and \eqref{Def:cQ} imply that
\[\label{Lem4.4:Eq1}%$$
\cT^*_x\left(\cQ^M_\hbar \!(F)\,\psi \right)(w)=\tfrac1{(2\pi\hbar)^n}\int_{T^*_yM\oplus T_yM}  \e^{-\frac{\bi}{\hbar}\<p,v\>}\;\cT^V_{y,x}\Big[\big(\cT^{*}_{y}F\big)(v/2,p)\cdot\big(\cT^*_y\psi\big)(v)\Big]\;\chi(y,v)\ud p\,\ud v,
\]
where $y=\exp_x\!w$.
The equality $\exp_y\!v=\exp_x(w+\cT_{y,x}v)$ and Eq.\ \eqref{def:pull-back} lead to
\[\label{Lem4.4:Eq2}%$$
\left(\cT^*_{y}\psi\right)(v)=\cT^V_{x,y}\,\left[\left(\cT_x^*\psi\right)(w+\cT_{y,x}v)\right].
\]%$$
%by applying  formula \eqref{def:pull-back} to
%$\left(\cT^*_{y}\psi\right)(v)$ and $\left(\cT_x^*\psi\right)(w+\cT_{y,x}v)$.
Similarly, using the identity
 $\exp_y\!\frac{v}{2}=\exp_x(w+\frac{1}{2}\cT_{y,x}v)$ and
Eq.\ \eqref{def:pull-backT}, we have
\[\label{Lem4.4:Eq3}%$$
\cT^V_{y,x}\circ\left[\left(\cT^*_{y}F\right)(v/2,p)\right]=\left[\left(\cT_x^*F\right)(w+\tfrac{1}{2}\cT_{y,x}v,\cT_{y,x}p)\right]\circ\cT^V_{y,x}.
\]%$$
As parallel transport in the fibers of $TM\oplus T^*M$
preserves the duality pairing $\<p,v\>$ and the measure $\ud p\ud v$,
by Eqns \eqref{Lem4.4:Eq1}-\eqref{Lem4.4:Eq3} and the change of variables
$(\cT_{y,x}v,\cT_{y,x}p)\mapsto (v,p)$, we obtain
\begin{multline*}
\cT^*_x\left(\cQ^M_\hbar\!(F)\,\psi\right)(w)= \tfrac1{(2\pi\hbar)^n}\int_{T^*_xM\oplus T_xM}  \e^{-\frac{\bi}{\hbar}\<p,v\>}\;\big(\cT^{*}_{x}F\big)(w+v/2,p)\cdot\big(\cT^*_x\psi\big)(w+v)\\
\chi(\exp_x\!w,\cT_{x,\,\exp_x\!w}v)\ud p\,\ud v.
\end{multline*}
Since $F$ is polynomial in $p$,
its Fourier transform w.r.t.\ $p$ is a distribution supported at $0$.
As the function $v\mapsto\chi(x,v)$ is equal to $1$
in a neighborhood of the zero section,
the right hand side in the above equation reduces to
$\big(\cQ^{\bbR^n}_\hbar\!(\cT_x^*F)\,\cT^*_x\psi\big)(w)$
(see formula \eqref{WQ-Rn}).
The conclusion thus follows.
\end{proof}

As a straightforward consequence, we  have

\begin{prop}\label{prop:chi}
Given any  connections on $TM$ and $V$, there is
a unique Weyl quantization map:
$$\cQ^M_\hbar : \cO(T^*[2]M)\otimes_\Cinfty\Ga(\End V)\rightarrow\D(M,V).$$
That is, Formula \eqref{Def:cQ} does not depend on the choice
of cut-off functions $\chi$.
Moreover, $\cQ^M_\hbar$ is a linear isomorphism.
\end{prop}

In what follows, we assume a choice of connections on $TM$ and $V$ is made
and we refer to the corresponding map $\cQ^M_\hbar$ as the Weyl
quantization map.

Let $V$ %:=V_0\otimes\Vol^{1/2}$,  with $V_0$
be a vector bundle endowed with a %fiberwise
pseudo-Hermitian pairing %. It induces a pairing,
$\langle\cdot,\cdot\rangle_V:\Ga(V)\times\Ga(V)\rightarrow\Ga(\Vol)\otimes\bbC$
 as in \eqref{PairingS}.
Then  $(\phi,\psi):=\int_M\langle\phi,\psi\rangle_V$
defines a pseudo-Hermitian scalar product on $\Ga_c(V)$.
Denote by $*:\D(M,V)\to\D(M,V)$ the adjoint operation
associated to $(\cdot,\cdot)$ and by
$*_V:\Ga(\End V)\to\Ga(\End V)$ the adjoint operation associated to
$\langle\cdot,\cdot\rangle_V$.
Abusing notation, we also denote by $*_V$ the obvious extension to
$\cO(T^*[2]M)\otimes_\Cinfty\Ga(\End V)$
%\ni F_0\otimes F_1$, 
defined by $(F_0\otimes F_1)^{*_V}=F_0\otimes F_1^{*_V}$,
$\forall F_0\otimes F_1\in  \cO(T^*[2]M)\otimes_\Cinfty\Ga(\End V)$.

\begin{prop}\label{prop:Qadjoint}
Assume that the following equality holds
\[\label{nabla-pairing}
\nabla_X\langle\phi,\psi\rangle_V=\langle\conn_X^V\phi,\psi\rangle_V+\langle\phi,\conn_X^V\psi\rangle_V,
\qquad\forall X\in{\mathfrak X}(M),\phi,\psi\in\Ga(V).
\]
Then the Weyl quantization $\cQ^M_\hbar$ satisfies the property
$$%\[\label{AdjointcQ}
\cQ^M_\hbar\!(F)^*=\cQ^M_\hbar\!(F^{*_V}), %\qquad \forall F\in\cO(T^*[2]M)\otimes_\Cinfty\Ga(\End V).
$$%\]
for all $F\in\cO(T^*[2]M)\otimes_\Cinfty\Ga(\End V)$.
\end{prop}

\begin{proof}
Eqns \eqref{def:pull-back}--\eqref{Def:cQ} yield the formula
$$
(\cQ^M_\hbar\!(F)\, \psi)(x)=\tfrac1{(2\pi\hbar)^n}\int_{T^*_xM\oplus T_xM}  \e^{-\frac{\bi}{\hbar}\<p,v\>}\;\cT^{V}_{y,x} \Big[F(y,\cT_{x,y}p)\cdot\cT^V_{z,y}\big(\psi(z)\big)\Big]\chi(x,v)\;\ud p\,\ud v,
$$
where $y=\exp_x (v/2)$ and $z=\exp_x(v)$.
By Eq.\ \eqref{nabla-pairing},
the parallel transports $\cT$ and $\cT^V$
satisfy the relation
$$
\cT_{y,x}\big(\left<\phi(y),\psi(y)\right>_V\big)=\left<\cT^V_{y,x}\big(\phi(y)\big),\cT^V_{y,x}\big(\psi(y)\big)\right>_V,
\qquad \forall \phi,\psi\in\Ga(V),
$$
which implies
$$
\left<\cT^{V}_{y,x}\Big[F(y,\cT_{x,y}p)^{*_V}\cdot\cT^V_{z,y}\big(\phi(z)\big)\Big],\psi(x)\right>_V
= \cT_{z,x} \left<\phi(z),\cT^{V}_{y,z}\Big[F(y,\cT_{x,y}p)\cdot\cT^V_{x,y}\big(\psi(x)\big)\Big]\right>_V.
$$
We deduce that
\begin{multline}
(\cQ^M_\hbar\!(F^{*_V})\,\phi,\psi)= \tfrac1{(2\pi\hbar)^n}\times\\ \label{Eq:Qphipsi}
\int_M
\int_{T_xM\oplus T^*_xM}\cT_{z,x}\left<\phi(z),\cT^{V}_{y,z}\Big[F(y,\cT_{x,y}p)\cdot\cT^V_{x,y}\big(\psi(x)\big)\Big]\right>_V\;\e^{\frac{\bi}{\hbar}\<p,v\>}\chi(x,v)\;\ud p\,\ud v.
\end{multline}
%Let us recall that, for integration of $1$-densities
%on $TM\oplus T^*M$, the change of variables formula reads as   %\JP{New sentence}
%$$
%\int_{TM\oplus T^*M} (\Phi^*G)(X)=\int_{TM\oplus T^*M}\cT_{Z,X} \left(G\circ\Phi(X)\right)=\int_{TM\oplus T^*M} G(Z),
%$$
%where $\Phi:X\mapsto Z$ is a diffeomorphism
%and $G$ a $1$-density.
%We apply the above formula to Eq.\ \eqref{Eq:Qphipsi}
Consider the change of variables $(x,v,p)\mapsto (z,v',p')$
 in Eq.\ \eqref{Eq:Qphipsi},
where $z=\exp_x v$, $v'=\cT_{x,z}v$ and $p'=\cT_{x,z}p$.
Since  $\ud p\,\ud v=\cT_{z,x}(\ud p'\,\ud v')$ and $\<p,v\>=\<p',v'\>$,
we obtain
\begin{multline*}
(\cQ^M_\hbar\!(F^{*_V})\,\phi,\psi)= \tfrac1{(2\pi\hbar)^n}\times\\
\int_M \left<\phi(z), \int_{T_zM\oplus T^*_zM}
\cT^{V}_{y,z}\Big[F(y,\cT_{z,y}p')\cdot\cT^V_{x,y}\big(\psi(x)\big)\Big]\;\e^{-\frac{\bi}{\hbar}\<p',v'\>}\chi(x,\cT_{z,x}v')\;\ud p'\,\ud v'\right>_V,
\end{multline*}
where $y=\exp_z(v'/2)$  and $x=\exp_z(v')$.
Using Eq.\ \eqref{Def:cQ}, we conclude that
$(\cQ^M_\hbar\!(F^{*_V})\,\phi,\psi)=(\phi,\cQ^M_\hbar\!(F)\,\psi)$, $\forall \phi,\psi\in\Ga_c(V)$,
and the result follows.
\end{proof}

Finally, let us describe an explicit formula for
the Weyl quantization $\cQ^M_\hbar$ in low degrees.
By abuse of notation, we also denote by $\nabla$
the induced connection on $TM\otimes \End V$.
For a vector bundle $W$ over $M$, let
\[\label{def:trace}
  \Tr :\Ga(W\otimes TM\otimes T^*M)\longrightarrow\Ga(W)
\]
be  the trace map.
%with $W=\End V$.

\begin{prop}\label{Qexplicit}
The Weyl quantization $\cQ_\hbar^M$ satisfies
\[\label{Formula:Q}
\cQ_\hbar^M\!(F)= F \quad\text{and}\quad
 \cQ_\hbar^M\!( X\otimes F) = \frac{\hbar}{\bi}\left[F\conn_{X}^V +\tfrac12\Tr\,\nabla\big( X\otimes F\big)\right],
\]
for all $F\in\Ga(\End V)$ and $X\in\cO_2(T^*[2]M)\cong{\mathfrak X}(M)$.
\end{prop}

\begin{proof}
The first case is trivial. For the second equation,
choose an open neighborhood $U_x$ of
$0\in T_xM$ such that $\exp_x:U_x\rightarrow\exp_x(U_x)$
is a diffeomorphism.
Pulling back a Cartesian coordinate system on $T_xM\times V_x$ by $\cT_x^{-1}$,
we obtain normal coordinates $(x^i)$ centered at~$x$
and a trivialization of $V$ over $\exp_x(U_x)\subset M$.
%provided by the frame parallel along the geodesic radii emanating from $x$.
Denote by $(x^i,p_i)$ the induced canonical coordinates on $T_x^*(\exp_x(U_x))$.
Under such  coordinates, we can write $X(x)=X^i(x)p_i$, with
 $X^i\in\iCinfty(\exp_x(U_x))$.
By Lemma~\ref{prop:cQ:Q} and Eq.\ \eqref{WQ-ND},
we have
\[\label{Eq:localQ}
  \left(\cQ_\hbar^M\!(X\otimes F)\,\psi\right)(x) =
\left(F(x)X^i(x)\frac{\hbar}{\bi}\frac{\partial}{\partial x^i} +\frac{\hbar}{2\bi}\frac{\partial}{\partial x^i} \Big(F(x)X^i(x)\Big)\right)\psi\,(x),
\qquad \forall \psi\in\Ga(V).
\]
Since $(x^i)$ are normal coordinates at $x$,
the partial derivatives $\frac{\partial}{\partial x^i}$
coincide with covariant derivatives at the point $x$.
Hence, the second  equation of  \eqref{Formula:Q} holds  at~$x$.
The point $x$ being arbitrary, this implies the result.
\end{proof}

%%% 4.4 %%%%%%%%%%%%%%%%%%%%%%%%%%%%%%%%%%%%%%%
\subsection{Weyl quantization on symplectic graded manifolds of degree $2$}
%%%%%%%%%%%%%%%%%%%%%%%%%%%%%%%%%%%%%%%%%%%%%%%

Let $(E,\mg)$ be a pseudo-Euclidean vector bundle
equipped with a metric connection $\nabla^E$.
Assume that $(E,\mg)$ admits a spinor bundle $S$.
Choose a connection $\nabla$ on $TM$,
and choose a spinor connection on $S$,
compatible with $\nabla^E$ (see Eq.\ \eqref{compatibilitySE}).
They induce a Weyl quantization $\cQ^M_\hbar$, valued in $\D(M,S)$.

\begin{defn}\label{def:WQ}
The Weyl quantization~$\WQ_\hbar$ on the symplectic graded manifold
$(T^*[2]M\oplus E[1],\om_{\mg,\nabla^E})$ is defined
by the following commutative diagram
\[\label{WQ}
\xymatrix{
 \gsmooth^\bbC(T^*[2]M\oplus E[1])\ar[d]_{\Qodd_\hbar}\ar[rr]^{\WQ_\hbar} && \D(M,S) \\
 \gsmooth(T^*[2]M)\otimes_\Cinfty\Ga(\End S)\ar[urr]_{\cQ^M_\hbar} &&
}
\]
The vertical map reads as
\[\label{Def:gamma}
\g_\hbar := \Id\otimes\Big(\frac{\hbar}{\bi}\Big)^{\k/2}\g, \quad \text{ on } \cO^\bbC(T^*[2]M)\otimes_\Cinfty\Ga(\wedge^\k E),
\]
where $\g$ is the standard Clifford quantization map (see \eqref{Chevalley}).
\end{defn}

The Weyl quantization $\WQ_\hbar$ extends those studied  in \cite{Get83,Vor99},
with $E$ being $TM$ and $TM\oplus T^*M$, respectively.

\begin{thm}\label{thm:WQ}
%The Weyl quantization $\WQ_\hbar$ satisfies the following properties.
For any $F\in\cO^\bbC_k(T^*[2]M\oplus E[1])$, $k\in\bbN$, and
 $G\in\cO^\bbC(T^*[2]M\oplus E[1])$, we have
\begin{enumerate}[label=(\roman*)]%[(i)]
\item\label{p:1} $\WQ_\hbar(F)=(\hbar/\bi)^{k/2}\;\WQ_\bi(F)$;
\item\label{p:ii}  $\sigma_k\circ \WQ_\hbar(F) = \left({\hbar}/{\bi}\right)^{k/2} F$
(see Definition \ref{def:symbols} for the principal symbol map  $\sigma_k$);
% \eqref{Diagram:sigma};
\item\label{p:iv}
$\WQ_\hbar(F)\in\D^+(M,S)$ if $k$ is even
and $\WQ_\hbar(F)\in\D^-(M,S)$ if $k$ is odd (see ~\eqref{ParityD});
\item\label{p:iii} $[\WQ_\hbar(F),\WQ_\hbar(G)]=\frac{\hbar}{\bi}\WQ_\hbar(\{F,G\})+O(\hbar^2)$.
\end{enumerate}
\end{thm}

\begin{proof}
As $\cO^\bbC_k(T^*[2]M\oplus E[1])=
\bigoplus_{2\ell+\k=k}\cO_{\ell}(T^*[2]M)\otimes_\Cinfty\cO_\k(E[1])$,
it suffices to check Properties \ref{p:1}, \ref{p:ii} and \ref{p:iv}
on homogeneous functions $F\in\cO_{\ell}(T^*[2]M)\otimes_\Cinfty\cO_\k(E[1])$.
This can be easily checked using Lemma \ref{prop:cQ:Q} and Eqns \eqref{WQ-ND}, \eqref{Def:gamma}.

Let $F\in\cO^\bbC_k(T^*[2]M\oplus E[1])$ and
$G\in\cO^\bbC_l(T^*[2]M\oplus E[1])$. Set
$$
H=\WQ_\bi^{-1}\left( [\WQ_\bi(F),\WQ_\bi(G)]-\WQ_\bi(\{F,G\})\right).
$$
Using Proposition~\ref{Thm:Symbol} together with \ref{p:ii} and \ref{p:iv},
we deduce that $H$ has degree $k+l-4$. The claim \ref{p:iii} thus
follows  from \ref{p:1}.
\end{proof}

Next we provide an explicit expression for
the Weyl quantization $\WQ_\hbar$ in low degrees.
By abuse of notation, we denote by $\nabla$ the connection
on $\wedge E\otimes TM$ induced by the connections on $E$ and $TM$.
From the definition of $\WQ_\hbar$ and Proposition~\ref{Qexplicit}, we deduce

\begin{prop}\label{c:pqf}
The Weyl quantization $\WQ_\hbar$ satisfies
\begin{align}\nonumber
 \WQ_\hbar(F)&= \left(\tfrac{\hbar}{\bi}\right)^{\k/2}\g(F),\\ \label{WQT}
%\quad \text{and}\quad
\WQ_\hbar( F X) &=\left(\tfrac{\hbar}{\bi}\right)^{1+\k/2}\left[ \g(F)\conn_{X}^S +\tfrac12\g\big(\Tr\, \nabla (XF)\big)\right],
\end{align}
for all $F\in\cO^\bbC_\k(E[1])$ and $X\in\cO_2(T^*[2]M)\cong{\mathfrak X}(M)$.
Here, the trace map $\Tr$ is defined as in \eqref{def:trace}, with $W=\wedge^\k E\otimes\bbC$.
\end{prop}

To proceed further, we assume that $(E,\mg)$
admits a twisted spinor bundle
$$\sS:=S\otimes(\det S^*)^{1/N}\otimes\Vol^{1/2}, $$
with $N$ being the rank of $S$.
A connection $\nabla$
on $TM$ and a spinor connection on $S$
 induce a spinor connection $\nabla^\sS$ on $\sS$
and a Weyl quantization map
\[\label{WQsS}
\WQ_\hbar:\gsmooth^\bbC(T^*[2]M\oplus E[1])\longrightarrow\D(M,\sS).
\]
The latter satisfies in particular Theorem \ref{thm:WQ}
and Proposition \ref{c:pqf}. Obviously, one should replace
$\nabla^S$ by $\nabla^\sS$ in Eq.\ \eqref{WQT}.
Considering $\sS$ allows for further structures
on $\D(M,\sS)$ (see Section \ref{ParReal}).
We study below their interplay with the map $\WQ_\hbar$.
For all $f\in\iCinfty(M)$, $\xi\in\Ga(E)$ and $X\in{\mathfrak X}(M)$,
let
$$
\tau(f)=f,\quad \tau(\xi)=\bi\xi,\quad  \tau(X)=X.
$$
It is easy to see that $\tau$ extends to
a $\C$-antilinear involutive antiautomorphism of
$\cO^\bbC(T^*[2]M\oplus E[1])$.
 On $\cO^\bbC_\k(E[1])$, we have
$\tau(\xi^{i_1}\cdots\xi^{i_\kappa})=i^{\kappa}(-1)^{\frac{\kappa(\kappa-1)}{2}}\xi^{i_1}\cdots\xi^{i_\kappa}$.
Hence, for a homogeneous element $F\in\cO_k^\bbC(T^*[2]M\oplus E[1])$ of degree $k$, the involution $\tau$ satisfies
\[\label{Eq:tau}
\tau(F)=\begin{cases}  F  \quad\text{for }k\text{ even},\\
\bi F  \quad\text{for }k\text{ odd.}
\end{cases}
\]

\begin{prop}\label{p:v}
Let $\hbar\in\bbR^\times$. The Weyl quantization map as in  \eqref{WQsS}
 satisfies
$$\WQ_\hbar \circ\tau(\cdot)=\WQ_\hbar(\cdot)^*,$$
where $*$ is the adjoint operation defined in Lemma \ref{lem:*exist}.
\end{prop}
\begin{proof}
Given the local nature of the maps $\WQ_\hbar$, $\tau$ and $*$,
it suffices to work over a contractible open subset $U\subset M$.
Let $\<\cdot,\cdot\>_U$ be a spinor pairing (see Proposition \ref{Prop:pairing_S})
%pseudo-Hermitian pairing on $\sS|_{U}$, as in Proposition \ref{Prop:pairing_S}.
and  $*_\sS:\Gamma(\End \sS|_{U})\to \Gamma(\End \sS|_{U})$ its adjoint operation.
As $\hbar$ is real, we have $\g_\hbar(\xi)^{*_\sS}=\bi\g_\hbar(\xi)= \g_\hbar(\tau(\xi))$
for all $\xi\in\Ga(E|_{U})$.
Since both maps $*_\sS$ and $\tau$ are $\bbC$-antilinear antiautomorphisms,
 we have $\g_\hbar (\tau (F))=\big(\g_\hbar (F)\big)^{*_\sS}$,
for all $F\in\cO^\bbC(E[1]|_{U})$. 
Extending  $*_\sS$ to the tensor product
$\cO^\bbC(T^*[2]U)\otimes_\Cinfty \cO^\bbC(E[1]|_{U})$ in an obvious way,
we obtain that
$$%\[\label{AdjointGamma}
\g_\hbar (\tau (F))=\big(\g_\hbar (F)\big)^{*_\sS},
$$%\]
for all $F\in\cO^\bbC(T^*[2]U\oplus E[1]|_{U})$.
According to  Proposition \ref{Prop:pairing_S},
the pair of connections $(\nabla,\nabla^\sS)$ satisfies  Eq.\ \eqref{nabla-pairing},
with $V=\sS$.
Applying Proposition~\ref{prop:Qadjoint} yields the result.
\end{proof}
The case $\hbar =\bi$ plays a peculiar role.

\begin{prop}\label{cor:WQ}
Setting $\WQ:=\WQ_\bi$, we have
\begin{enumerate} [label=(\roman*)]
\item\label{p:2} $\WQ$ restricts to an $\bbR$-linear isomorphism
$\WQ : \cO(T^*[2]M\oplus E[1])\to \D(M,\sS)_\bbR$,
where the real algebra $\D(M,\sS)_\bbR$ is defined in Proposition \ref{Dreal},
\item\label{p:3} $\WQ(F)^*=(-1)^{\left\lfloor k/2\right\rfloor}\WQ(F)$,
for all $F\in\cO_k(T^*[2]M\oplus E[1])$, with $k\in\bbN$
and $\left\lfloor k/2\right\rfloor$ the integer part of $k/2$.
\end{enumerate}
\end{prop}

\begin{proof}
Complex conjugation in $\bbC$ extends naturally to
the algebras $\cO^\bbC(T^*[2]M\oplus E[1])$,
$\cO(T^*[2]M)\otimes_\Cinfty\Ga(\bbCl(E))$ and  $\D(M,\sS)$
(see Proposition \ref{Prop:Conjugation}).
To establish \ref{p:2}, it suffices to show that
$$
\WQ(\overline{F})=\overline{\WQ(F)},
$$
for all $F\in\cO^\bbC(T^*[2]M\oplus E[1])$.
Given the local nature of this property,
we can restrict ourselves to a contractible
 open neighborhood $U$ of a point $x\in M$.
Then we have
\begin{align*}
\WQ(\overline{F})&= \left(\cT_x^{-1}\right)^*\circ\cQ^{\bbR^n}_\bi\!\Big(\cT_x^*\big(\g_\bi(\overline{F})\big)\Big)\circ\cT_x^*,
& \text{(by Lemma \ref{prop:cQ:Q})},\\
&= \left(\cT_x^{-1}\right)^*\circ\cQ^{\bbR^n}_\bi\!\Big(\cT_x^*\big(\overline{\g_\bi(F)}\big)\Big)\circ\cT_x^*,
& \quad\text{(by the equality }\g_\bi=\Id\otimes\g).
\end{align*}
%The equality $\g_\bi=\Id\otimes\g$ ensures that
%$\g_\bi(\overline{F})=\overline{\g_\bi(F)}$.
Since $\cT_x$ is built from
connections on the real vector bundles $TM$ and $E$,
we deduce that $\cT_x^*\big(\overline{\g_\bi(F)}\big)=\overline{\cT_x^*\big(\g_\bi(F)\big)}$.
Therefore, from  \eqref{WQ-ND}, it follows that
$$
\cQ^{\bbR^n}_\bi\!\Big(\cT_x^*\big(\overline{\g_\bi(F)}\big)\Big)=
\overline{\cQ^{\bbR^n}_\bi\!\Big(\cT_x^*\big(\g_\bi(F)\big)\Big)}.
$$
Thus
$$
\WQ(\overline{F})=\left(\cT_x^{-1}\right)^*\circ\left[
\overline{\cT_x^*\circ\WQ(F)\circ\left(\cT_x^{-1}\right)^*}\right]
\circ\cT_x^*.
$$
From the equation above,
we deduce that the map $\WQ(F)\mapsto\WQ(\overline{F})$
is an algebra automorphism of $\D(M,\sS)$. By Proposition \ref{c:pqf},
we see that this map indeed  satisfies all the conditions of the map $\overline{\, \cdot \,}$
as in Proposition \ref{Prop:Conjugation}.
% \eqref{Eq:Conjugation}.
Hence, according to the uniqueness in
 Proposition \ref{Prop:Conjugation}, the map
$\WQ(F)\mapsto\WQ(\overline{F})$ must be of the form $\WQ(F)\mapsto\overline{\WQ(F)}$.
Therefore, we have $\WQ(\overline{F})=\overline{\WQ(F)}$.

Let $\hbar\in\bbR^\times$ and $F\in\cO_k^\bbC(T^*[2]M\oplus E[1])$.
Using Theorem~\ref{thm:WQ}   \ref{p:1} and Proposition \ref{p:v},
 we have 
 $$
 \WQ(F)^*=\left(\left(\tfrac{\hbar}{\bi}\right)^{-\frac{k}{2}} \WQ_\hbar(F)\right)^*
 =\left(\tfrac{\hbar}{\bi}\right)^{-\frac{k}{2}} \bi^{-k}\,\WQ_\hbar(\tau(F))
 =\bi^{-k}\,\WQ(\tau(F)).
 $$
The second assertion  of the proposition
thus  follows  from  \eqref{Eq:tau}.
\end{proof}

The Weyl quantization map $\WQ_\hbar$  in  \eqref{WQsS} induces
a star-product $\star_\hbar$ on the symplectic manifold $T^*[2]M\oplus E[1]$.
The properties of the map $\WQ_\hbar$, exhibited above,
can be rephrased in terms of properties of $\star_\hbar$.

\begin{cor}\label{cor:MW}
The product defined on $\cO^\bbC(T^*[2]M\oplus E[1])$ by
$$
F\star_\hbar G:=(\WQ_\hbar)^{-1}\big(\WQ_\hbar(F)\circ\WQ_\hbar(G)\big),
$$
is a symmetric star-product, explicitly given by
$$
F\star_\hbar G=\sum_{k=0}^\infty \left(\tfrac{\hbar}{\bi}\right)^k B_{2k}(F,G),
$$
such that, for each $k\in\bbN$,
\begin{enumerate}[label=(\roman*)]
\item $B_{2k}$ is a real bidifferential operator of degree $-2k$ independent of $\hbar$,
\item $B_{2k}$ is symmetric if $k$ is even and skew-symmetric if $k$ is odd,
\item $B_0$ is the multiplication and $B_2$ is the Poisson bracket on $(T^*[2]M\oplus E[1],\om_{\mg,\nabla^E})$.
\end{enumerate}
\end{cor}

\begin{proof}
From  Theorem~\ref{thm:WQ} ~\ref{p:1}, we deduce that
$$
F\star_\hbar G=\sum_{k=0}^\infty \left(\tfrac{\hbar}{\bi}\right)^{k/2} B_k(F,G),
$$
where, for each $k\geq 0$,  $B_k$ is a  bilinear operator of degree $-k$.
By Lemma~\ref{prop:cQ:Q} and Eq.  \eqref{WQ-ND}, we see  that
$B_k$ is a  bidifferential operator. Since $\WQ_\hbar$ preserves the parity,
we deduce that $B_k$ vanishes if $k$ is odd.
From Theorem~\ref{thm:WQ} ~\ref{p:ii} and   ~\ref{p:iii},
it follows  that  $B_0 (F, G)=FG$ and 
$B_1 (F, G)=\{F, G\}$.
%The identifications of $B_0$ and~$B_2$, with the multiplication
%and the Poisson bracket respectively, follow from
%the points~\ref{p:ii} and \ref{p:iii} in Theorem~\ref{thm:WQ}.

Set $\hbar=\bi$. According to  Proposition~\ref{cor:WQ},
the operator $B_{2k}$ is real and (skew-)symmetric
according to the parity of $k$. This concludes the proof of the corollary.
\end{proof}

%% 5 %%%%%%%%%%%%%%%%%%%%%%%%%%%%%%%%%%%%%%%%%%%%%%%%%%%
%%%%%%%%%%%%%%%%%%%%%%%%%%%%%%%%%%%%%%%%%%%%%%%%%%%%%%%%
\section{Applications to Courant algebroids}
%%%%%%%%%%%%%%%%%%%%%%%%%%%%%%%%%%%%%%%%%%%%%%%%%%%%%%%%
%%%%%%%%%%%%%%%%%%%%%%%%%%%%%%%%%%%%%%%%%%%%%%%%%%%%%%%%

%%% 5.1 %%%%%%%%%%%%%%%%%%%%%%%%%%%%%%%%%%%%%%%%%%%%%%%%
\subsection{Definition of Courant algebroids}
%%%%%%%%%%%%%%%%%%%%%%%%%%%%%%%%%%%%%%%%%%%%%%%%%%%%%%%%

A pre-Courant algebroid is a pseudo-Euclidean vector bundle  $(E,\mg)$
over a smooth manifold $M$, together with a vector bundle morphism
$\rho :E \to TM$, called the anchor, and an $\bbR$-bilinear operation
$\llbracket\cdot,\cdot\rrbracket$ on $\Ga(E)$, called the Dorfman bracket,
subject to the following rules:
\begin{align} \nonumber
&  \llbracket\xi,f\cdot\eta\rrbracket = \rho(\xi)[f]\cdot\eta + f\cdot\llbracket\xi,\eta\rrbracket,  \\ \label{Eq:preCourant}
&  \llbracket\xi,\xi\rrbracket = \half \rho^*\ud\left(\mg(\xi,\xi)\right),   \\ \nonumber
&  \rho(\xi)[\mg(\eta,\eta)] = 2\mg\big(\llbracket\xi,\eta\rrbracket,\eta\big),
\end{align}
for all $f\in\iCinfty(M)$ and $\xi,\eta\in\Ga(E)$.
In the second equation above, $\ud$ stands for
the de Rham differential and $\rho^*:T^*M\rightarrow E^*\cong E$
is the dual map of $\rho$. Moreover, if the bracket  satisfies the Jacobi identity
\[\label{Eq:Courant}
\llbracket\xi,\llbracket\eta_1,\eta_2\rrbracket\rrbracket = \llbracket\llbracket\xi,\eta_1\rrbracket,\eta_2\rrbracket + \llbracket\eta_1,\llbracket\xi,\eta_2\rrbracket\rrbracket,
\]
for all $\xi,\eta_1,\eta_2\in\Ga(E)$, then
$(E,\mg,\rho,\llbracket\cdot,\cdot\rrbracket)$ is called
a Courant algebroid \cite{LWX97,Roy99}.

\begin{ex}
For any smooth manifold $M$, the vector bundle $E=TM\oplus T^*M$
admits a standard Courant algebroid structure, where
the anchor is the projection onto the first component
and the pairing and Dorfman bracket are given, respectively by
 %\JP{$\tfrac{1}{2}$ added to match introduction and section 6.}
\begin{align*}
\mg( X+\alpha, Y+\beta)&=\langle \alpha,Y \rangle+\langle \beta,X \rangle,\\
\llbracket X+\alpha, Y+\beta\rrbracket &= [X,Y] + \Lie_X\beta -\ins_Y\ud\alpha,
\end{align*}
for all $X,Y\in{\mathfrak X}(M)\cong\Ga(TM)$ and $\a,\b\in\Ga(T^*M)$.
One can also twist the above bracket
by a closed 3-form $H\in\Omega^3(M)$  \cite{Sev7, Sev01},
$$
\llbracket X+\alpha, Y+\beta\rrbracket = [X,Y] + \Lie_X\beta -\ins_Y\ud\alpha+H(X,Y,\cdot).
$$
For more examples, see e.g.\ \cite{LWX97}.
\end{ex}

%%% 5.2 %%%%%%%%%%%%%%%%%%%%%%%%%%%%%%%%%%%%%%%%%%%%%%%%
\subsection{Courant algebroids and symplectic graded manifolds of degree $2$}
%%%%%%%%%%%%%%%%%%%%%%%%%%%%%%%%%%%%%%%%%%%%%%%%%%%%%%%%

We assume, from now on, that $(E,\mg)$ is endowed with
a metric connection $\nabla^E$ so that the minimal symplectic
realization of the Poisson manifold $E[1]$ is given by
$(T^*[2]M\oplus E[1],\om_{\mg,\nabla^E})$ as in Proposition
\ref{PropRoth} (see Section~\ref{Par:SymplGM}).
The Poisson bracket on $T^*[2]M\oplus E[1]$ is denoted by
$\{\cdot,\cdot\}$ which is of degree~$-2$.
We will use the following identifications without further comments:
\begin{align*}
\cO_0(T^*[2]M\oplus E[1])&\cong\iCinfty(M),\\
\cO_1(T^*[2]M\oplus E[1])&\cong\Gamma(E),\\
\cO_2(T^*[2]M\oplus E[1])&\cong\Gamma(TM\oplus\wedge^2 E).
\end{align*}

Every degree~$3$ function $\Theta\in\cO_3(T^*[2]M\oplus E[1])$
induces a pre-Courant algebroid structure on~$(E,\mg)$
by setting
\[\label{DerivedPoisson}
\begin{split}
  \rho(\xi)[f] &:= \{\{\Theta,\xi\},f\},  \\
   \llbracket\xi,\eta \rrbracket&:= \{\{\Theta,\xi\},\eta\},
 \end{split}
\qquad \forall f\in\iCinfty(M), \xi,\eta \in\Ga(E).
\]
%where $f\in\iCinfty(M)$ and $\xi_1,\xi_2\in\Ga(E)$.
The structural identities \eqref{Eq:preCourant} follow from
the fact that $\{\cdot,\cdot\}$ is a Poisson bracket
which satisfies the relation $\{\xi,\eta\}=\mg(\xi,\eta)$ for all $\xi,\eta\in\Ga(E)$.
Moreover, if $\{\Theta,\Theta\}=0$, then $\llbracket\cdot,\cdot\rrbracket$
satisfies the Jacobi identity \eqref{Eq:Courant} and therefore
$E$ becomes a Courant algebroid.

Conversely, one can construct a degree $3$ function
$\Theta\in\cO_3(T^*[2]M\oplus E[1])$ out of a pre-Courant algebroid
$(E,\mg,\rho,\llbracket\cdot,\cdot\rrbracket)$ as follows.
The anchor map $\rho$, being a bundle map, can be identified with
a section in $\Ga(TM\otimes E^*)$, which is a function of degree $3$
on $T^*[2]M\oplus E[1]$. From $\mg$ and $\llbracket\cdot,\cdot\rrbracket$,
one can define the torsion map $\CTor : \Ga(\wedge^3 E)\rightarrow\iCinfty(M)$ 
by
\[\label{torsion}
\CTor(\xi_1,\xi_2,\xi_3)=\frac{1}{2}\cycl_{123}\;
\mg\left(\frac{1}{3}\big(\llbracket\xi_1,\xi_2\rrbracket-\llbracket\xi_2,\xi_1\rrbracket\big)-\big(\conn_{\rho(\xi_1)}^E\xi_2-\conn_{\rho(\xi_2)}^E\xi_1\big),
\ \xi_3\right),
\]
where  $\xi_1,\xi_2,\xi_3\in\Ga(E)$ and $\cycl_{123}$ denotes
the sum over cyclic permutations. As proved in~\cite{AX03},
the identities \eqref{Eq:preCourant} ensure that
$\CTor$ is $\iCinfty(M)$-multilinear, so that it can be identified with
a section in $\Ga(\wedge^3 E)$. Set $\Theta=\rho-\CTor$.
Then, $\Theta$ is a degree $3$ function which satisfies Eq.~\eqref{DerivedPoisson}.
Moreover, if $\llbracket\cdot,\cdot\rrbracket$
satisfies the Jacobi identity \eqref{Eq:Courant}, 
we have $\{\Theta,\Theta\}=0$.
Thus we recover the following

\begin{thm}[\cite{Roy02}]\label{thm:Theta}
Let $(E,\mg)$ be a pseudo-Euclidean vector bundle over $M$.
There is a bijection between
pre-Courant algebroids $(E,\mg,\rho,\llbracket\cdot,\cdot\rrbracket)$
and degree $3$  functions $\Theta\in\cO_3(T^*[2]M\oplus E[1])$.
They are related via $\Theta=\rho-\CTor$.% and Eq.\ \eqref{DerivedPoisson}.

Moreover, $(E,\mg,\rho,\llbracket\cdot,\cdot\rrbracket)$ is
a Courant algebroid if and only if $\{\Theta,\Theta\}=0$.
\end{thm}

The above function $\Theta$ is called the
\emph{Hamiltonian generating function} of the Courant algebroid.

\begin{rem} %\JP{New. Check please!}
Both $\Theta$ and the Poisson bracket on $T^*[2]M\oplus E[1]$
depend on the choice of a metric connection on $E$.
Via the symplectic diffeomorphism $\Xi_\nabla:\cM\to T^*[2]M\oplus E[1]$,
both of them  admit an intrinsic version on $\cM$.
They satisfy again Eq.\ \eqref{DerivedPoisson} (see \cite{Roy02}).
\end{rem}

%%% 5.3 %%%%%%%%%%%%%%%%%%%%%%%%%%%%%%%%%%%%%%%%%%%%%%%%
\subsection{Dirac generating operators of Courant algebroids}%Courant algebroids and spinor differential operators
%%%%%%%%%%%%%%%%%%%%%%%%%%%%%%%%%%%%%%%%%%%%%%%%%%%%%%%%

There is another approach to generate a Courant algebroid,
via the so-called \emph{Dirac generating operators}  \cite{AX03,CS09}.

From now on, we assume that $(E,\mg)$ admits a spinor bundle $S$
and that $\det(S^*)^{1/N}$ exists, with $N$ being the rank of $S$.
Then, the twisted spinor bundle,
$
\sS:=S\otimes (\det S^*)^{1/N}\otimes\Vol^{1/2},
$
is well-defined.
The algebra of real differential operators $\D(M,\sS)_\bbR$, defined in Proposition \ref{Dreal},
inherits a $\bbZ_2$-grading and a filtration from $\D(M,\sS)$
(cf.\  \eqref{ParityD} and~\eqref{filtration}).
In particular, we have
\begin{gather*}
  \D_0(M,\sS)_\bbR\cong\iCinfty(M)\cong\cO_0(T^*[2]M\oplus E[1]),  \\
  \D_1^{-}(M,\sS)_\bbR\cong\Ga(E)\cong\cO_1(T^*[2]M\oplus E[1]).
\end{gather*}
Since the Weyl quantization map $\WQ:\cO(T^*[2]M\oplus E[1])\to\D(M,\sS)_\bbR$
(see Proposition~\ref{cor:WQ}) preserves  the
parity (see Theorem \ref{thm:WQ}), it induces the following isomorphism:
$$%\[\label{D3odd}
\WQ^{-1}:\D^{-}_3(M,\sS)_\bbR\stackrel{\sim}{\longrightarrow}\cO_3(T^*[2]M\oplus E[1])\oplus\cO_1(T^*[2]M\oplus E[1]).
$$%\]

The Dorfman bracket can be obtained as a derived bracket
of the commutator in~$\D(M,\sS)_\bbR$, for a well-chosen
generating operator $D\in \D^{-}_3(M,\sS)_\bbR$, as shown in \cite{AX03}.
This approach provides a quantum analog
to the Hamiltonian picture for Courant algebroids,
presented in the previous section.  Namely, as the commutator
lowers the order by $2$ and preserves the parity, we can define
\[\label{DerivedCommutator}
 \begin{split}
  \rho(\xi)[f] &:= [[D,\g(\xi)],f]\in\D_0(M,\sS)_\bbR\cong\iCinfty(M),\\
  \llbracket\xi, \eta \rrbracket &:= [[D,\g(\xi)],\g(\eta )]\in\D_1^{-}(M,\sS)_\bbR\cong\Ga(E),
\end{split}
\qquad \forall f\in\iCinfty(M), \xi,\eta \in\Ga(E).
\]
%$\forall f\in\iCinfty(M)$ and $\xi_1,\xi_2\in\Ga(E)$.
%They satisfy the following properties.

We have the following

\begin{prop}\label{prop:DCourant}
Let  $D\in \D^{-}_3(M,\sS)_\bbR$ such that $\sigma_3(D)\neq 0$. Then
\begin{enumerate}[label=(\roman*)]
\item the map $\rho$ and the bracket $\llbracket\cdot,\cdot\rrbracket$
given by  \eqref{DerivedCommutator} define a pre-Courant algebroid
structure on $(E,\mg)$;
\item the pre-Courant algebroid defined by $D$ as in  \eqref{DerivedCommutator}
coincides with the one defined by its principal symbol $\sigma_3(D)$ as in
  \eqref{DerivedPoisson};
\item the operator $D$ generates a Courant algebroid if and only if  $D^2\in\D_2(M,\sS)_\bbR$.
\end{enumerate}
\end{prop}

\begin{proof}
Using Eq.\ \eqref{SymbolPoisson},  we obtain that
the map $\rho$ and the bracket $\llbracket\cdot,\cdot\rrbracket$,
defined by $D$ via  \eqref{DerivedCommutator},
indeed coincide with the ones determined by
$\sigma_3(D)$ via  \eqref{DerivedPoisson}.
Hence, they define a pre-Courant algebroid structure on $(E,\mg)$.

It remains to prove the last assertion. Since $D$ is odd, we have $D^2=\frac{1}{2}[D,D]$.
As the commutator lowers the order by $2$ and preserves the parity,
the operator $[D,D]$ is of order $4$, $2$ or $0$.
By Eq.\ \eqref{SymbolPoisson}, its principal symbol satisfies
$\sigma_4([D,D])=\{\sigma_3(D),\sigma_3(D)\}$.
Hence, $D^2\in\D_2(M,\sS)_\bbR$ if and only if $\{\sigma_3(D),\sigma_3(D)\}$=0.
%The latter is equivalent to
%$(E,\mg,\rho,\llbracket\cdot,\cdot\rrbracket)$ is a Courant algebroid,
The result thus follows  from Theorem~\ref{thm:Theta}.
\end{proof}

In \cite{AX03}, a stronger condition of Dirac generating operators was introduced.

\begin{defn}\cite{AX03}
A \emph{Dirac generating operator} is
an operator $D\in \D^{-}_3(M,\sS)_\bbR$,
such that $\sigma_3(D)\neq 0$ and $D^2\in\D_0(M,\sS)_\bbR$.%\cong\iCinfty(M)$.
\end{defn}

According to Proposition \ref{prop:DCourant},
a Dirac generating operator indeed generates
a Courant algebroid structure on $(E,\mg)$
via \eqref{DerivedCommutator}.

The existence of a Dirac generating operator
for a given Courant algebroid is nontrivial and
was one of the main results of \cite{AX03}.
 Note that Dirac generating operators are not unique:
two Dirac generating operators $D$ and $D'$
defining the same Courant algebroid
have the same principal symbol but may differ by an element
$\g(\xi)\in\D_1^{-}(M,\sS)_\bbR$ satisfying $\{\sigma_3(D),\xi\}=0$.
As an application of our theory,
in what follows,  we construct a Dirac generating operator $D$
via the Weyl quantization map $\WQ$.
It turns out that such  a Dirac generating operator $D$
is independent of any choice of geometric data
and  is therefore   unique. Below, 
following an idea of \v{S}evera \cite{Sev7}, 
we  describe a completely intrinsic characterization of such an
operator  $D$ 
in terms of the  adjoint operation defined in \eqref{adjoint}.
%, which is completely intrinsic.
%The definition of the map $\WQ$ requires to choose
%a connection on $TM$ and a spinor connection on $S$,
%compatible with $\nabla^E$.

\begin{thm}\label{thm:DiracGenOp}
Let $(E,\mg,\rho,\llbracket\cdot,\cdot\rrbracket)$
be a Courant algebroid admitting a twisted spinor bundle~$\sS$.
There exists a unique Dirac generating operator $D\in\D_3^-(M,\sS)_\bbR$ satisfying:
\begin{enumerate}[label=(\roman*)]
\item $D$ generates the given Courant algebroid, and
\item$D^*=-D$, where $*$ is the adjoint operation defined in \eqref{adjoint}.
\end{enumerate}
In fact, we have $D=\WQ(\Theta)$, where $\Theta\in\cO_3(T^*[2]M\oplus E[1])$
is the Hamiltonian generating function of $(E,\mg,\rho,\llbracket\cdot,\cdot\rrbracket)$.
%Contrary to the map $\WQ$, it does not depend on any connection.
%but only on the Clifford action on the twisted spinor bundle $\sS$.
\end{thm}

Note that  the adjoint operation $*$ does not depend on any
 choice of connections according to Proposition
\ref{prop:*}.
We need a lemma first.

\begin{lem}\label{lem:Theta-D}
The Weyl quantization map $\WQ$ induces a linear isomorphism
\[\label{WQ:iso}
\WQ : \cO_3(T^*[2]M\oplus E[1])\stackrel{\sim}{\longrightarrow}\{D\in\D_3^{-}(M,\sS)_\bbR|\;D^*=-D\},
\]
%Moreover the map $\WQ$
which establishes a bijection between
Hamiltonian generating functions and
skew-symmetric Dirac generating operators.
\end{lem}

\begin{proof}
According to Proposition~\ref{cor:WQ}, we have
$\WQ(F)^*=(-1)^{\left\lfloor k/2\right\rfloor}\WQ(F)$
for all $F\in\cO_k(T^*[2]M\oplus E[1])$. Hence,
the map $\WQ$ sends $\cO_1(T^*[2]M\oplus E[1])$
to the space of symmetric (or self-adjoint) operators in $\D^{-}_3(M,\sS)_\bbR$
and $\cO_3(T^*[2]M\oplus E[1])$ to the space of
skew-symmetric operators in $\D^{-}_3(M,\sS)_\bbR$.
Thus, the map \eqref{WQ:iso} is indeed an isomorphism.

Let $\Theta$ be a function of degree $3$ and  $D=\WQ(\Theta)$.
By Corollary~\ref{cor:MW}, we have
$$
2(\WQ)^{-1}(D^2)=\{\Theta,\Theta\}+B_4(\Theta,\Theta)+B_6(\Theta,\Theta),
$$
where $B_6(\Theta,\Theta)$ is of degree $0$ and
$B_4(\cdot,\cdot)$ is a symmetric bidifferential operator.
As $\Theta$ is of odd degree, we have $B_4(\Theta,\Theta)=0$.
As a consequence, we conclude that
$\{\Theta,\Theta\}=0$ is equivalent to $D^2\in\D_0(M,\sS)_\bbR$.
\end{proof}

\begin{proof}[Proof of Theorem \ref{thm:DiracGenOp}]
Let $\Theta\in\cO_3(T^*[2]M\oplus E[1])$ be the
Hamiltonian generating function of $(E,\mg,\rho,\llbracket\cdot,\cdot\rrbracket)$.
By Proposition \ref{prop:DCourant}, a Dirac generating operator $D\in\D^{-}_3(M,\sS)_\bbR$
of $(E,\mg,\rho,\llbracket\cdot,\cdot\rrbracket)$
must satisfy $\sigma_3(D)=\Theta$.
According to Lemma \ref{lem:Theta-D},
there exists a unique  such $D$
satisfying in addition $D^*=-D$,
which  is indeed given by $D=\WQ(\Theta)$.
\end{proof}

Next we will describe an explicit formula for the skew-symmetric
Dirac generating operator $D$. By uniqueness, $D$ does not depend
on any choice of connections. However, the Weyl quantization
map $\WQ$ does.
This is reflected in the formula \eqref{Formula:D} below,
which is written in terms of the
connections $\nabla^\sS$ on $\sS$  and $\nabla$
on $E\otimes TM$.  The latter are induced by a connection
on $TM$ and compatible  connections on $S$ and $E$.
In addition, the formula \eqref{Formula:D} involves the torsion $\CTor\in\Ga(\wedge^3 E)$
(see Eq.\ \eqref{torsion}), a  local frame $(\xi^a)$ of $E$ and
the metric components $\mg_{ab}$,
defined by $\mg^{bc}=\mg(\xi^b,\xi^c)$ and $\mg_{ab}\mg^{bc}=\delta^c_a$.

\begin{cor}\label{c:DGO1}
Let $(E,\mg,\rho,\llbracket\cdot,\cdot\rrbracket)$
be a Courant algebroid. The skew-symmetric Dirac
generating operator is given by the following formula
\[\label{Formula:D}
D=\mg_{ab}\,\g(\xi^a)\conn_{\rho(\xi^b)}^\sS-\g(\CTor)+\frac{1}{2}\g\left(\Tr\,\nabla\rho\right),
\]
where $\nabla\rho\in\Ga(E\otimes TM\otimes T^*M)$ and
$\Tr\,\nabla\rho \in \Ga(E)$ denotes its trace  as defined in \eqref{def:trace}.
\end{cor}

\begin{proof}
By Theorem \ref{thm:DiracGenOp} and Theorem \ref{thm:Theta},
we know that $D=\WQ(\Theta)$ with $\Theta=\rho-\CTor$.
The result follows from Proposition~\ref{c:pqf}.
\end{proof}

%%%5.4%%%%%%%%%%%%%%%%%%%%%%%%%%%%%%%%%%%%%%%%%%%%%%%%%%%%%
\subsection{An alternative formula for Dirac generating operators}
%%%%%%%%%%%%%%%%%%%%%%%%%%%%%%%%%%%%%%%%%%%%%%%%%%%%%%%%

We now prove that the Dirac generating operator
constructed by Alekseev and Xu in \cite{AX03}
coincides with the operator given by Eq.\ \eqref{Formula:D}.
We need to introduce several notations, directly borrowed from \cite{AX03}.

\begin{defn}[\cite{AX03}]
Let $(E,\mg,\rho,\llbracket\cdot,\cdot\rrbracket)$
be a Courant algebroid and $V$ a vector bundle over $M$.
An $\bbR$-bilinear  map
$\Enabla : \Gamma(E) \times \Gamma(V) \to \Gamma(V)$
is called an $E$-connection on $V$ if it satisfies the following conditions:
\begin{align*}
\Econn_{f\xi} v &= f \Econn_\xi v,  \quad \text{and}\\
\Econn_\xi (fv) &= f \Econn_\xi v + \rho(\xi)[f]\cdot v,
\end{align*}
for all $\xi \in \Gamma(E), v \in \Gamma(V)$ and $f \in \iCinfty(M)$.
\end{defn}

\begin{rem}
An ordinary linear connection $\nabla$ on the vector bundle $V$ induces
an $E$-connection as follows:
\[ \label{eq:nablarho}
\Econn_\xi v := \conn_{\rho(\xi)} v,\qquad \forall\xi\in\Ga(E),v\in\Ga(V).
\]
\end{rem}

Let $\Enabla^E$ be an $E$-connection on $E$.
According to \cite{AX03},
\[\label{TorsionEconn}
 \CTor[\Enabla](\xi_1,\xi_2,\xi_3)
  =\frac{1}{2}\cycl_{123}\;\mg\left(\frac{1}{3}\big(\llbracket\xi_1,\xi_2\rrbracket-\llbracket\xi_2,\xi_1\rrbracket \big)-\big(\Enabla^E_{\xi_1}\xi_2-\Enabla^E_{\xi_2}\xi_1\big), \ \xi_3\right)
\]
defines a section in $\Ga(\wedge^3 E)$,
where $\cycl_{123}$ denotes the sum over cyclic permutations.
This is called the torsion of $E$ with respect to $\Enabla^E$.
If $\Enabla^E$ is induced from an ordinary connection $\nabla^E$
as in Eq.\ \eqref{eq:nablarho}, then $\CTor[\Enabla]$ coincides
with the torsion $\CTor$, introduced previously in  \eqref{torsion}.

Let  $(\xi^a)$ be a  local frame of $E$.
For any $\xi\in \Ga (E)$, the function
$\Div\xi\in \iCinfty(M)$, defined by
\[\label{Eq:Div}
\Div\xi=\sum_{a}\mg(\Enabla^E_{\xi^{a}} \xi,\  \xi^{a}),
\]
does not depend on the chosen local frame $(\xi^a)$.
We call $\Div\xi$ the divergence of $\xi$.
Then,
\[\label{EnablaLambdaDefn}
\Enabla^\Lambda_{\xi}s := \Lie_{\rho (\xi)}s -(\Div\xi)s
\]
defines an $E$-connection on the line bundle
 $\wedge^\Top T^{*}M$.  Here $\Lie_{\rho(\zeta)}s$ stands for
the Lie derivative of $s\in\Gamma(\wedge^\Top T^*M)$ along the vector field
$\rho(\zeta)\in\Vect(M)$.

\begin{lem}\label{lem:Enabla_Lambda}
Let $\nabla^E$ be a metric connection on $E$.
It induces an $E$-connection on $E$ as in Eq.\ \eqref{eq:nablarho}
and an $E$-connection $\Enabla^\Lambda$
on $\wedge^\Top T^{*}M$ via Eq.\ \eqref{EnablaLambdaDefn}.
For any linear connection $\nabla^{TM}$ on $TM$,
we have the relation
$$%\[\label{eq:Enabla_Lambda}
\Enabla^\Lambda_\xi=\conn_{\rho(\xi)}^\Lambda+\mg(\xi,\Tr\,\nabla\rho),\qquad\forall\xi\in\Ga(E),
$$%\]
where $\nabla^\Lambda$ and $\nabla$ are induced connections on $\wedge^\Top T^{*}M$
and $E\otimes TM$ respectively.
\end{lem}

\begin{proof}
The Lie derivative on $\wedge^\Top T^{*}M$ satisfies the relation:
\[\label{Eq:LieVol}
\Lie_{\rho(\xi)}=\conn_{\rho(\xi)}^\Lambda+\Tr\,\nabla^{TM}\!\left(\rho(\xi)\right).
\]
%where $\nabla^{TM}$ is the connection on $TM$.
Since $\nabla^E$ preserves the metric, we have
$\nabla^{TM}\!\left( \rho(\xi)\right)=
\mg(\nabla\rho,\xi)+\mg(\rho,\nabla^E\xi)$, which
 is an equality in $\Ga(TM\otimes T^*M)$. By  taking the trace, we obtain
\[\label{Eq:Tr}
\Tr\,\nabla^{TM}\!\left( \rho(\xi)\right)=\mg(\Tr\,\nabla\rho,\xi)+\Div\xi,
\]
where the divergence of $\xi$ is computed for the $E$-connection on $E$
induced by $\nabla^E$ as in \eqref{Eq:Div}.
The conclusion thus follows by combining  Eqns \eqref{EnablaLambdaDefn},
\eqref{Eq:LieVol} and \eqref{Eq:Tr}.
\end{proof}

An $E$-connection on $E$ preserves the metric if
$$
\rho(\xi_1)[\mg(\xi_2,\xi_3)]=\mg(\Econn_{\xi_1}^E\xi_2,\xi_3)+\mg(\xi_2,\Econn_{\xi_1}^E \xi_3),
$$
for all $\xi_1,\xi_2,\xi_3\in\Ga(E)$. An $E$-connection
$\Enabla^E$ on $E$ and an $E$-connection $\Enabla^S$ on $S$
are compatible if $\Enabla^E$ preserves the metric,
% and both $\Enabla^E$ and  $\Enabla^S$ induce the same $E$-connection
and their induced connections on the Clifford algebra bundle
coincide under the isomorphism  $\bbCl(E)\cong\End S$.
By Eq.\ \eqref{EnablaLambdaDefn}, we see that such a pair of compatible
$E$-connections induces an $E$-connection on
the twisted spinor bundle $\sS=S\otimes (\det S^*)^{1/N}\otimes\Vol^{1/2}$.

\begin{prop}\label{p:DGO2}
Let $(E,\mg,\rho,\llbracket\cdot,\cdot\rrbracket)$
be a Courant algebroid  with twisted spinor bundle~$\sS$.
The unique skew-symmetric Dirac generating operator is given by
\[\label{Eq:Dbis}
D=\mg_{ab}\,\g(\xi^a)\Econn_{\xi^b}^\sS-\g(\CTor[\Enabla]),
\]
where $\Enabla^\sS$ is the $E$-connection on $\sS$
induced by compatible $E$-connections on $S$ and $E$.
Here, $(\xi^a)$ is a  local  frame  of $E$ and
$\mg_{ab}$ are the metric components.
\end{prop}
\begin{proof}
According to \cite[Theorem 4.3]{AX03}, the formula \eqref{Eq:Dbis}
does not depend on the choice of $E$-connections on $E$ and $S$.
Hence, we can start with $E$-connections induced by ordinary linear connections.
In this case, we have  $\CTor[\Enabla]=\CTor$ and
$\Econn_\xi^\sS=\conn_{\rho(\xi)}^\sS+\frac{1}{2}\mg(\xi,\Tr\,\nabla\rho)$
 by Lemma~\ref{lem:Enabla_Lambda}.
Substituting both equalities in Eq.\ \eqref{Eq:Dbis} yields Eq.\ \eqref{Formula:D}.
The result thus follows.
\end{proof}

%%%5.5%%%%%%%%%%%%%%%%%%%%%%%%%%%%%%%%%%%%%%%%%%%%%%%%%%%%%
\subsection{ A Courant algebroid invariant}
%%%%%%%%%%%%%%%%%%%%%%%%%%%%%%%%%%%%%%%%%%%%%%%%%%%%%%%%

Since a skew-symmetric Dirac generating operator $D$ is unique,
then $D^2\in\iCinfty(M)$ is an invariant of the Courant algebroid.
It is natural to ask how this invariant can be described geometrically.

Let  $(E,\mg,\rho,\llbracket\cdot,\cdot\rrbracket)$ be a Courant algebroid.
By abuse of notation, we denote by the same symbol
 $\mg$ the metric on $\wedge E$
induced by the metric on $E$ as follows: %\JP{convention for $\mg$ on $\wedge E$ added}
$\mg(\xi,\eta)=\det(\mg(\xi_i,\eta_j))$ for all
$\xi=\xi_1\wedge\cdots\wedge\xi_k$ and
$\eta=\eta_1\wedge\cdots\wedge\eta_k$
in $\Ga(\wedge^k E)$.
Using a metric connection $\nabla^E$ on $E$
and a linear connection $\nabla^{TM}$ on $TM$, we define $f_E\in\iCinfty(M)$ by
\[\label{Inv fE}
  f_E  = \mg(\CTor,\CTor)-\frac{1}{4}\mg\left(\Tr\,\nabla\rho,\Tr\,\nabla\rho\right)-\frac{1}{2}\Div\left(\Tr\,\nabla\rho\right),
\]
where $\CTor$ is the torsion defined as in \eqref{torsion},
$\nabla$ denotes the induced connection on $E\otimes TM$
and $\Div\left(\Tr\,\nabla\rho\right)=\sum_a \mg\left(\xi^a,\conn_{\rho(\xi^a)}^E\left(\Tr\,\nabla\rho\right)\right)$
is the divergence of $\Tr\,\nabla\rho\in\Ga(E)$ (see Eq.\ \eqref{Eq:Div}).
%w.r.t.\ the $E$-connection on $E$
%induced by $\nabla^E$ (see Eq.\ \eqref{eq:nablarho}).
%Here,
%with $(\xi^a)$ is a local basis of sections of $E$.

\begin{thm}\label{cor:Inv}
The function $f_E$ 
%defined by Eq.\ \eqref{Inv fE}
 satisfies the following properties:
\begin{enumerate}[label=(\roman*)]
\item $f_E$ is independent of the choice of the connections $\nabla^{TM}$
and $\nabla^E$;
\item $\{\Theta,f_E\}=0$, where $\Theta=\rho-\CTor$
is the Hamiltonian generating function of $(E,\mg,\rho,\llbracket\cdot,\cdot\rrbracket)$;
\item if $(E,\mg)$ admits a twisted spinor bundle,
then $f_E=-D^2$, where $D$ is the skew-symmetric
Dirac generating operator of $(E,\mg,\rho,\llbracket\cdot,\cdot\rrbracket)$
as in  Theorem \ref{thm:DiracGenOp}.
\end{enumerate}
\end{thm}
\begin{proof}
We prove the last assertion first.
In view of \eqref{Inv fE},
it suffices to prove this locally
over a contractible open subset $U\subset M$.
Being restricted to $U$, a pseudo-Euclidean vector bundle $(E|_U,\mg)$
always admits a twisted spinor bundle $\sS|_U$.
According to Theorem~\ref{thm:DiracGenOp},  we obtain a skew-symmetric
Dirac generating operator $D\in\D(U,\sS|_U)_\bbR$.
By a straightforward computation using Eq.\ \eqref{Formula:D},
we check that $f_E=-D^2$.

Since the skew-symmetric
Dirac generating operator $D$ is   independent of the
 choice of the connections $\nabla$ and $\nabla^E$ according
to Theorem~\ref{thm:DiracGenOp}, so does
$f_E=-D^2$. This proves the assertion (i).

Finally, we have $\{\Theta,f_E\}=-\sigma_2([D,D^2])=0$.
Hence the assertion (ii) follows.
\end{proof}

As a consequence, the function $f_E$ is indeed
 an invariant of the Courant algebroid.
Our next goal is to provide an intrinsic formula for this function.
We first recall some standard facts on the structure of
 a Courant algebroid $(E,\mg,\rho,\llbracket\cdot,\cdot\rrbracket)$.
It is well-known that $(\ker\rho)^\perp\subset\ker\rho$ (see e.g.\ \cite{Roy02}).
Moreover, $(\ker\rho)^\perp$ and $\ker\rho$ are two-sided ideals in $E$
for the Dorfman bracket. Hence, the quotient
$$
  \GG:=\ker\rho/(\ker\rho)^\perp
$$
is a bundle of quadratic Lie algebras, whose fiberwise Lie brackets
 $[\cdot,\cdot]^\mathcal{G}$
 and non-degenerate $\ad$-invariant scalar product $(\cdot,\cdot)^\mathcal{G}$
are inherited, respectively,
 from the Dorfman bracket $\llbracket\cdot,\cdot\rrbracket$
and the metric~$\mg$ (for further details, see \cite{CSX09}).
Since fibers of $\GG$ are quadratic Lie algebras, we have an analogue
of Cartan $3$-form:
% $\CTor[\GG]\in\Ga(\wedge^3 \GG)$:
$$
  \CTor[\GG](\br_1,\br_2,\br_3)
  := \left([\br_1,\br_2]^\GG,\br_3\right)^\GG \in \wedge^3 \GG|_x ,
$$
for all $\br_1,\br_2,\br_3\in \GG|_x, \ \forall x\in M$.
 If $\rho$ is of constant rank
over a neighborhood of a point $x\in M$, we say that the Courant algebroid
is regular at the point $x$. In this case, $\CTor[\GG]$ is smooth in a neighborhood of $x$.
Note that regular points form a dense open subset $M_{\text{reg}}$
 of the base manifold $M$.

\begin{thm}
Let $(E,\mg,\rho,\llbracket\cdot,\cdot\rrbracket)$ be a Courant algebroid.
On $M_{\text{reg}}$, we have
$$
f_E= \mg(\CTor[\GG],\CTor[\GG]).
$$
\end{thm}
\begin{proof}
As both sides of the above equality depend only on the local structure of
the Courant algebroid, we can assume that $E$ is a regular Courant algebroid.

We recall some useful constructions regarding
regular Courant algebroids  \cite{CSX09}.
The vector bundle $E$ admits a splitting:
\[\label{Eregular}
  E\cong F^*\oplus \GG\oplus F,
\]
where $F=\rho(E)\subseteq TM$ is an integrable distribution.
Under the above isomorphism, the anchor map $\rho$
becomes the  projection $\rho:E\to F$ and each section $\xi\in\Ga(E)$
decomposes as $\xi=\eta+\br+x$, with $\eta\in\Ga(F^*)$,
$\br\in\Ga(\GG)$ and $x\in\Ga(F)$.
In the remainder of the proof, we use such a decomposition
$\xi_i=\eta_i+\br_i+x_i$, for all $\xi_i\in\Ga(E)$.
Let $\nabla^{TM}$ be a torsion-free  connection on $TM$ and denote
by $\nabla^F$ and $\nabla^{F^*}$ the induced connections
on $F$ and $F^*$ respectively.
According to \cite[Proposition 1.1]{CSX09}
and \cite[Lemma 1.4]{CSX09},
there exists
a metric $F$-connection $\nabla^{\GG} :\Ga(F)\otimes\Ga(\GG)\to\Ga(\GG)$
and a section $\calH\in\Ga(\wedge^3 F^*)$
such that the formula
\[\label{Eq:Enabla'}
  \Econn_{\xi_1}'\xi_2 = \big(\conn_{x_1}^{F^*}\eta_2-\frac13\calH(x_1,x_2,\cdot)\big)
+\big(\conn_{x_1}^\GG\br_2+\frac23[\br_1,\br_2]^\GG\big)+\conn_{x_1}^F x_2 , \quad \forall \xi_1,\xi_2\in\Ga(E),
%  \Econn_{\theta+\br+x}'(\eta+\s+y) = \big(\conn_x^{F^*}\eta-\frac13\calH(x,y,\cdot)\big)
%+\big(\conn_x^\GG\s+\frac23[\br,\s]^\GG\big)+\conn_x^F y ,
\]
defines an
 $E$-connection $\Enabla'$ on $E$.
As shown in \cite[Formula~(3.1)]{CSX09},  the torsion of $\Enabla'$
reads as
\[\label{CtorEnabla'}
  \CTor[\Enabla'](\xi_1,\xi_2,\xi_3)
=\calH(x_1,x_2,x_3)+\cycl_{123} \left(R^\GG(x_1,x_2),\br_3\right)^\GG -\CTor[\GG](\br_1,\br_2,\br_3),
\]
for all $\xi_i\in\Ga(E)$, $i=1,2,3$, where $R^\GG :\wedge^2 F\to\GG$ is a certain  bundle map.

To compute $f_E$ via Eq.\ \eqref{Inv fE},
we need a metric connection on $E$ and a linear connection on $TM$.
For the latter, we choose the above connection $\nabla^{TM}$.
As for the metric connection on $E$, we define it as follows.
Pick a Riemannian metric on $M$ and denote by $F^\perp\subseteq TM$
the distribution orthogonal to $F$.

Choose a metric connection $\tilde{\nabla}^\GG$ on the vector bundle $\GG$
and define a new metric connection on $\GG$ by setting
$$
\conn_X^\GG\br:=\conn_{X_F}^\GG (\br)+\tilde{\nabla}^\GG_{\! X_{F^\perp}}(\br), \qquad\forall X\in{\mathfrak X}(M),\br\in\Ga(\GG),
$$
where the splitting $X=X_F+X_{F^\perp}$ corresponds to the decomposition
$TM\cong F\oplus F^\perp$.
Then a  metric connection on $E$ can  be defined by  the following formula
\[\label{nablaE}
  \nabla^E: = \nabla^{F^*}\!\oplus\nabla^\GG\oplus\nabla^F.
\]
Since the anchor map $\rho$ is identified with the projection onto $F$
through the decomposition \eqref{Eregular},
it is simple to check that
$$
  \nabla\rho = \nabla^{TM}\circ\rho-\rho\circ\nabla^E=0.
$$
Therefore, according to Eq.\ \eqref{Inv fE},  we have
\[\label{Formula1-fE}
  f_E = \mg(\CTor,\CTor),
\]
with $\CTor$ being the torsion of $E$ with respect to $\nabla^E$.

Our next goal is to compute $\CTor$.
The linear connection $\nabla^E$ defined in  \eqref{nablaE}
induces an $E$-connection on $E$ given explicitly by
$$
  \Econn_{\xi_1}\xi_2 = \conn_{x_1}^{F^*}\eta_2+\conn_{x_1}^\GG\br_2+\conn_{x_1}^F x_2, \qquad \forall\xi_1,\xi_2\in\Ga(E).
%  \Econn_{\theta+\br+x}(\eta+\s+y) = \conn_x^{F^*}\eta+\conn_x^\GG\s+\conn_x^F y,
$$
In view of Eq.\ \eqref{Eq:Enabla'},
the $E$-connections $\Enabla$ and $\Enabla'$ are related as follows:
$$
  \Econn_{\xi_1}' \xi_2 -\Econn_{\xi_1}\xi_2  = -\frac13\calH(x_1,x_2,\cdot) +\frac23[\br_1,\br_2]^\GG.
$$
Using the fact that $\CTor[\Enabla]=\CTor$ and Eq.\ \eqref{TorsionEconn},
we deduce that
\begin{align*}
  \CTor(\xi_1,\xi_2,\xi_3) -\CTor[\Enabla'](\xi_1,\xi_2,\xi_3)
  &= \cycl_{123}\; \mg\left(-\frac13\calH(x_1,x_2,\cdot)+\frac{2}{3}[\br_1,\br_2]^\GG, \ \xi_3\right),\\
  &= -\frac13\cycl_{123} \left(\calH(x_1,x_2,x_3)\right)+\frac{2}{3}\cycl_{123}\; \mg\left([\br_1,\br_2]^\GG,\br_3\right),\\
  &= -\calH(x_1,x_2,x_3)+2\CTor[\GG](\br_1,\br_2,\br_3),
\end{align*}
for all
  $\xi_i=\eta_i+\br_i+x_i \in\Ga(E)$, $i=1,2,3$.
Therefore, by Eq.\ \eqref{CtorEnabla'}, we have
$$%\[\label{Ctor}
  \CTor(\xi_1,\xi_2,\xi_3)
  = \CTor[\GG](\br_1,\br_2,\br_3)
   +\cycl_{123} \left(R^\GG (x_1,x_2), \ \br_3\right)^\GG.\\
$$%\]
This means that
$  \CTor-\CTor[\GG]\in\Ga(\wedge^2 F^*\otimes \GG)$.
Using the fact that $\mg(F^*,F^*\oplus\GG)=0$ and Eq.\ \eqref{Formula1-fE},
we obtain $f_E=\mg(\CTor[\GG],\CTor[\GG])$. This
concludes the proof of the theorem.
\end{proof}

\begin{rem}
By definition----see   Eq. \eqref{Inv fE}, the function $f_E$ is always 
smooth, even for a non-regular Courant algebroid $E$.
Therefore, one can consider
the function $f_E$ as the smooth extension of $\mg(C_\GG,C_\GG)$
from $M_{\text{reg}}$  to $M$.
There is no reason to expect that the function $\mg(C_\GG,C_\GG)$
itself is smooth at  non-regular points of $E$. 
However we were unable
to find an example of Courant algebroid
such that $\mg(C_\GG,C_\GG)\neq f_E$.% is not smooth.
\end{rem}

%%%% 6 %%%%%%%%%%%%%%%%%%%%%%%%%%%%%%%%%%%%%%%%%%%%%%%%%%%%
%%%%%%%%%%%%%%%%%%%%%%%%%%%%%%%%%%%%%%%%%%%%%%%%%%%%%%%%
\section{Applications to Lie bialgebroids}%The splittable case
%%%%%%%%%%%%%%%%%%%%%%%%%%%%%%%%%%%%%%%%%%%%%%%%%%%%%%%%
%%%%%%%%%%%%%%%%%%%%%%%%%%%%%%%%%%%%%%%%%%%%%%%%%%%%%%%%

%%%%6.1%%%%%%%%%%%%%%%%%%%%%%%%%%%%%%%%%%%%%%%%%%%%%%%%%%%%
\subsection{Weyl quantization in the splittable case}\label{Sec:WQsplit}
%%%%%%%%%%%%%%%%%%%%%%%%%%%%%%%%%%%%%%%%%%%%%%%%%%%%%%%%

From now on, we assume  that the pseudo-Euclidean
vector bundle $(E,\mg)$ splits as $E:=A\oplus A^*$
and the pseudo-metric on $E$  is given by
\[\label{Eq:gAA}
\mg( \zeta_1+\eta_1, \zeta_2+\eta_2)=\langle \zeta_1,\eta_2 \rangle+\langle \zeta_2,\eta_1 \rangle,
\qquad \forall\zeta_1,\zeta_2\in\Ga(A), \eta_1,\eta_2\in\Ga(A^*),
\]
where $\langle\cdot,\cdot\rangle$ is the duality pairing between $A$ and $A^*$.
%Here, $\zeta_1,\zeta_2\in\Ga(A)$ and $\eta_1,\eta_2\in\Ga(A^*)$.
The previous constructions can then be described more explicitly.

According to Example \ref{ex:wedgeA}, $S_\bbR:=\wedge A^*$
is a real spinor bundle of $(E,\mg)$. Under the identification
$$
\Ga(S_\bbR)=\cO(A[1]),
$$
the Clifford action of $\Ga(E)$, given in \eqref{Eq:CliffordAction1},
reads as
\[\label{Eq:CliffordAction}
\g(\zeta)\phi=\langle\zeta,\eta^a\rangle\left(\frac{\partial}{\partial \eta^a} \phi\right) \quad\text{and}\quad \g(\eta)\phi=\eta\phi, %\qquad \forall \phi\in\cO(A[1]),
\qquad \forall \zeta\in\Ga(A),\eta\in\Ga(A^*),\phi\in\cO(A[1]),
\]
where  $(\eta^a)$ is a  local frame of~$A^*$.
We will need the following relations later on:
\[\label{Eq:Cliffordbis}
\g(\psi)\phi=\psi\phi \quad\text{and}\quad
\g(\psi\zeta)\phi=\left(\langle\zeta,\eta^a\rangle\psi\frac{\partial}{\partial \eta^a} +\frac{(-1)^\k}{2}\langle\zeta,\psi\rangle\right)\phi,
\]
for all $\zeta\in\Ga(A)$, $\psi\in\cO_\k(A[1])$, $\phi\in\cO(A[1])$.
Assume that the line bundle $\ccL:=\wedge^\Top A\otimes\wedge^\Top T^*M $
admits a square root, denoted $\ccL^{1/2}$,
 and set $\sS_\bbR:=\wedge A^*\otimes\ccL^{1/2}$.
Since $\det S_\bbR\cong(\wedge^{\mathrm{top}}A^*)^{\frac{N}{2}}$,
with $N$ being the rank of $\wedge A^*$, the $\bbC$-vector bundle
$\sS:=\sS_\bbR\otimes\bbC$ defines a twisted spinor bundle
(see Section \ref{ParReal}). Accordingly, the $\bbR$-vector bundle
$\sS_\bbR$ is called a \emph{real twisted spinor bundle}.
By pull-back along $\pi:A[1]\to M$, the line bundle $\pi^*\ccL$
can be identified with the Berezinian line bundle 
$\Ber_\aA\rightarrow A[1]$
(see e.g.\ \cite{MR2275685, Gra12, KMo02}).
 Therefore, we obtain an isomorphism of $\cO(A[1])$-modules
\[\label{S:Ber}
\upsilon:\Ga(\sS_\bbR)\xto{\ \sim\ }\Ga\big(\Ber_{\aA}^{\,1/2}\big).
\]

A linear connection on  the vector bundle
$A$ induces a metric connection $\nabla^E$ on $E$
and a compatible spinor connection $\nabla^S$ on~$S_\bbR$.
It in turn, together with  a connection on $TM$,
 induces a spinor connection $\nabla^\sS$ on $\sS_\bbR$,
compatible with $\nabla^E$,
and a   connection $\nabla$ on
the Berezinian line bundle $\Ber_{\aA}^{\,1/2}\rightarrow A[1]$.
The latter can be defined as
the pull back connection, via the canonical projection
$\pi: A[1]\to M$, of the induced connection
on the line bundle $\ccL^{1/2}\to M$. It can be expressed explicitly  
 as follows.
First, note that the space of vector fields $\mathfrak{X}(A[1])$ is linearly spanned
by the vector fields $X$ which,
as derivations of $\cO(A[1])=\Ga(S_\bbR)$,
take the form
\[\label{Eq:Decomp}
X=f\cdot\conn_{X_0}^S+g\cdot\g(\zeta),
\]
with $f,g\in\cO(A[1])$, $X_0\in{\mathfrak X}(M)$ and $\zeta\in\Ga(A)$.
By setting
\begin{equation}
\label{eq:CDG}
\conn_X:=\upsilon\circ\left(f\cdot \conn_{X_0}^\sS +g\cdot\g(\zeta) \right)\circ\upsilon^{-1},
\end{equation}
one obtains a well-defined linear connection $\nabla$ on
 $\Ber_{\aA}^{\,1/2}\to A[1]$.

The algebra  $\D\big(A[1],\Ber_{\aA}^{\,1/2}\big)$,
of differential operators on $\Ber_{\aA}^{\,1/2}$, is a subalgebra of
$\End\big(\Ga(\Ber_{\aA}^{\,1/2})\big)$ generated by multiplication operators by functions on $A[1]$
and covariant derivatives $\conn_X$, $\forall X\in\mathfrak{X}(A[1])$---see
\eqref{eq:CDG}.

\begin{lem}
The map defined by $\Upsilon_A(D):=\upsilon\circ D\circ\upsilon^{-1}$, $\forall D\in \D(M,\sS_\bbR)$,
provides an isomorphism of algebras
$$
\Upsilon_A:\D(M,\sS_\bbR)\xto{\ \sim\ }\D\big(A[1],\Ber_{\aA}^{\,1/2}\big).
$$
\end{lem}

\begin{proof}
By \ \eqref{S:Ber}, the map
$\Upsilon_A(\cdot )=\upsilon\circ \,\cdot\,\circ\upsilon^{-1}$
is an algebra isomorphism between $\End(\Ga(\sS_\bbR))$
and $\End(\Ga(\Ber_{\aA}^{\,1/2}))$.
Since generators of $\D(M,\sS_\bbR)$ and
$\D\big(A[1],\Ber_{\aA}^{\,1/2}\big)$ are clearly in bijection
via $\Upsilon_A$, the result follows.
\end{proof}

Since the metric $\mg$ in Eq.\ \eqref{Eq:gAA} is the duality pairing,
by Lemma~\ref{lem:splitting}, we have a symplectic diffeomorphism
$\widetilde{\Xi}_{\nabla}$  between $\big(T^*[2](A[1]), \om_{\text{can}}\big)$
 and
$\big(T^*[2]M\oplus E[1],\om_{\mg,\nabla^E}\big)$.
Therefore, the map
$\Phi_A :=(\widetilde{\Xi}_{\nabla})^{-1}$
defines a symplectic diffeomorphism:
\begin{equation}
\label{eq:A}
\Phi_A :T^*[2]M\oplus E[1]\xto{\ \sim\ } T^*[2](A[1]).
\end{equation}
Consider the Weyl quantization map
$\WQ:\cO(T^*[2]M\oplus E[1])\rightarrow\D(M,\sS)_\bbR$
(see  Proposition \ref{cor:WQ}).
According to Remark~\ref{rem:realspinor},
the restriction from $\sS$ to $\sS_\bbR$
induces an algebra isomorphism
$$
\mathcal{R}:\D(M,\sS)_\bbR\xto{\ \sim\ }\D(M,\sS_\bbR).
$$
%Using the symplectic diffeomorphism $\Phi$,
%We obtain a quantization of $T^*[2]A[1]$:

\begin{defn}
\label{def:WQ1}
The Weyl quantization on $T^*[2](A[1])$ is the map
$$
\WQ^A : \cO\big(T^*[2](A[1])\big)\longrightarrow\D\big(A[1],\Ber_{\aA}^{\,1/2}\big),
$$
 given by 
$$\WQ^A:=\Upsilon_A\circ\mathcal{R}\circ\WQ\circ(\Phi_A )^*.$$
\end{defn}

Note that $\WQ^A$ is a linear isomorphism.

\begin{rem}
Connections on $A$ and $TM$ give rise to a connection $\nabla$
on the graded vector bundle on
$T(A[1])\to A[1]$. 
It would be interesting to compare the quantization map $\WQ^A$
with the Weyl quantization constructed directly from the induced connections
on $T(A[1])\to A[1]$ and $\Ber_{\aA}^{\,1/2}\to A[1]$ via Eq.\ \eqref{Def:cQ}.
\end{rem}

There is a notion of Lie derivative on the space
of sections $\Ga(\Ber_{\aA}^{\,1/2})$ (see e.g.\ \cite{KMo02}).
We will need its expression in local affine coordinates $(x^i,\eta^a)$ on $A[1]$.
Assume that $X\in\mathfrak{X}(A[1])$ is a vector field of degree $\k$.
%i.e., $X$ maps $\cO_k(A[1])$ into $\cO_{k+\k}(A[1])$ for all $k\in\bbN$.
Using the local expression
$X=X^i\frac{\partial}{\partial x^i}+X^a\frac{\partial}{\partial \eta^a}$,
where $X^i$ and $X^a$ are local functions on $A[1]$,
the Lie derivative on $\Ber_{\aA}^{\,1/2}$ reads as
\[\label{Local_Lie}
\Lie_X=X^i\frac{\partial}{\partial x^i}+X^a\frac{\partial}{\partial \eta^a}
+\frac{1}{2}\left(\frac{\partial}{\partial x^i}X^i
+(-1)^{\k+1}\frac{\partial}{\partial \eta^a}X^a \right),
\]
in the trivialization of
 $\Ber_{\aA}^{\,1/2}$ provided by $( \prod\zeta_a \otimes \wedge dx^i )^{1/2}$,
with $(\zeta_a)$ being the dual frame of~$(\eta^a)$.

\begin{prop}\label{prop:WQLX}
For any vector field $X\in\mathfrak{X}(A[1])$, we have
\[\label{WQLX}
\WQ^A(F_X)=\Lie_X,
\]
where $F_X$ is the fiberwise linear function on $T^*[2](A[1])$
corresponding to $X$.
\end{prop}
\begin{proof}
Let $X\in\mathfrak{X}(A[1])$ be any vector field of degree $\k$.
It suffices to prove Eq.~\eqref{WQLX} at each point $x\in M$.
%As in the proof of Proposition \ref{Qexplicit},
Consider local affine coordinates $(x^i,\eta^a)$ on $A[1]$
obtained by pull-back of a Cartesian coordinate system on $T_xM\times A_x$,
via the map $\cT_x$ defined in ~\eqref{map:T}.
This system of coordinates induces fiberwise coordinates
$(\zeta_a)$ of $A$ and $(p_i)$ of $T^*M$.
Then, $(x^i,\eta^a,p_i,\zeta_a)$ is a local coordinate system of
both $T^*[2](A[1])$ and %$(x^i,p_i,\zeta^a,\eta_a)$ is a local coordinate system of
$T^*[2]M\oplus E[1]$.
%which are also
By Eq.\ \eqref{Eq:fxietaX}, the map
 $\Phi_A :T^*[2]M\oplus E[1]\to T^*[2](A[1])$ in
\eqref{eq:A} satisfies
$$
(\Phi_A )^*x^i=x^i,\quad(\Phi_A )^*p^\nabla_i=p_i,\quad(\Phi_A )^*\eta^a=\eta^a, \quad (\Phi_A )^*\zeta_a=\zeta_a,
$$
where $p_i^\nabla \in C^\infty (T^*[2](A[1])) $ 
is the fiberwise linear function on $T^*[2](A[1])$
corresponding to the vector field $\conn_{i}^A$ on $A[1]$. %\frac{\partial}{\partial x^i}
Note that
$\conn_{i}^A=\frac{\partial}{\partial x^i}+\Ga_{ib}^a \eta^b \frac{\partial}{\partial \eta^a}$,
%Hence, we have $p_i^\nabla=
where the real function $\Ga_{ib}^a$ vanishes at $x$
by definition of the coordinates  $(x^i,\eta^a)$.
Using the local expression
$X=X^i\frac{\partial}{\partial x^i}+X^a\frac{\partial}{\partial \eta^a}$,
we obtain
$$
(\Phi_A )^*F_X= X^i\cdot\big(p_i+\Ga_{ib}^a\eta^b\zeta_a\big)+ X^a\zeta_a.
$$
%The identification maps $\mathcal{R}$ and $\Upsilon_A$
%entering in the definition of $\WQ^A$ can be considered as trivial around $x$.
From Eqns \eqref{Eq:localQ}-\eqref{Def:gamma}
and $\Ga_{ib}^a(x)=0$, we deduce that
$$
  \left(\WQ^A(F_X)\,\psi\right)(x,\eta) =
\left(\g\Big(X^i(x,\eta)\Big)\frac{\partial}{\partial x^i} +\frac{1}{2}\,\g\Big(\frac{\partial}{\partial x^i} X^i(x,\eta)\Big)
+\g\Big( X^a(x,\eta) \zeta_a\Big)\right)\psi\,(x,\eta),
$$
in the trivialization of $\Ber_{\aA}^{\,1/2}$ provided by
 $( \prod\zeta_a \otimes \wedge dx^i )^{1/2}$.
%for all $\psi\in\cO(A[1])$.
The result thus follows
from Eqns \eqref{Eq:Cliffordbis} and \eqref{Local_Lie}.
\end{proof}

\begin{rem}
Let $X\in\mathfrak{X}(A[1])$ be  a vector field
as in Eq.\ \eqref{Eq:Decomp}. According to Eqns \eqref{WQT} and  \eqref{WQLX},
the Lie derivative on $\Ber_{\aA}^{\,1/2}$ is given by
$$
\Lie_X=\nabla_X+\frac{1}{2}\Big(\Tr \nabla (fX_0)+\langle \zeta,g\rangle\Big).
$$
\end{rem}

\begin{rem}
Let $\a':\cO^\bbC(A[1])\rightarrow\cO^\bbC(A[1])$ be the
$\iCinfty(M)$-linear antiautomorphism
satisfying $\a'(\phi_0)=\phi_0$ if $\phi_0\in\cO^\bbC_1(A[1])$.
%By definition, $\a(\phi_0\psi_0)=\a(\psi_0)\a(\phi_0)$
%for all $\phi_0,\psi_0\in\cO^\bbC(A[1])$. Via Eq.\ \eqref{S:Ber},
The map $\a'$ extends naturally to
$\Ga\big(\Ber_{\aA}^{\,1/2}\big)\otimes\bbC\cong\cO^\bbC(A[1])\otimes_\Cinfty\Ga(\ccL^{1/2})$
by setting $\a:=\a'\otimes\Id$.
Let
$$
\ve=\begin{cases}  1  \quad\text{if } \rk A \equiv0,1 \mod 4,\\
   \bi \quad\text{ if } \rk A \equiv2,3 \mod 4.
 \end{cases}
$$
Via the Berezin integration over the supermanifold $A[1]$,
we define a pseudo-Hermitian scalar product
$$
(\phi,\psi)=\ve\int_{A[1]}\overline{\a(\phi)}\psi,
$$
on the space of compactly supported sections
$\Ga_c(\Ber_{\aA}^{\,1/2})\otimes\bbC$.
Under the isomorphism
$\upsilon\otimes \Id_\bbC:\Ga(\sS)\xto{\ \sim\ }\Ga(\Ber_{\aA}^{\,1/2})\otimes\bbC$,
the above scalar product on $\Ber_{\aA}^{\,1/2}\otimes\bbC$
furnishes a globalization of the local spinor scalar product on $\sS$
given in Eq.\ \eqref{HermitianProduct}.
Accordingly, the adjoint operation on $\D(M,\sS)$,
defined in Proposition \ref{prop:*},
coincides with the one defined on
$\D\big(A[1],\Ber_{\aA}^{\,1/2}\otimes\bbC\big)$
by the above scalar product.
\end{rem}

%%%%6.2%%%%%%%%%%%%%%%%%%%%%%%%%%%%%%%%%%%%%%%%%%%%%%%%%%%%
\subsection{Preliminaries on Lie algebroids}
%%%%%%%%%%%%%%%%%%%%%%%%%%%%%%%%%%%%%%%%%%%%%%%%%%%%%%%%

Now let $ (A,[\cdot,\cdot],\rho)$ be a Lie algebroid.
According to \cite{ELW99}, the line bundle
$\ccL=\wedge^\Top A\otimes\wedge^\Top T^*M$
is an $A$-module with the $A$-action being given by
$$
D_\zeta (v\otimes s):=[\zeta, v]\otimes s+v\otimes \Lie_{\rho(\zeta)}s,
$$
for any 
 $\zeta\in\Gamma(A)$, $v\in\Gamma(\wedge^\Top A)$ and $s\in\Gamma(\wedge^\Top T^*M)$.
Here the bracket stands for the Schouten bracket on $\Ga(\wedge A)$,
and $\Lie_{\rho(\zeta)}s$ stands for
the Lie derivative of $s\in\Gamma(\wedge^\Top T^*M)$ along the vector field
$\rho(\zeta)\in\Vect(M)$.  Since the map $\zeta\mapsto D_\zeta$ is $\iCinfty(M)$-linear,
for each section $v\otimes s\in\Ga(\ccL)$,
there exists a unique $\eta_0\in\Ga(A^*)$ such that
\[\label{ModularCocycle}
D_\zeta (v\otimes s)=\<\eta_0,\zeta\>\, v\otimes s, \qquad\forall \zeta\in\Ga(A).
\]
The element $\eta_0$ is called
the modular 1-cocycle of the Lie algebroid $(A,[\cdot,\cdot],\rho)$
associated to the section $v\otimes s\in\Ga(\ccL)$.
Assume the square root of $\ccL$ exists. Then, $\ccL^{1/2}$ admits an action
of the Lie algebroid $A$, defined by
\[\label{Eq:actionBer}
\tilde{D}_\zeta \big((v\otimes s)^{1/2}\big):=\frac{1}{2}\<\eta_0,\zeta\>\,(v\otimes s)^{1/2}.
\]

The Chevalley-Eilenberg differential $Q: \Ga(\wedge^kA^*)\to \Ga(\wedge^{k+1}A^*)$ of the  Lie algebroid $A$ is given by
 %the differential of the Chevalley-Eilenberg cochain  complex $\bigoplus_k\Ga(\wedge^kA^*)$ given by
\begin{equation*}\begin{split}
  Q(\phi_0)(\zeta_0,\ldots,\zeta_k) =\,& \sum_{a=0}^{\rk A} (-1)^a \rho(\zeta_a)
\left[\phi_0(\zeta_0,\ldots,\hat{\zeta}_a,\ldots,\zeta_k)\right]\\
  &+\sum_{a<b}(-1)^{a+b}\phi_0([\zeta_a,\zeta_b],\zeta_0,\ldots,\hat{\zeta}_a,\ldots,\hat{\zeta}_b,\ldots,\zeta_k),
\end{split}\end{equation*}
where $\phi_0\in\Ga(\wedge^kA^*)$
and $\zeta_0,\ldots,\zeta_k\in\Ga(A)$.
Equivalently, $Q$ can be considered as
 a homological vector field on~$A[1]$ \cite{MR1480150}.
 %$Q\in\mathfrak{X}(A[1])$.
In local coordinates $(x^i,\eta^a)$, % on $A[1]$,
%We denote by $(\zeta_a)$ a frame on $A$ dual to $(\eta^a)$.
%under the duality pairing between $A$ and $A^*$.
 $Q$ is expressed as follows:
\[\label{Eq:Q}
Q=\rho_a^i\,\eta^a\frac{\partial}{\partial x^i}+\frac{1}{2}C^a_{bc}\,\eta^c\eta^b\frac{\partial}{\partial \eta^a},
\]
where, by definition,
$$
\rho_a^i=\<\rho(\zeta_a),dx^i\>\quad\text{and}\quad C_{bc}^a=\<[\zeta_b,\zeta_c],\eta^a\>
$$
are smooth functions on the base manifold and $(\zeta_a)$ is the dual frame to $(\eta^a)$.
From  Eq.\ \eqref{Eq:actionBer},
we obtain another differential $\widetilde{Q}$,
 defined on the complex
$\bigoplus_k\Ga(\wedge^kA^*\otimes\ccL^{1/2})\cong\Ga(\Ber_{\aA}^{\,1/2})$:
\begin{equation*}\begin{split}
  \widetilde{Q}(\phi)(\zeta_0,\ldots,\zeta_k) =\,& \sum_{a=0}^{\rk A} (-1)^a\tilde{D}_{\zeta_a}
\left(\phi(\zeta_0,\ldots,\hat{\zeta}_a,\ldots,\zeta_k)\right)\\
  &+\sum_{a<b}(-1)^{a+b}\phi([\zeta_a,\zeta_b],\zeta_0,\ldots,\hat{\zeta}_a,\ldots,\hat{\zeta}_b,\ldots,\zeta_k),
\end{split}\end{equation*}
where $\phi\in\Ga(\wedge^kA^*\otimes\ccL^{1/2})$
and $\zeta_0,\ldots,\zeta_k\in\Ga(A)$.

\begin{prop}\label{prop:Ld}
Let $Q\in\mathfrak{X}(A[1])$ be the homological vector field
of the Lie algebroid $A$.
%We have
As an operator on $\Ga(\Ber_{\aA}^{\,1/2})$, the differential $\widetilde{Q}$ satisfies
\[\label{eq:LQQ}
\widetilde{Q}=\Lie_Q.
\]
If
$\phi=\phi_0\otimes(v\otimes s)^{1/2}\in\Ga(\wedge A^*\otimes\ccL^{1/2})$,
we have
\[\label{localLieQ}
\Lie_Q(\phi)=\big(Q(\phi_0)+\tfrac12\eta_0\phi_0\big)\otimes (v\otimes s)^{1/2},
\]
where $\eta_0\in\Ga(A^*)$ is the modular 1-cocycle
of the Lie algebroid $A$ associated to the section $v\otimes s$.
\end{prop}

\begin{proof}
Let $\phi=\phi_0\otimes(v\otimes s)^{1/2}\in\Ga(\wedge A^*\otimes\ccL^{1/2})$.
A direct computation shows that
$$\widetilde{Q}\phi=\big(Q+\tfrac12\eta_0\big)\phi_0\otimes(v\otimes s)^{1/2}.
$$
It suffices  to prove that $\Lie_Q=\widetilde{Q}$ in a local coordinate system
$(x^i,\eta^a)$ on $A[1]$.
We work in the trivialization of $\ccL^{1/2}$ provided by the local section
$\left(\prod\zeta_a \otimes \wedge dx^i\right)^{1/2}$,
 with $(\zeta_a)$ being the dual frame to $(\eta^a)$.
By Eqns \eqref{Eq:CliffordAction}, \eqref{Local_Lie} and \eqref{Eq:Q}, we have
$$
\Lie_Q=Q+\frac{1}{2}\left(\frac{\partial}{\partial x^i}\rho^i_b+C^a_{ba}\right)\eta^b.
$$
By Eq.\ \eqref{ModularCocycle},
the modular 1-cocycle with respect to  the local section
$\left(\prod\zeta_a \otimes\wedge dx^i\right)$
satisfies $\eta_0=\frac{\partial}{\partial x^i}\rho^i_b+C^a_{ba}$.
The result thus follows.
\end{proof}

\begin{rem}
Each nowhere vanishing  section $v\otimes s\in\Ga(\ccL)$ defines
a nowhere vanishing section of the Berezinian bundle $1\otimes(v\otimes s)\in\Ga\big(\Ber_\aA\big)$
and then a divergence by the formula:
$\Div X=\frac{\Lie_X (1\otimes v\otimes s)}{1\otimes v\otimes s}$ (see e.g.\ \cite{KMo02}).
The above proposition shows that $\frac{1}{2}\Div Q=\eta_0$.
This result was proved in the more general situation of skew-algebroids in \cite{Gra12}.
\end{rem}

%%%%6.3%%%%%%%%%%%%%%%%%%%%%%%%%%%%%%%%%%%%%%%%%%%%%%%%%%%%
\subsection{Dirac generating operators for Lie bialgebroids}
%%%%%%%%%%%%%%%%%%%%%%%%%%%%%%%%%%%%%%%%%%%%%%%%%%%%%%%%

Assume that both
 $\big(A,[\cdot,\cdot],\rho\big)$ and $\big(A^*,[\cdot,\cdot]_*,\rho_*\big)$
are Lie algebroids, with homological vector fields
$Q\in{\mathfrak X}(A[1])$ and $Q_*\in{\mathfrak X}(A^*[1])$ respectively.
  %$Q\in\mathfrak{X}(A[1])$ and $Q_*\in\Ga\big(T(A^*[1])\big)$ respectively.
The Lie algebroid brackets, extended by the Leibniz rule,
turn $\big(\cO(A^*[1]),[\cdot,\cdot]\big)$
and $\big(\cO(A[1]),[\cdot,\cdot]_*\big)$ into
 Gerstenhaber algebras  \cite{YKS98,Xu99}.
The pair $(A,A^*)$ is a Lie bialgebroid if $Q_*$ is a derivation
of the Gerstenhaber algebra $\big(\cO(A^*[1]),[\cdot,\cdot]\big)$,
or equivalently if $Q$ is a derivation of the Gerstenhaber algebra
$\big(\cO(A[1]),[\cdot,\cdot]_*\big)$.

The duality pairing between $A$ and $A^*$ extends to
their exterior powers and leads to an isomorphism  \cite{CS09}:
$$
\beta_k : \wedge^k A\otimes (\wedge^n A^*\otimes \wedge^\Top T^*M)^{1/2}
\longrightarrow \wedge^{n-k}A^*\otimes\ccL^{1/2},%\qquad\forall k, \;0\leq k\leq n,
$$
where $n$ is the rank of $A$ and $0\leq k\leq n$. In turn, the maps $(\beta_k)$ induce an isomorphism
$$
\beta : \Ga\big(\Ber_{\aaA}^{\,1/2}\big)\xto{\ \sim\ } \Ga\big(\Ber_{\aA}^{\,1/2}\big).
$$
By symmetry in $A$ and $A^*$, the maps $\Phi_A$, $\Upsilon_A$ and $\WQ^A$
introduced in Section \ref{Sec:WQsplit} have counterparts
$\Phi_{A^*}$, $\Upsilon_{A^*}$ and $\WQ^{A^*}$, respectively.
We also denote by $\Lie_X$ the Lie derivative
 on $\Ga\big(\Ber_{\aaA}^{\,1/2}\big)$ along a vector field $X\in\mathfrak{X}(A^*[1])$.

As a consequence of Theorem~\ref{thm:DiracGenOp}, we recover
the following  theorem in \cite{CS09}.

\begin{thm}[\cite{CS09}]
\label{Thm:CS09}
Let $(A,A^*)$ be a pair of Lie algebroids with homological
vector fields $Q\in\mathfrak{X}(A[1])$ and $Q_*\in\mathfrak{X}(A^*[1])$, respectively.
Let $F_Q\in\cO\big(T^*[2](A[1])\big)$ and $F_{Q_{*}}\in\cO\big(T^*[2](A^*[1])\big)$
be their  corresponding Hamiltonian functions.
Set $\Theta=F_Q+(\Phi_{A^*}\circ\Phi_A^{-1})^*F_{Q_*}$.
The  following statements are equivalent:
\begin{enumerate}[label=(\roman*)]
\item  $(A,A^*)$ is a Lie bialgebroid;
\item  $\{\Theta,\Theta\}=0$;
\item $\left(\Lie_Q +\beta\circ \Lie_{Q_*}\circ\beta^{-1}\right)^2\in\iCinfty(M)$.
\end{enumerate}
\end{thm}

\begin{proof}
The equivalence between the first two assertions was proved in \cite{Roy99}.

We prove the equivalence between the last two assertions.
According to  Theorem~\ref{thm:DiracGenOp}, $\WQ((\Phi_A)^*\Theta)$
is a Dirac generating operator if and only if
$\{(\Phi_A)^*\Theta,(\Phi_A)^*\Theta\}=0$, as a function on $T^*[2]M\oplus E[1]$.
Since $\Upsilon_A$ is an algebra isomorphism and $\Phi_A$ is a symplectic diffeomorphism,
we deduce that
$\{\Theta,\Theta\}=0$ if and only if $\big(\WQ^A(\Theta) \big)^2\in\iCinfty(M)$,
where $\WQ^A$ is the quantization map as in Definition \ref{def:WQ1}.
The conclusion follows from the lemma below.
\end{proof}
\begin{lem}
Under the same  hypothesis of Theorem \ref{Thm:CS09}, we have
$\WQ^A(\Theta)=\Lie_Q +\beta\circ \Lie_{Q_*}\circ\beta^{-1}$.
\end{lem}
\begin{proof}
According to Proposition \ref{prop:WQLX}, we have $\WQ^A(F_Q)=\Lie_Q$.
%By symmetry in $A$ and $A^*$
Similarly,  we also have $\WQ^{A^*}(F_{Q_*})=\Lie_{Q_*}$.
Since $\Upsilon_A(\cdot)=\beta\circ\Upsilon_{A^*}(\cdot)\circ\beta^{-1}$, the 
conclusion thus follows.
\end{proof}

According to Proposition \ref{prop:Ld}, the above Dirac generating operator
is also equal to
$$
\WQ^A(\Theta)=\widetilde{Q} +\beta\circ \widetilde{Q_*}\circ\beta^{-1},
$$
where $\widetilde{Q}$ is  the differential operator
 on $\Ga(\Ber_{\aA}^{\,1/2})$
and $\widetilde{Q_*}$ is the differential  operator
 on $\Ga\big(\Ber_{\aaA}^{\,1/2}\big)$.
This expression is the one obtained in \cite{CS09}.

Choosing $s\in\Ga(\wedge^\Top T^*M)$,
$v\in\Ga(\wedge^\Top A)$ and $w\in\Ga(\wedge^\Top A^*)$,
%satisfying locally $\langle v,w\rangle=1$ for the duality pairing,
we obtain modular 1-cocycles $\eta_0\in\Ga(A^*)$ and $\zeta_0\in\Ga(A)$
of the Lie algebroids $A$ and $A^*$, respectively.
We work locally and assume that
%the section $v$ is non-vanishing and the normalization condition ... is satisfied.
$\langle v,w\rangle=1$.
Then, $v$ provides a local isomorphism
$v^\sharp :\Ga(\wedge^k A^*)\rightarrow\Ga(\wedge^{n-k}A)$,
 and  we define
$\widehat{Q}_*:=(-1)^k (v^\sharp)^{-1}\circ Q_* \circ v^\sharp$ on
$\Ga(\wedge^k A^*)$, for all $k\in\bbN$.
According to \cite{CS09}, we have
$$
\WQ^A(\Theta)=\left(Q-\widehat{Q}_*+\frac{1}{2}\eta_0 +\frac{1}{2}\g(\zeta_0)\right)\otimes\Id_{\ccL^{1/2}},
$$
under the local trivialization of $\ccL^{1/2}$ provided by the section
$(v\otimes s)^{1/2}$. As a consequence of Proposition~\ref{cor:Inv}, we obtain

\begin{cor}
The function
$2\big(\WQ^A(\Theta) \big)^2=\frac{1}{2}\langle\zeta_0,\eta_0\rangle-\widehat{Q}_*\eta_0$
is an invariant of the Lie bialgebroid $(A,A^*)$.
\end{cor}

% ==================================================================
% BIBLIOGRAPHIE
% Gestion améliorée de la bibliographie
\newcommand\eprint[1]{\href{http://arXiv.org/abs/#1}{arXiv: #1}}
\newcommand\doi[1]{\href{http://dx.doi.org/#1}{DOI: #1}}
\newcommand\eucl[1]{\href{http://projecteuclid.org/getRecord?id=euclid.#1}{Eucl: #1}}
\newcommand\proquest[2][]{\href{http://gateway.proquest.com/openurl?#2}{ProQuest: #1}}
\bibliographystyle{plain}
\bibliography{wq-Biblio}
\end{document}